\documentclass{amsart}
\usepackage{array}
\usepackage{multirow}
\usepackage{graphicx}
\usepackage{epsfig}
\usepackage{amsfonts,amsmath}
\usepackage{amssymb}
\usepackage{latexsym}
\usepackage{epsf}
\usepackage{graphicx,color,graphics}
\usepackage{amsmath,amssymb}
\usepackage{dsfont}
\usepackage[latin1]{inputenc}  
\newtheorem{theorem}{Theorem}[section]

\newtheorem{remark}{Remark}

\numberwithin{equation}{section}

\begin{document}

\title[Bresse and Timoshenko type systems with Gurtin-Pipkin's law]{Study on the stability of thermoelastic Bresse and Timoshenko type systems with Gurtin-Pipkin's law via the vertical displacements}
\author[A. Guesmia]{Aissa Guesmia}
\maketitle

\pagenumbering{arabic}

\begin{center}
Institut Elie Cartan de Lorraine, UMR 7502, Universit\'e de Lorraine\\
3 Rue Augustin Fresnel, BP 45112, 57073 Metz Cedex 03, France
\end{center}
\begin{abstract}
The objective of this paper is to study the stability of a linear one-dimensional thermoelastic Bresse system in a bounded domain, where the coupling is given through the first component of the Bresse model with the heat conduction of Gurtin-Pipkin type. Two kinds of coupling are considered; the first coupling is of order one with respect to space variable, and the second one is of order zero. We state the well-posedness and show the polynomial stability of the systems, where the decay rates depend on the smoothness of initial data. Moreover, in case of coupling of order one, we prove the equivalence between the exponential stability and some new conditions on the parameters of the system. However, when the coupling is of order zero, we prove the non-exponential stability independently of the parameters of the system. Applications to the corresponding particular Timoshenko models are also given, where we prove that both couplings lead to the exponential stability if and only if some conditions on the parameters of the systems are satisfied, and both couplings guarantee the polynomial stability independently of the parameters of the systems. The proof is based on the semigroup theory and a combination of the energy method and the frequency domain approach.
\end{abstract}

\vskip0,1truecm
asymptotic behavior, semigroup theory, energy method, frequency domain approach.
\vskip0,1truecm
{\bf AMS Classification.} 35B40, 35L45, 74H40, 93D20, 93D15.

\renewcommand{\thefootnote}{}
\footnotetext{E-mail addresse: aissa.guesmia@univ-lorraine.fr. }

\section{Introduction}

This paper is concerned with the stability of two mathematical Bresse type thermoelastic models with heat conduction given by Gurtin-Pipkin's law and working only on the vertical displacements, as well as the related Timoshenko type thermoelastic models. The first considered system is the following: 
\begin{equation}
\left\{
\begin{array}{ll}
\rho _{1}\varphi _{tt}-k_1\left( \varphi _{x}+\psi +l\,w\right)
_{x}-lk_{3}\left( w_{x}-l\varphi \right) +\delta \theta _{x}=0 & \text{in }
\left( 0,L\right) \times \left( 0,\infty \right) , \vspace{0.2cm}\\
\rho _{2}\psi _{tt}-k_2 \psi _{xx}+k_1\left( \varphi _{x}+\psi +l\,w\right) =0 &\text{in }\left( 0,L\right) \times \left( 0,\infty \right) , \vspace{0.2cm}\\
\rho _{3}w_{tt}-k_{3}\left( w_{x}-l\varphi \right) _{x}+lk_1\left( \varphi
_{x}+\psi +l\,w\right) =0 & \text{in }\left( 0,L\right) \times \left(
0,\infty \right) , \vspace{0.2cm}\\
\rho _{4}\theta _{t} -\displaystyle\int_0^{\infty}f(s)\theta _{xx} (x,t-s)ds +\delta \varphi _{xt} =0 & \text{in }\left(
0,L\right) \times \left( 0,\infty \right)
\end{array}
\right. \label{syst1}
\end{equation}
along with 
the homogeneous Dirichlet-Neumann boundary conditions
\begin{equation}
\left\{
\begin{array}{ll}
\varphi_x \left( 0,t\right) =\,\psi \left( 0,t\right) =\,w\left(
0,t\right) =\theta \left( 0,t\right) =0 & \text{in }\left( 0,\infty \right), \vspace{0.2cm}\\
\varphi_x \left( L,t\right) =\,\psi \left( L,t\right) =w\left(
L,t\right) =\theta\left( L,t\right) =0 & \text{in }\left( 0,\infty \right)
\end{array}
\right.  \label{cdt_100}
\end{equation}
and the initial data
\begin{equation}
\left\{
\begin{array}{ll}
\varphi \left( x,0\right) =\varphi _{0}\left( x\right) ,\,\varphi _{t}\left(
x,0\right) =\varphi _{1}\left( x\right) & \text{in }\left( 0,L\right) , \vspace{0.2cm}\\
\psi \left( x,0\right) =\psi _{0}\left( x\right) ,\,\psi _{t}\left(
x,0\right) =\psi _{1}\left( x\right) & \text{in }\left( 0,L\right) , \vspace{0.2cm}\\
w\left( x,0\right) =w_{0}\left( x\right) ,\,w_{t}\left( x,0\right)
=w_{1}\left( x\right) & \text{in }\left( 0,L\right) , \vspace{0.2cm}\\
\,\theta \left( x,-t\right) =\theta_{0}\left( x,t\right) & \text{in }\left( 0,L \right)\times \left( 0,\infty \right),
\end{array}
\right.  \label{cdt_10}
\end{equation}
where $\rho_{j},\, k_j,\, l$ and $L$ are positive real constants, $\delta$ is  a real constant different from zero, $f:\,\mathbb{R}_+ \to \mathbb{R}_+$ is a given function, the unknowns $\varphi,\,\psi,\, w$ and $\theta$ are functions on $(x,t)\in \left(0,L\right) \times \left( 0,\infty \right)$ and represent, respectively, the vertical displacements, the shear angle displacements, the longitudinal displacements and the temperature, 
$\varphi_j,\,\psi_j,\, w_j$ and $\theta_0$ are fixed initial data, and the subscripts $t$ and $x$ denote, respectively, the derivative with respect to the time variable $t$ and the space variable $x$. 
\vskip0,1truecm
The coupling terms $\delta \theta_x$ and $\delta \varphi_{xt}$ between the Bresse type system and the Gurtin-Pipkin's law are of order one with respect to  $x$. The second considered Bresse type thermoelastic model in this paper is the one where these coupling terms are of order zero; more precisely  
\begin{equation}
\left\{
\begin{array}{ll}
\rho _{1}\varphi _{tt}-k_1\left( \varphi _{x}+\psi +l\,w\right)
_{x}-lk_{3}\left( w_{x}-l\varphi \right) +\delta \theta =0 & \text{in }
\left( 0,L\right) \times \left( 0,\infty \right) , \vspace{0.2cm}\\
\rho _{2}\psi _{tt}-k_2 \psi _{xx}+k_1\left( \varphi _{x}+\psi +l\,w\right) =0 &\text{in }\left( 0,L\right) \times \left( 0,\infty \right) , \vspace{0.2cm}\\
\rho _{3}w_{tt}-k_{3}\left( w_{x}-l\varphi \right) _{x}+lk_1\left( \varphi
_{x}+\psi +l\,w\right) =0 & \text{in }\left( 0,L\right) \times \left(
0,\infty \right) , \vspace{0.2cm}\\
\rho _{4}\theta _{t} -\displaystyle\int_0^{\infty}f(s)\theta _{xx} (x,t-s)ds -\delta \varphi _{t} =0 & \text{in }\left(
0,L\right) \times \left( 0,\infty \right)
\end{array}
\right. \label{syst11}
\end{equation}
along with \eqref{cdt_10} and
the homogeneous Dirichlet-Neumann boundary conditions
\begin{equation}
\left\{
\begin{array}{ll}
 \varphi \left( 0,t\right) =\,\psi _{x}\left( 0,t\right) =\,w_{x}\left(
0,t\right) =\theta \left( 0,t\right) =0 & \text{in }\left( 0,\infty \right), \vspace{0.2cm}\\
\varphi \left( L,t\right) =\,\psi_x \left( L,t\right) =\,w_x\left(
L,t\right) =\theta\left( L,t\right) =0 & \text{in }\left( 0,\infty \right).
\end{array}
\right.  \label{cdt_1}
\end{equation}
\vskip0,1truecm
The Bresse type models \cite{bres} are known as a circular arch problem, while the Timoshenko type models \cite{timo} are known as a beam problem; see, for example, \cite{ag6} and the references therein. In this paper, we consider also the Timoshenko type thermoelastic models related to \eqref{syst1} and \eqref{syst11} that are corresponding to a beam with negligible longitudinal displacements; that is, 
\begin{equation}
w=l=0. \label{wl0}
\end{equation}  
\vskip0,1truecm
The stability of Bresse and Timoshenko type models have been widely studied in the literature using various controls, like frictional dampings, memories, heat conduction and boundary feedbacks. Several stability and non-stability  results depending on the considered controls and some connections between the coefficients have been established. We focus our attention on the known  results via the heat conduction, which is the subject of our paper. For more details about other kind of controls, as well as in what concerns the mathematical modeling of the thermoelasticity, we refer the readers to 
\cite{ag1, chan, ag4, gree1, gree2, 8, ag2, 9, lagn2, lagn1, liu2, 10, 12}.

\subsection{Bresse type thermoelastic models} A general Bresse type model with heat conduction can be presented in the form 
\begin{equation}
\left\{
\begin{array}{ll}
\rho _{1}\varphi _{tt}-k_1\left( \varphi _{x}+\psi +l\,w\right)
_{x}-lk_{3}\left( w_{x}-l\varphi \right) +A_1\theta=0 & \text{in }
\left( 0,L\right) \times \left( 0,\infty \right) , \vspace{0.2cm}\\
\rho _{2}\psi _{tt}-k_2 \psi _{xx}+k_1\left( \varphi _{x}+\psi +l\,w\right) ++A_2\theta =0 &\text{in }\left( 0,L\right) \times \left( 0,\infty \right) , \vspace{0.2cm}\\
\rho _{3}w_{tt}-k_{3}\left( w_{x}-l\varphi \right) _{x}+lk_1\left( \varphi
_{x}+\psi +l\,w\right) +A_3\theta=0 & \text{in }\left( 0,L\right) \times \left(
0,\infty \right),
\end{array}
\right.  \label{brssegen}
\end{equation}
where the temperature variation $\theta$ from an equilibruim reference satisfies either   
\begin{equation}
\rho _{4}\theta _{t} +q_x + B_1 \varphi +B_2 \psi +B_3 w=0 \quad \text{in }\left(0,L\right) \times \left( 0,\infty \right)\label{heatI}
\end{equation}
(known as the thermoelasticity of type I) or  
\begin{equation}
\rho _{4}\theta _{tt} -k_4 \theta _{xxt} +q_x + B_1 \varphi +B_2 \psi +B_3 w=0 \quad \text{in }\left(0,L\right) \times \left( 0,\infty \right)\label{heatIII}
\end{equation}
(known as the thermoelasticity of type III), $A_j$ and $B_j$ are given operators and $q$ is the heat flux. In order to make \eqref{brssegen}-\eqref{heatIII} determined, different additional connections between $\theta$ and $q$ were considered in the literature, where in the classical theory of thermoelasticity, $q$ is expressed in therm of $\theta$ through the Fourier's law 
\begin{equation}
q=-k_5\theta_x\label{fourier}
\end{equation} 
or the Cattaneo's law (known also as the thermoelasticity of second sound) 
\begin{equation}
\rho_5 q_t +q=-k_5\theta_x\label{cattaneo}
\end{equation}
or the Gurtin-Pipkin's law 
\begin{equation}
q=-\displaystyle\int_0^{\infty}f(s)\theta_x (x,t-s)ds,\label{gurtinpipkin}
\end{equation}
the positive real numbers $k_5$ and $\rho_5$ represent, respectively, the coefficient of the thermal conductivity and the time lag in the response of the heat flux to the temperature gradient, and the kernel $f$ describes the memory effect, so \eqref{gurtinpipkin} allows to take in consideration the history of the temperature gradient.    
\vskip0,1truecm 
The authors of \cite{fato1} considered the case \eqref{brssegen}-\eqref{heatI} with \eqref{fourier} and
\begin{equation}
A_1 =A_3 =B_1 =B_3 =0,\quad A_2\theta =\delta \theta_x \quad\hbox{and}\quad
B_2\psi =\delta \psi_{xt} \label{ajbj2}
\end{equation}
and proved the exponential and polynomially stability depending on some relationships between the coefficients. The results of \cite{fato1} were extended in \cite{nawe} to the local dissipation case; that is $\delta$ is a function on $x$ and vanishes on some part of $(0,L)$.
\vskip0,1truecm
The case \eqref{brssegen}-\eqref{heatI} with \eqref{cattaneo} and \eqref{ajbj2} was treated in \cite{kama}, where similar exponential and polynomially stability results to the ones of \cite{fato1} were proved under some restrictions on the coefficients. When the heat conduction is effective on the longitudinal displacements; that is, 
\begin{equation}
A_1 =A_2 =B_1 =B_2 =0,\quad A_3\theta =\delta \theta_x \quad\hbox{and}\quad
B_3 w =\delta w_{xt}, \label{ajbj3}
\end{equation}
the exponential and polynomially stability of \eqref{brssegen}-\eqref{heatI} with \eqref{cattaneo} and \eqref{ajbj3} were proved in \cite{afas1}.
\vskip0,1truecm
The subject of \cite{ag3} was the study of the exponential and polynomially stability of \eqref{brssegen} with \eqref{fourier} in both cases \eqref{heatI} and \eqref{heatIII}, where the thermoelastic effect is effecive on the longitudinal displacements:
\begin{equation}
A_1 =A_2 =B_1 =B_2 =0,\quad B_3 w =\delta w_{xt} \quad\hbox{and}\quad
A_3\theta =\left\{
\begin{array}{ll}
\delta \theta_x & \text{in case \eqref{heatI}}, \vspace{0.2cm}\\
\delta \theta_{xt} & \text{in case \eqref{heatIII}}.
\end{array}
\right.\label{ajbj30}
\end{equation}
As in \cite{afas1, fato1, kama}, the type of stability property proved in \cite{ag3} is related to the values of the coefficients. 
\vskip0,1truecm 
The author of the present paper treated in \cite{ag6} the stability of \eqref{brssegen} with \eqref{fourier} in both cases \eqref{heatI} and \eqref{heatIII}, where the thermoelastic effect is effecive on the vertical displacements:
\begin{equation}
A_2 =A_3 =B_2 =B_3 =0,\quad B_1 \varphi =\delta \varphi_{xt} \quad\hbox{and}\quad A_1\theta =\left\{
\begin{array}{ll}
\delta \theta_x & \text{in case \eqref{heatI}}, \vspace{0.2cm}\\
\delta \theta_{xt} & \text{in case \eqref{heatIII}}.
\end{array}
\right. \label{ajbj1}
\end{equation}
The author of \cite{ag6} proved that the case \eqref{ajbj1} is deeply different from the ones \eqref{ajbj2}-\eqref{ajbj30} in the sense that, independently of the values of the coefficients, the exponential stablility does not hold, but the polynomial stability is satisfied with a decay rate depending on the smoothness of the initial data.
\vskip0,1truecm  
The first objective of this paper is to complete from the mathematical view point the work \cite{ag6} by considering the Gurtin-Pipkin's law \eqref{gurtinpipkin}. We will prove that, independently of the coefficients, \eqref{cdt_10}-\eqref{cdt_1} is not exponential stable, but it is at least polynomially stable. However, for \eqref{syst1}-\eqref{cdt_10}, we prove the exponential and non-exponential stability depending on some connections between the coefficients, as well as the polynomial stability in general.
\vskip0,1truecm
The proof of the well-posedness is based on the semigroup theory. However, the stability results are proved using the energy method combining with the frequency domain approach.  

\subsection{Timoshenko type thermoelastic models} Concerning the case of Timoshenko type models related to \eqref{brssegen}-\eqref{heatIII}; that is,
\begin{equation}
w=l =A_3 =B_3 =0, \label{wla3b3}
\end{equation} 
the stability question has attracted the attention of many researchers in the last three decades. 
\vskip0,1truecm
The authors of \cite{15} studied the case \eqref{heatI} with \eqref{fourier}
and
\begin{equation}
A_1 =B_1 =0,\quad A_2\theta =\delta \theta_x \quad\hbox{and}\quad
B_2 \psi =\delta \psi_{xt} \label{ajbj2+}
\end{equation}
and proved that the exponential stability is equivalent to
\begin{equation}
\frac{k_1}{\rho_1}=\frac{k_2}{\rho_2} . \label{equalspeeds}
\end{equation}
\vskip0,1truecm
Under the condition \eqref{equalspeeds}, the exponential stability was proved in \cite{14} for the case \eqref{heatIII}-\eqref{fourier} with
\begin{equation}
A_1 =B_1 =0,\quad A_2\theta =\delta \theta_{xt} \quad\hbox{and}\quad
B_2 \psi =\delta \psi_{xt} . \label{ajbj20+}
\end{equation} 
When \eqref{equalspeeds} is not satisfied, the authors of
\cite{13} proved that the polynomial stability holds. The results of \cite{13, 14} were extended in \cite{19} to the case \eqref{heatIII} with \eqref{fourier},
\begin{equation}
A_1 \theta =\delta \theta_{xt} ,\quad A_2\theta =-\delta \theta_t ,\quad B_1 \varphi =\delta \varphi_{xt} \quad\hbox{and}\quad B_2 \psi =\delta \psi_{t}. \label{ajbj12}
\end{equation}   
\vskip0,1truecm 
The case \eqref{heatI} with \eqref{cattaneo} and \eqref{ajbj2+} was considered in \cite{4}, where the authors proved that the system is not exponentially stable even if \eqref{equalspeeds} is satisfied. The result of \cite{4} was completed in \cite{19} by proving that the exponential stability holds if 
\begin{equation}
\left(\rho_5 -\frac{\rho_1}{\rho_4 k_1}\right)
\left(\rho_2 -\frac{\rho_1 k_2}{k_1}\right)-\frac{\rho_5 \rho_1 \delta^2}{\rho_4 k_1} =0. \label{chi0}
\end{equation}  
\vskip0,1truecm 
Concerning the analysis of the stability in case 
\eqref{heatI} with \eqref{ajbj2+} and the Gurtin-Pipkin's law \eqref{gurtinpipkin}, we mention the work \cite{2}, where it was proved that the exponential stability is equivalent to
\begin{equation}
\left(\frac{\rho_1}{\rho_4 k_1} -\frac{1}{f(0)}\right)
\left(\frac{\rho_1}{k_1} -\frac{\rho_2}{k_2}\right)-\frac{\rho_1 \delta^2}{f(0)\rho_4 k_1 k_2} =0. \label{chif}
\end{equation}   
\vskip0,1truecm
For other types of Timoshenko thermoelastic models, as well as the wave equation with Gurtin-Pipkin's law \eqref{gurtinpipkin}, we refer the readers  to \cite{1, 3, 11, 16, 17}.
\vskip0,1truecm
The second objective of this paper is to complete the work \cite{2} by considering the Gurtin-Pipkin's law \eqref{gurtinpipkin} on the vertical displacements; that is,
\begin{equation}
A_1 \theta =\delta \theta_{x} ,\quad B_1 \varphi =\delta \varphi_{xt} \quad\hbox{and}\quad A_2 =B_2 =0, \label{ajbj1+}
\end{equation}
and extend its result to the case
\begin{equation}
A_1 \theta =\delta \theta ,\quad B_1 \varphi =-\delta \varphi_{t} \quad\hbox{and}\quad A_2 =B_2 =0. \label{ajbj1++}
\end{equation} 
Using the same arguments of proofs as in the case of Bresse models, we give, for both couplings, necessary and sufficient conditions for the exponential 
stability of the corresponding Timoshenko type systems \eqref{wl0}, and prove the same polynomial stability result independently of the parameters of the systems.
\vskip0,1truecm
The paper is organized as follows: in section 2, we prove the well-posedness of \eqref{syst1}-\eqref{cdt_10} and \eqref{cdt_10}-\eqref{cdt_1}. In section 3, we show the non-exponential stability of \eqref{syst1}-\eqref{cdt_10} and \eqref{cdt_10}-\eqref{cdt_1}. In section 4, we prove the exponential stability of \eqref{syst1}-\eqref{cdt_10}. Section 5 will be devoted to the proof of the polynomial stability of \eqref{syst1}-\eqref{cdt_10} and 
\eqref{cdt_10}-\eqref{cdt_1}. Finally, we end our paper by proving the non-exponential, exponential and polynomial stability of the Timosheko type thermoelastic systems in section 6.

\section{Formulation of the Bresse models}

In this section, and under appropriate assumptions on $f$ and $l$, we give a bref idea on the proof of the well-posedness of 
\eqref{syst1}-\eqref{cdt_10} and \eqref{cdt_10}-\eqref{cdt_1}
based on the semigroup theory. 
\vskip0,1truecm
First, in oder to simplify the computations, we do not indicate the variables $x,\,t$ and $s$ except when it is necessary to avoid ambiguity, and without lose of generality, we take
\begin{equation}
\rho_1 =\rho_2 =\rho_3 =\rho_4 =L=1. \label{rhojL} 
\end{equation} 
\vskip0,1truecm
Second, we assume that 
\begin{equation}
l\ne m\pi,\quad\forall m\in \mathbb{N} \label{l} 
\end{equation}
and $f$ satisfies
\begin{equation}
f\in C^3 (\mathbb{R}_+),\quad f^{\prime}\leq 0,\quad f(0) >0,\quad lim_{s\to\infty} f(s)=0 \quad\hbox{and}\quad \displaystyle\int_0^{\infty} \vert f^{\prime \prime \prime}(s)\vert ds <\infty , \label{f} 
\end{equation}
and there exist positive real numbers $\mu_1$ and $\mu_2$ such that 
\begin{equation}
-\mu_1 f^{\prime}\leq f^{\prime \prime}\leq -\mu_2 f^{\prime}. \label{f+} 
\end{equation} 
As a simple class of $f$ satisfying \eqref{f}-\eqref{f+}, one can consider $f(s)=d_0 e^{-d_1 s}$ with $d_0 ,\,d_1\in (0,\infty)$.
\vskip0,1truecm
Third, we put $g=-f^{\prime}$, so $g:\,\mathbb{R}_+\to \mathbb{R}_+$, 
$g\in C^2 (\mathbb{R}_+)$, 
\begin{equation}
-\mu_2 g\leq g^{\prime}\leq -\mu_1 g, \label{g} 
\end{equation}
\begin{equation}
g_0 :=\displaystyle\int_0^{\infty} g(s)ds=f(0)>0\quad\hbox{and}\quad \displaystyle\int_0^{\infty} \vert g^{\prime \prime}(s)\vert ds<\infty. \label{g0} 
\end{equation}
We observe that \eqref{g}-\eqref{g0} imply that $g^{\prime}\leq 0$, $g(0)>0$ and, by integrating,
\begin{equation}
\left\{
\begin{array}{ll}
g(0)e^{-\mu_2 s}\leq g(s)\leq g(0)e^{-\mu_1 s}, \quad \forall s\in \mathbb{R}_+ ,\vspace{0.2cm}\\
g_m :=\displaystyle\int_0^{\infty} s^m g(s)ds\in (0,\infty),\quad\forall m\in \mathbb{N} .
\end{array}
\right. 
\label{gexp} 
\end{equation}
\vskip0,1truecm
Fourth, using the idea of \cite{dafe}, we consider the varibale $\eta$ and its initial data $\eta_0$ given by
\begin{equation}\label{etaname}
\eta (x,t,s)= \displaystyle\int_{t-s}^t \theta(x,\tau)\,d\tau\quad\hbox{and}\quad\eta_0 (x, s)=\displaystyle\int_{0}^s \theta_0 (x,\tau)\,d\tau.
\end{equation}
Direct computations and the use of \eqref{cdt_10} and \eqref{cdt_1} show that the functional 
$\eta$ satisfies
\begin{equation}
\left\{
\begin{array}{ll}
\eta_t (x,t,s)+\eta_s (x,t,s) =\theta (x,t) & \text{in }
\left( 0,1\right) \times \left( 0,\infty \right)\times\left( 0,\infty \right) , \vspace{0.2cm}\\
\eta (0,t,s) =\eta (1,t,s) =\eta (x,t,0) =0 & \text{in } \left( 0,1\right) \times \left( 0,\infty \right)\times\left( 0,\infty \right),
\end{array}
\right. \label{eta}
\end{equation}
where the subscript $s$ denotes the derivative with respect to $s$. According to \eqref{f} and \eqref{eta}$_2$, and integrating with respect to 
$s$, the integral in \eqref{syst1} and \eqref{syst11} can be expressed in term of $\eta$ as follows:  
\begin{equation*}
\displaystyle\int_0^{\infty}\,g (s) \eta_{xx} (x,t,s)\,ds=-\displaystyle\int_0^{\infty}\,f^{\prime} (s) \eta_{xx} (x,t,s) \,ds
=\displaystyle\int_0^{\infty}\,f (s) \eta_{xxs} (x,t,s) \,ds.
\end{equation*}
On the other hand, from the definition of $\eta$, we see that 
$\eta_s (x,t,s) =\theta(x,t-s)$. Consequently
\begin{equation}
\displaystyle\int_0^{\infty}\,g (s) \eta_{xx} (x,t,s)\,ds=\displaystyle\int_0^{\infty}\,f (s) \theta_{xx} (x,t-s) \,ds. \label{etatheta}
\end{equation}
\vskip0,1truecm
Fifth, we put 
\begin{equation}
k=\left\{
\begin{array}{ll}
1\quad &\hbox{in case \eqref{syst1}-\eqref{cdt_10}}, \vspace{0.2cm}\\
0\quad &\hbox{in case \eqref{cdt_10}-\eqref{cdt_1}},
\end{array}
\right. \label{k}
\end{equation}
we consider the Hilbert spaces
\begin{equation*}
L_0 =L^{2}\left( 0,1\right),\quad L_1 =\left\{ v\in L^{2}\left( 0,1\right): 
\displaystyle\int_0^{1} v(x)dx=0\right\},\quad
H_0 = H_0^{1}\left( 0,1\right),\quad H_1 = H^{1}\left( 0,1\right)\cap L_1  
\end{equation*}
and
\begin{equation*}
L_g =\left\{ v:\,\mathbb{R}_+\to H_0^{1}\left( 0,1\right): 
\Vert v\Vert_{L_g}^2 := \displaystyle\int_0^{\infty} g\Vert v_x\Vert^2 ds<\infty\right\},
\end{equation*}
and we introduce the energy space
\begin{equation*}
\mathcal{H} = H_k\times L_k \times H_{1-k}\times L_{1-k}
\times H_{1-k}\times L_{1-k}\times L^{2}\left( 0,1\right)\times L_g , 
\end{equation*}
where $L^{2}\left( 0,1\right)$, $H^{1}\left( 0,1\right)$ and  
$H_0^{1}\left( 0,1\right)$ are the classical Sobolev spaces, and $L^{2}\left( 0,1\right)$ is equipped with its standard inner product $\left\langle \cdot,\cdot\right\rangle$ and generated norm $\Vert\cdot\Vert$. For 
\begin{equation*}
\Phi_j =(\varphi_j ,\,\tilde{\varphi}_j,\,\psi_j ,\,\tilde{\psi}_j, \,w_j ,\,\tilde{w}_j ,\,\theta_j ,\,\eta_j )^T ,
\end{equation*}
we introduce on $\mathcal{H}$ the inner product  
\begin{equation*}
\left\langle \Phi_1,\Phi_2\right\rangle _{\mathcal{H}} = k_1\left\langle \varphi _{1,x}+\psi_{1}+l\,w_{1} ,
\varphi _{2,x}+\psi _{2}+l\,w_{2} \right\rangle +k_2\left\langle \psi _{1,x},\psi _{2,x}\right\rangle +k_{3}\left\langle w_{1,x}-l\varphi _{1} ,w_{2,x}-l\varphi _{2} \right\rangle 
\end{equation*}
\begin{equation}
+\left\langle \tilde{\varphi}_{1},\tilde{\varphi}_{2}\right\rangle 
+\langle \tilde{\psi}_{1},\tilde{\psi}_{2}\rangle +\left\langle \tilde{w}_{1},\tilde{w}_{2}\right\rangle + \left\langle \theta _{1},\theta _{2}\right\rangle + \left\langle \eta _{1},\eta _{2}\right\rangle_{L_g}. \label{innerproductH} 
\end{equation}
We see that, if 
\begin{equation*}
(\varphi ,\psi, w)\in \mathcal{H}_k :=H_k\times H_{1-k}\times H_{1-k}
\end{equation*} 
satisfying 
\begin{equation}
\mathcal{N} (\varphi ,\psi, w):=k_1\left\Vert \varphi_{x}+\psi +l\,w \right\Vert^2 +k_2\left\Vert \psi _{x} \right\Vert^2 +k_{3}\left\Vert w_{x}-l\varphi \right\Vert^2 =0,\label{Nvarphipsiw}
\end{equation}
then
\begin{equation*}
\varphi_{x}+\psi +l\,w =\psi _{x} =w_{x}-l\varphi =0.
\end{equation*}
Therefore, from the definition of $H_{k}$, we obtain
\begin{equation}
\psi =0,\quad w=\frac{-1}{l}\varphi_{x}\quad\hbox{and}\quad \varphi_{xx}+l^2\varphi =0. \label{varphipsiw0}
\end{equation}
By integrating the last equation in \eqref{varphipsiw0}, we get, for some constants $c_1$ and $c_2$,
\begin{equation}
\varphi (x)=c_1\cos \,(lx)+c_2\sin\,(lx).\label{varphi} 
\end{equation}
For case $k=0$, and from the definition of $H_0$, we have 
$\varphi(0)=\varphi(1)=0$, so we get $c_1 =0$ and 
$c_2\sin\,l=0$. Using \eqref{l}, we find that $c_2 =0$, and then \eqref{varphipsiw0}-\eqref{varphi} lead to $(\varphi ,\psi ,w)=(0,0,0)$ 
in $\mathcal{H}_k$. 
\vskip0,1truecm
For case $k=1$, and from \eqref{varphipsiw0}-\eqref{varphi}, we see that 
\begin{equation}
w (x)=c_1\sin \,(lx)-c_2\cos\,(lx),\label{varphiw} 
\end{equation}
then, similarly, because $w(0)=w(1)=0$ and sing \eqref{l}, we obtain  
$c_1 =c_2 =0$, and then \eqref{varphipsiw0}-\eqref{varphiw} imply that 
$(\varphi ,\psi ,w)=(0,0,0)$ in $\mathcal{H}_k$.
\vskip0,1truecm
We conclude that $\mathcal{H}_k$, endowed with the inner product 
$\left\langle ,\right\rangle_{\mathcal{H}_k}$ that generates $\mathcal{N}$, is a Hilbert space. Hence, $\mathcal{H}$ is also a Hilbert space.
\vskip0,1truecm
The definition of $L_1$ allows to apply Poincar\'e's inequality and it is justified by the fact that, using a change of variable, the property 
\begin{equation*}
\int_0^1 v(x)dx =0
\end{equation*}
can be obtained, where 
$v=\varphi$ in case \eqref{syst1}-\eqref{cdt_10}, and $v\in\{\psi ,w\}$ in case \eqref{cdt_10}-\eqref{cdt_1}. Indeed, for \eqref{syst1}-\eqref{cdt_10} (with \eqref{rhojL}), we put 
\begin{equation}
{\hat\varphi}(t)=\int_0^1 \varphi (x,t)dx,\quad {\hat\varphi}_0 =\int_0^1 \varphi_0 (x)dx\quad\hbox{and}\quad {\hat\varphi}_1 =\int_0^1 \varphi_1 (x)dx, \label{hatvarphi0} 
\end{equation}
integrate \eqref{syst1}$_1$ over $(0,1)$ and use \eqref{cdt_100}, we get  
\begin{equation*}
{\hat\varphi}_{tt} +l^2 k_3 {\hat\varphi}=0,
\end{equation*}
then, by integrating with respect to $t$, we find, for some constants $c_1$ and $c_2$,
\begin{equation*}
{\hat\varphi} (t)=c_1\cos\,\left(l{\sqrt{k_3}}t\right)+c_2\sin\,\left(l{\sqrt{k_3}}t\right),
\end{equation*}
therefore, using \eqref{cdt_10}, we see that
\begin{equation}
{\hat\varphi} (t)={\hat\varphi}_0\cos\,\left(l{\sqrt{k_3}}t\right)+\frac{{\hat\varphi}_1}{l{\sqrt{k_3}}}\sin\,\left(l{\sqrt{k_3}}t\right). \label{hatvarphi} 
\end{equation}  
Puting ${\bar\varphi}=\varphi-{\hat\varphi}$, we observe that \eqref{hatvarphi0} implies that 
\begin{equation*}
\int_0^1 {\bar\varphi} (x,t)dx=0,
\end{equation*}
and moreover, using \eqref{hatvarphi}, system
\eqref{syst1}-\eqref{cdt_10} is still satisfied with
${\bar\varphi}, \varphi_0 -{\hat\varphi}_0$ and 
$\varphi_1 -{\hat\varphi}_1$ instead of 
$\varphi ,\varphi_0$ and $\varphi_1$, respectively. Similarily, for \eqref{cdt_10}-\eqref{cdt_1}, integrating \eqref{syst11}$_2$ and \eqref{syst11}$_3$ over $(0,1)$, using 
\eqref{cdt_10}, \eqref{cdt_1} and the same arguments as before, similar change of variables ${\bar\psi}$ and ${\bar w}$ can be done for $\psi$ and 
$w$, respectively. 
\vskip0,1truecm
Sixth, we put
\begin{equation*}
\left\{
\begin{array}{ll}
\tilde{\varphi}=\varphi_{t},\quad \tilde{\psi}=\psi_{t},\quad \tilde{w}=w_{t} ,\vspace{0.2cm}\\
\Phi = \left( \varphi ,\,\tilde{\varphi},\,\psi ,\,\tilde{\psi},\,w,\,\tilde{w},\,\theta,\,\eta\right)^{T} ,\vspace{0.2cm}\\
\Phi_0 = \left( \varphi_0 ,\,\varphi_1,\,\psi_0 ,\,\psi_1,\,w_0,\,w_1,\,\theta_0,\,\eta_0\right)^{T} .
\end{array}
\right.  
\end{equation*}
So, exploiting \eqref{eta}$_1$ and \eqref{etatheta}, systems \eqref{syst1}-\eqref{cdt_10} and \eqref{cdt_10}-\eqref{cdt_1} can be formulated in the form
\begin{equation}
\left\{
\begin{array}{ll}
\Phi_{t}=\mathcal{A}\Phi\quad\quad\quad\hbox{in}\,\, \left( 0,\infty \right) ,\vspace{0.2cm}\\
\Phi \left( 0\right) =\Phi_{0} ,
\end{array}
\right.  \label{syst_2}
\end{equation}
where $\mathcal{A}$ is a linear operator defined by
\begin{equation}
\mathcal{A}\Phi =\left(
\begin{array}{c}
\tilde{\varphi} \vspace{0.2cm}\\
k_1\left( \varphi _{x}+\psi +l\,w\right) _{x}+lk_{3}\left( w_{x}-l\varphi \right) -\delta \frac{\partial^k}{\partial x^k}\theta
\vspace{0.2cm}\\
\tilde{\psi} \vspace{0.2cm}\\
k_2\psi _{xx}-k_1\left( \varphi _{x}+\psi +l\,w\right) \vspace{0.2cm}\\
\tilde{w} \vspace{0.2cm}\\
k_{3}\left( w_{x}-l\varphi \right) _{x}-lk_1\left( \varphi _{x}+\psi +l\,w\right) \vspace{0.2cm}\\
\displaystyle\int_0^{\infty}g\eta_{xx} ds +(-1)^k\delta \frac{\partial^k}{\partial x^k}\tilde{\varphi} \vspace{0.2cm}\\
\theta-\eta_s
\end{array}
\right) \label{A}
\end{equation} 
with domain given by
\begin{equation*}
D\left( \mathcal{A}\right) =\left\{
\begin{array}{c}
\Phi \in \mathcal{H}\mid \,\varphi ,\,\psi ,\,w\in \,H^{2}\left( 0,1\right),\,\tilde{\varphi}\in H_k,\,\theta\in H_{0}^{1}\left( 0,1\right),\\
\tilde{\psi},\,\tilde{w}\in H_{1-k},\,\eta_s\in L_g ,\,\displaystyle\int_0^{\infty}g\eta_{xx}\in L^2 \left( 0,1\right),  \\
\frac{\partial^k}{\partial x^k}\varphi (0)=\frac{\partial^k}{\partial x^k}
\varphi (1)=\frac{\partial^{1-k}}{\partial x^{1-k}}\psi\left( 0\right) =\frac{\partial^{1-k}}{\partial x^{1-k}}\psi\left( 1\right) =\frac{\partial^{1-k}}{\partial x^{1-k}}w\left( 0\right) =\frac{\partial^{1-k}}{\partial x^{1-k}}w\left( 1\right)=\eta (0)=0
\end{array}
\right\}.
\end{equation*}
Last, we prove that the operator $\mathcal{A}$ generates a linear 
$C_0$-semigroup of contractions on $\mathcal{H}$. Using \eqref{innerproductH} and \eqref{A}, integrating with respect to $x$ and using the boundary conditions \eqref{cdt_100} and \eqref{cdt_1}, we get
\begin{equation*}
\left\langle \mathcal{A}\Phi ,\Phi \right\rangle _{\mathcal{H}}=-
\displaystyle\int_0^{\infty}g\left\langle\eta_{xs},\eta_{x}\right\rangle ds=
-\frac{1}{2}\displaystyle\int_0^{\infty}g\left(\Vert\eta_{x}\Vert^2\right)_s ds,
\end{equation*}
therefore, because $\lim_{s\to\infty} g(s)=\eta_{x} (0)=0$ (thanks to \eqref{gexp}$_1$ and \eqref{eta}$_2$), we arrive at
\begin{equation}
\left\langle \mathcal{A}\Phi ,\Phi \right\rangle_{\mathcal{H}}=\frac{1}{2} 
\displaystyle\int_0^{\infty}g^{\prime}\Vert\eta_{x}\Vert^2 ds\leq 0, \label{dissp}
\end{equation}
since $g^{\prime}\leq 0$, hence $\mathcal{A}$ is dissipative in 
$\mathcal{H}$. On the other hand, we show that 
$0\in\rho\left( \mathcal{A}\right)$ ($\rho\left(\mathcal{A}\right)$ is the resolvent of $\mathcal{A}$); that is, for any 
\begin{equation*}
F :=(f_1 ,\cdots,f_8 )^T\in \mathcal{H}, 
\end{equation*}
there exists $\Phi \in D\left( \mathcal{A}\right)$ satisfying
\begin{equation}
\mathcal{A} \Phi =F.  \label{ZF}
\end{equation}
We start by noting that \eqref{A} implies that \eqref{ZF}$_1$, \eqref{ZF}$_3$ and \eqref{ZF}$_5$ are equivalent to
\begin{equation}
\tilde{\varphi} =f_1 ,\quad\tilde{\psi} =f_3 \quad\hbox{and}\quad\tilde{w} =f_5 ,  \label{z1f1}
\end{equation}
and then
\begin{equation}
\tilde{\varphi} \in H_k \quad\hbox{and}\quad \tilde{\psi} ,\,\tilde{w} \in H_{1-k}. \label{z2f2}
\end{equation}
After, we remark that  
\begin{equation}
\eta (s)=s\theta -\displaystyle\int_0^{s} f_8 (\tau)d\tau \label{etathetaf8}
\end{equation}
is the unique unknown satisfying $\eta (0)=0$ and $\eta_s =\theta-f_8$, so \eqref{ZF}$_8$ holds. 
\vskip0,1truecm 
Next, from \eqref{k}, we have 
\begin{equation}
f_9 := -f_7 +(-1)^k \delta\frac{\partial^k}{\partial x^k}f_1 \in L^2 (0,1),
\label{f1f8f9} 
\end{equation}
then the equation $-z_{xx}=f_9$ has a unique solution 
$z\in H^2 (0,1)\cap H_0^1 (0,1)$. We take
\begin{equation}
\theta =\frac{1}{g_1}z +\frac{1}{g_1}\displaystyle\int_0^{\infty} g(s) 
\displaystyle\int_0^{s}f_8 (\tau)d\tau ds, \label{thetaf8}
\end{equation}
where $g_1$ is defined in \eqref{gexp}$_2$. Hence \eqref{z1f1}, \eqref{etathetaf8}, \eqref{f1f8f9}
and \eqref{thetaf8} imply that
\begin{equation*}
\displaystyle\int_0^{\infty} g\eta ds=z\quad\hbox{and}\quad \displaystyle\int_0^{\infty} g\eta_{xx} ds=z_{xx} =-f_9 =f_7 -(-1)^k \delta\frac{\partial^k}{\partial x^k}{\tilde\varphi}\in L^2 (0,1),
\end{equation*} 
and so \eqref{ZF}$_7$ is satisfied.
\vskip0,1truecm 
Now, for $\mu_0\in (0,\mu_1 )$, we have, by applying H\"older's inequality and Fubini theorem,
\begin{equation}\label{fubini}
\begin{array}{lll}
\displaystyle\int_{0}^{\infty}g (s)\left\Vert \displaystyle\int_0^{s} f_{8x} (\tau)d\tau\right\Vert^{2} \,ds
& = & \displaystyle\int_{0}^{\infty}g (s)\left\Vert \displaystyle\int_0^{s} e^{\frac{\mu_0}{2}\tau}e^{-\frac{\mu_0}{2}\tau}f_{8x} (\tau)d\tau\right\Vert^{2} \,ds \\
& \leq & \displaystyle\int_{0}^{\infty}g (s)\left(\displaystyle\int_{0}^{s}e^{\mu_0 y}dy\right)\displaystyle\int_{0}^{s}e^{-\mu_0\tau}\Vert  f_{8x}(\tau)\Vert^{2}\, d\tau\,ds \\
\\
& \leq & \frac{1}{\mu_0}\displaystyle\int_{0}^{\infty}\Vert f_{8x}(\tau)\Vert^{2}e^{-\mu_0\tau} \displaystyle\int_{\tau}^{\infty}g (s) \left(e^{\mu_0 s}-1\right)\,ds\, d\tau.
\end{array}
\end{equation}
Let consider the function 
\begin{equation*}
J(\tau)=\displaystyle\int_{\tau}^{\infty}g (s) \left(e^{\mu_0 s}-1\right)\,ds-\frac{1}{\mu_1 -\mu_0}e^{\mu_0\tau}g(\tau).
\end{equation*}
We have, using \eqref{g} and \eqref{gexp}$_1$, $\lim_{\tau\to\infty} J(\tau)=0$ and
\begin{equation*}
\begin{array}{lll}
J^{\prime}(\tau) & = & -g (\tau) \left(e^{\mu_0 \tau}-1\right)-\frac{\mu_0}{\mu_1 -\mu_0}e^{\mu_0 \tau}g(\tau)-\frac{1}{\mu_1 -\mu_0}e^{\mu_0 \tau}g^{\prime}(\tau) \\
\\
& \geq & -g (\tau) \left(e^{\mu_0 \tau}-1\right)-\frac{\mu_0}{\mu_1 -\mu_0}e^{\mu_0 \tau}g(\tau)+\frac{\mu_1}{\mu_1 -\mu_0}e^{\mu_0 \tau}g(\tau) \\
\\
& \geq & g(\tau)\geq 0;
\end{array} 
\end{equation*} 
that is $J$ is non-decreasing, then $J(\tau)\leq 0$, so 
\begin{equation*}
e^{-\mu_0 \tau}\displaystyle\int_{\tau}^{\infty}g (s) \left(e^{\mu_0 s}-1\right)\,ds\leq \frac{1}{\mu_1 -\mu_0}g(\tau) ,
\end{equation*}
and consequently, \eqref{fubini} implies that 
\begin{equation*}
\displaystyle\int_{0}^{\infty}g (s)\left\Vert \displaystyle\int_0^{s} f_{8x} (\tau)d\tau\right\Vert^{2} \,ds\leq \frac{1}{\mu_0 (\mu_1 -\mu_0)} \Vert f_8\Vert_{L_g}^2 <\infty, 
\end{equation*}
since $f_8\in L_g$, which implies that
\begin{equation}
s\mapsto \displaystyle\int_0^{s} f_8 (\tau)d\tau\,\in L_g .\label{intf8} 
\end{equation}
On the other hand, we have, using again \eqref{gexp}$_2$, \eqref{intf8} and H\"older's inequality, 
\begin{equation*}
\begin{array}{lll}
\left\Vert\displaystyle\int_{0}^{\infty}g (s)\displaystyle\int_{0}^{s}f_{8x}(\tau)d\tau ds\right\Vert^2 & = & \left\Vert\displaystyle\int_{0}^{\infty}{\sqrt{g (s)}}{\sqrt{g (s)}}\displaystyle\int_{0}^{s}f_{8x}(\tau)d\tau ds\right\Vert^2\\
\\
& \leq & \left(\displaystyle\int_{0}^{\infty}g (s)ds\right)\displaystyle\int_{0}^{\infty}g (s) \left\Vert\displaystyle\int_0^{s} f_{8x} (\tau)d\tau\right\Vert^{2} ds \\
\\
& \leq & g_0\left\Vert s\mapsto \displaystyle\int_0^{s}f_{8}(\tau)d\tau\right\Vert_{L_g}^{2} <\infty,
\end{array}
\end{equation*}
then 
\begin{equation*}
\displaystyle\int_0^{\infty} g(s)\displaystyle\int_0^{s}f_8 (\tau)d\tau ds\in H_0^1 (0,1), 
\end{equation*}
and so, using \eqref{thetaf8}, we find that $\theta\in H_0^1 (0,1)$. Moreover, from \eqref{gexp}$_2$, we observe that
\begin{equation*}
\displaystyle\int_{0}^{\infty}g(s)\Vert s^m\theta_x\Vert^{2} ds=\Vert 
\theta_x\Vert^{2} \displaystyle\int_{0}^{\infty}s^{2m}g(s)ds=g_{2m}\Vert 
\theta_x\Vert^{2}<\infty,\quad m=0,1,
\end{equation*} 
thus
\begin{equation}
s\mapsto \theta,\,\,s\mapsto s\theta\,\in L_g .\label{thetas} 
\end{equation} 
The properties \eqref{etathetaf8}, \eqref{intf8} and \eqref{thetas} lead to 
$\eta,\,\eta_s\in L_g$.
\vskip0,1truecm
Finally, \eqref{ZF} has a solution $\Phi \in D\left( \mathcal{A}\right)$ if
there exists 
\begin{equation}
(\varphi ,\psi ,w)\in (H^2 (0,1)\cap H_k)\times(H^2 (0,1)\cap H_{1-k})\times(H^2 (0,1)\cap H_{1-k})\label{varphipsiw}
\end{equation}
satisfying 
\begin{equation}
\frac{\partial^k}{\partial x^k}\varphi (0)=\frac{\partial^k}{\partial x^k}
\varphi (1)=\frac{\partial^{1-k}}{\partial x^{1-k}}\psi (0)=\frac{\partial^{1-k}}{\partial x^{1-k}}\psi (1)=\frac{\partial^{1-k}}{\partial x^{1-k}}w (0)=\frac{\partial^{1-k}}{\partial x^{1-k}}w (1)=0\label{bc01}
\end{equation}
and the equations \eqref{ZF}$_2$, \eqref{ZF}$_4$ and \eqref{ZF}$_6$. Assuming that such unknown $(\varphi ,\psi ,w)$ exists, then, multiplying 
\eqref{ZF}$_2$, \eqref{ZF}$_4$ and \eqref{ZF}$_6$ by 
${\hat{\varphi}} \in H_k$ and ${\hat{\psi}},\,{\hat{w}} \in H_{1-k}$, respectively, inegrating by parts and using
\eqref{bc01}, we remark that $(\varphi ,\psi ,w)$ is a solution of the variational formulation
\begin{equation}
B \left((\varphi ,\psi ,w), ({\hat{\varphi}} ,{\hat{\psi}} ,{\hat{w}} )\right)=C ({\hat{\varphi}} ,{\hat{\psi}} ,{\hat{w}} ),\,\,\forall 
({\hat{\varphi}} ,{\hat{\psi}} ,{\hat{w}}) \in \mathcal{H}_k , \label{z7f7}
\end{equation}
where $B$ is a bilinear form over $\mathcal{H}_k\times \mathcal{H}_k$ given by 
\begin{equation*}
\begin{array}{lll}
B \left((\varphi ,\psi ,w), ({\hat{\varphi}} ,{\hat{\psi}} ,{\hat{w}} )\right) &=& k_1\left\langle \varphi_{x} +\psi +lw ,{\hat{\varphi}}_x +{\hat{\psi}}+l{\hat{w}} \right\rangle
+k_2\left\langle \psi_{x},{\hat{\psi}}_x \right\rangle +k_{3}\left\langle w_{x}-l\varphi, {\hat{w}}_x -l{\hat{\varphi}}\right\rangle\\
\\
&=& \left\langle (\varphi ,\psi ,w), ({\hat{\varphi}} ,{\hat{\psi}} ,{\hat{w}} )\right\rangle_{\mathcal{H}_k} 
\end{array}
\end{equation*}
and $C$ is a linear form over $\mathcal{H}_k$ defined by
\begin{equation*}
C ( {\hat{\varphi}} ,{\hat{\psi}} ,{\hat{w}} ) =-\left\langle \delta \frac{\partial^k}{\partial x^k}\theta +f_{2} ,{\hat{\varphi}}\right\rangle-\left\langle f_4 ,{\hat{\psi}} \right\rangle -\left\langle f_6 ,{\hat{w}} \right\rangle.
\end{equation*}
According to the fact that $\frac{\partial^k}{\partial x^k}\theta ,f_{2}, f_4, f_6 \in L^2 (0,1)$ and because $(\mathcal{H}_k ,\left\langle ,\right\rangle_{\mathcal{H}_k})$ is a Hilbert space, it is easy to see that $B$ is continuous and coercive, and $C$ is continuous. Then, the Lax-Milgram theorem implies that \eqref{z7f7} has a unique solution
\begin{equation*}
(\varphi ,\psi ,w)\in \mathcal{H}_k.
\end{equation*}
Therefore, using classical elliptic regularity arguments, we conclude that 
$(\varphi ,\psi ,w)$ satisfies \eqref{ZF}$_2$, \eqref{ZF}$_4$, 
\eqref{ZF}$_6$, \eqref{varphipsiw} and \eqref{bc01}. This proves that 
\eqref{ZF} admits a unique solution $\Phi\in D\left( \mathcal{A}\right)$. By the resolvent identity, we have $\lambda I -\mathcal{A}$ is surjective, for any $\lambda >0$ (see \cite{liu1}), where $I$ is the identity operator. Consequently, the Lumer-Phillips theorem implies that $\mathcal{A}$ is the infinitesimal generator of a linear $C_{0}$-semigroup of contractions on 
$\mathcal{H}$. The semigroup theory guarantees the next theorem (see \cite{pazy}).
\vskip0,1truecm
\begin{theorem}\label{Theorem 1.1}
Under assumptions \eqref{l}-\eqref{f+}, and for any $m\in \mathbb{N}$ and $\Phi_0 \in D(\mathcal{A}^m)$, system \eqref{syst_2} admits a unique solution
\begin{equation}
\Phi\in \cap_{j=0}^m C^{m-j} \left(\mathbb{R}_{+} ;D\left(\mathcal{A}^j\right)\right) . \label{exist}
\end{equation}
\end{theorem}

\section{Lack of exponential stability}

The subject of this section is showing that \eqref{syst_2} is not exponentailly stable in case $k=0$ (system \eqref{cdt_10}-\eqref{cdt_1}) independently of the values of the coefficients, and in case $k=1$ (system \eqref{syst1}-\eqref{cdt_10}) depending on the following connections:
\begin{equation}
\forall m\in  \mathbb{Z},\,\,l^2\ne
\frac{k_3 -k_2 }{k_3} \left(m\pi\right)^2 - \dfrac{k_1}{k_1 +k_3},\label{lpi}
\end{equation}
\begin{equation}
k_2 \ne k_3 \label{k2k3}
\end{equation}
and 
\begin{equation}
\delta^2 =\frac{(k_2 -k_1 )(k_2 -g_0 )}{k_2} =\frac{(k_3 -k_1 )(k_3 -g_0 )}{k_3} .\label{k23delta}
\end{equation}
\vskip0,1truecm
\begin{theorem}\label{Theorem 3.1}
Under the assumptions \eqref{l}-\eqref{f+}, we have the following exponential  stability results:
\vskip0,1truecm
1. System \eqref{syst_2} in case $k=0$ is not exponentially stable.
\vskip0,1truecm
2. System \eqref{syst_2} in case $k=1$ is not exponentially stable if
\eqref{lpi} or \eqref{k2k3} or \eqref{k23delta} is not satisfied.    
\end{theorem}
\vskip0,1truecm
\begin{proof} 
It is known that the exponential stability is equivalent to (see \cite{huan} and \cite{prus}) 
\begin{equation}
i\mathbb{R} \subset \rho\left( \mathcal{A}\right)\quad\hbox{and}\quad \sup_{\lambda\in \mathbb{R}}\left\Vert \left( i\lambda I-\mathcal{A}\right) ^{-1} \right\Vert_{\mathcal{L}\left( \mathcal{H}\right) } <\infty. \label{expon}
\end{equation}
We start by proving that the first condition in \eqref{expon} is equivalent to \eqref{lpi}. In the previous section, we have proved that 
$0\in\rho \left( \mathcal{A}\right)$. Moreover, $\mathcal{A}^{-1}$ is bounded and it is a bijection between $\mathcal{H}$ and $D(\mathcal{A})$. then, because $D(\mathcal{A})$ has a compact embedding into $\mathcal{H}$, 
$\mathcal{A}^{-1}$ is a compact operator, which implies that the spectrum of $\mathcal{A}$ is discrete and contains only eigenvalues. Let $\lambda\in\mathbb{R}^*$ and 
\begin{equation*}
\Phi = \left(\varphi ,{\tilde{\varphi}} ,\psi ,{\tilde{\psi}} ,w ,{\tilde{w}} ,\theta, \eta\right)^T\in D(\mathcal{A}) .
\end{equation*}
We have to prove that $i\lambda$ is not an eigenvalue of $\mathcal{A}$ if and only if \eqref{lpi} holds; that is, \eqref{lpi} is equivalent to the fact that $\Phi=0$ is the unique solution of the equation
\begin{equation}
\mathcal{A}\,\Phi =i\,\lambda\,\Phi.  \label{eq_3_4*}
\end{equation}
According to \eqref{A}, equation \eqref{eq_3_4*} is equivalent to
\begin{equation}
\left\{
\begin{array}{l}
{\tilde{\varphi}}-i\lambda\varphi = {\tilde{\psi}}-i\lambda\psi= {\tilde{w}}-i\lambda w =0, \vspace{0.2cm}\\
k_1\left( \varphi _{x}+\psi +l\,w\right) _{x}+lk_3\left( w_{x}-l\varphi \right) -\delta\frac{\partial^k}{\partial x^k}\theta =i\lambda {\tilde{\varphi}} ,\vspace{0.2cm}\\
k_2\psi _{xx}-k_1\left( \varphi _{x}+\psi
+l\,w\right) =i\lambda {\tilde{\psi}}, \vspace{0.2cm}\\
k_3\left( w_{x}-l\varphi \right) _{x}-lk_1\left( \varphi _{x}+\psi +l\,w\right) =i\lambda {\tilde{w}}, \vspace{0.2cm}\\
\displaystyle\int_0^{\infty}g\eta_{xx}ds+(-1)^k\delta \frac{\partial^k}{\partial x^k}{\tilde{\varphi}}=i\lambda\theta ,\vspace{0.2cm}\\
\theta -\eta_s =i\lambda\eta.
\end{array}
\right. \label{eq_3_5}
\end{equation}
Using \eqref{g}, \eqref{dissp} and \eqref{eq_3_4*}, we find 
\begin{equation*}
0=Re\,i\lambda\Vert\Phi\Vert_{\mathcal{H}}^2 =Re\,\left\langle i\lambda \Phi,\Phi \right\rangle _{\mathcal{H}} =Re\,\left\langle \mathcal{A}\Phi ,\Phi \right\rangle _{\mathcal{H}}= \frac{1}{2}\displaystyle\int_0^{\infty} g' \Vert\eta_x\Vert^{2} ds\leq 
-\frac{\mu_1}{2}\displaystyle\int_0^{\infty} g \Vert\eta_x\Vert^{2} ds=
-\frac{\mu_1}{2}\Vert\eta\Vert_{L_g}^{2} ,
\end{equation*}
so $\eta =0$. Therefore, \eqref{eq_3_5}$_6$ implies that $\theta =0$, and moreover, from \eqref{eq_3_5}$_5$ we get 
$\frac{\partial^k}{\partial x^k}{\tilde{\varphi}}=0$. This means that  
${\tilde{\varphi}}=0$ if $k=0$. And if $k=1$, it follows that 
${\tilde{\varphi}}$ is a constant, thus, using the definition of $H_1$,
${\tilde{\varphi}}=0$. Now, \eqref{eq_3_5}$_1$ leads to $\varphi=0$, and then        
\eqref{eq_3_5} is reduced to 
\begin{equation}
\left\{
\begin{array}{l}
{\tilde{\psi}}-i\lambda\psi ={\tilde{w}}-i\lambda w=0, \vspace{0.2cm}\\
k_1\psi _{x}+l\left( k_1+k_3\right) w_{x}=0, \vspace{0.2cm}\\
k_2\psi _{xx}-k_1\left( \psi +l\,w\right) =-\lambda^{2}\psi ,\vspace{0.2cm}\\
k_{3}w_{xx}-lk_1\left( \psi +l\,w\right) =-\lambda^{2}w.
\end{array}
\right.  \label{eq_3_1100}
\end{equation}
From \eqref{eq_3_1100}$_2$, we see that $k_1\psi +l\left( k_1 +k_{3}\right)w$ is a constant, then, thanks to the definition of $H_{1-k}$, we get 
\begin{equation}
\psi =-l\left( 1+\dfrac{k_{3}}{k_1}\right) w. \label{eq_3_1200}
\end{equation}
Combining \eqref{eq_3_1100}$_3$ and \eqref{eq_3_1100}$_4$, we arrive at
\begin{equation}
lk_2\psi _{xx}-k_{3}w_{xx}=-l\lambda^{2} \psi +\lambda^{2}w.  \label{eq_3_1300}
\end{equation}
Exploiting \eqref{eq_3_1200} and \eqref{eq_3_1300}, we find
\begin{equation*}
w_{xx}+\frac{l^2\left( k_1 +k_{3}\right) +k_1}{ k_2 l^2\left( k_1 +k_{3}\right) +k_1 k_{3}} \lambda^{2} w=0,
\end{equation*}
which implies that, for some constants $c_1$ and $c_2$,
\begin{equation*}
w(x)=c_1\cos\,\left(\lambda{\sqrt{\frac{l^2\left( k_1 +k_{3}\right) +k_1}{ k_2 l^2\left( k_1 +k_{3}\right) +k_1 k_{3}}}}x\right) + c_2\sin\,\left(\lambda{\sqrt{\frac{l^2\left( k_1 +k_{3}\right) +k_1}{ k_2 l^2\left( k_1 +k_{3}\right) +k_1 k_{3}}}}x\right).
\end{equation*}
The boundary conditions at $x=0$ in \eqref{cdt_100} and \eqref{cdt_1} lead to $c_1 =0$ in case \eqref{cdt_100}, and $c_2 =0$ in case \eqref{cdt_1}, and then
\begin{equation}
w(x)=c_2\sin\,\left(\lambda{\sqrt{\frac{l^2\left( k_1 +k_{3}\right) +k_1}{ k_2 l^2\left( k_1 +k_{3}\right) +k_1 k_{3}}}} x\right) \label{allam+000}
\end{equation}
if $k=1$, and 
\begin{equation}
w(x)=c_1\cos\,\left(\lambda{\sqrt{\frac{l^2\left( k_1 +k_{3}\right) +k_1}{ k_2 l^2\left( k_1 +k_{3}\right) +k_1 k_{3}}}}x\right) \label{allam+00}
\end{equation}
if $k=0$. Then the boundary conditions at $x=1$ in the definition of $D(\mathcal{A})$ lead to
\begin{equation}
c_1 =c_2 =0\quad\hbox{or}\quad\exists m\in\mathbb{Z} :\,\,\lambda{\sqrt{\frac{l^2\left( k_1 +k_{3}\right)+k_1}{ k_2 l^2\left( k_1 +k_{3}\right)+k_1 k_{3}}}} =m\pi .
\label{c2mpiZ}
\end{equation}
If $c_1 =c_2 =0$, then \eqref{eq_3_1200}, \eqref{allam+000} and 
\eqref{allam+00} imply that $\psi=w=0$, and therefore, we conclude from \eqref{eq_3_1100}$_1$ that ${\tilde\psi}={\tilde w}=0$. Consequently, 
$\Phi =0$. If $c_1 \ne 0$ and $c_2 \ne 0$, then \eqref{eq_3_1200}, 
\eqref{allam+000} and \eqref{allam+00} imply that \eqref{eq_3_1100}$_3$ and \eqref{eq_3_1100}$_4$ are equivalent to
\begin{equation}
\left(k_3 -k_2\right)\lambda^2 =\frac{k_3}{k_1 +k_3}\left[k_2l^2 \left(k_1 +k_3\right)+k_1 k_3\right]. \label{eq_3_14200}
\end{equation}   
By using the second assertion in \eqref{c2mpiZ}, it follows that \eqref{lpi} is not satisfied. Finally, we see that, if \eqref{lpi} holds, then necessarily $c_1 =c_2 =0$, which leads to $\Phi =0$. Otherwise, if \eqref{lpi} does not hold, one can define $\psi$ by \eqref{eq_3_1200}, and 
$w$ by \eqref{allam+000} and \eqref{allam+00}, for any constants $c_1$ and $c_2$, which implies that \eqref{eq_3_4*} has an infinite solutions. This ends the proof of the equivalence between \eqref{lpi} and the first condition in \eqref{expon}.   
\vskip0,1truecm
Now, we prove that the second condition in \eqref{expon} is not satisfied in case $k=0$, and in case $k=1$ if \eqref{k2k3} or \eqref{k23delta} is not satisfied. To do so, we prove that there exists a sequence 
$(\lambda_{n})_n \subset \mathbb{R}$ such that
\begin{equation*}
\lim_{n\to \infty}\left\Vert \left( i\lambda _{n}I-\mathcal{A}\right) ^{-1}\right\Vert_{\mathcal{L}\left( \mathcal{H}\right)} =\infty.
\end{equation*}
This is equivalent to prove that there exists a sequence 
$(F_{n})_n \subset \mathcal{H}$ satisfying
\begin{equation}
\left\Vert F_{n}\right\Vert_{\mathcal{H}}:=\left\Vert (f_{1,n},\cdots,f_{8,n})^T\right\Vert_{\mathcal{H}}\leq 1,\quad\forall n\in\mathbb{N} \label{Fn}
\end{equation}
and
\begin{equation}
\lim_{n\to \infty}\Vert\left( i\lambda _{n}I-\mathcal{A}\right)^{-1}F_{n}\Vert_{\mathcal{H}} = \infty.  \label{eq_3}
\end{equation}
For this purpose, let
\begin{equation*}
\Phi_{n}:=\left(\varphi_n ,{\tilde{\varphi}}_n ,\psi_n ,{\tilde{\psi}}_n ,w_n ,{\tilde{w}}_n ,\theta_n ,\eta_n\right)^T =\left( i\lambda _{n}I-\mathcal{A}\right)^{-1} F_n, \quad\forall n\in\mathbb{N}.
\end{equation*}
Then, we have to prove that there exists a sequence 
$(\Phi _{n})_n \subset D\left(\mathcal{A}\right)$ such that \eqref{Fn} holds, 
\begin{equation}
\lim_{n\to \infty}\Vert\Phi _{n}\Vert_{\mathcal{H}}=\infty\quad\hbox{and}\quad i\lambda _{n}\Phi _{n}-\mathcal{A}\Phi _{n}=F_{n} ,  \,\,\forall n\in\mathbb{N}. \label{eq_4}
\end{equation}
Therefore, from \eqref{A} and the second equality in \eqref{eq_4}, we have 
\begin{equation}
\left\{
\begin{array}{l}
i\lambda _{n}\varphi _{n}-{\tilde{\varphi}}_{n}=f_{1,n} ,\vspace{0.2cm}\\
i\lambda _{n}{\tilde{\varphi}}_{n}-k_1\left( \varphi _{n,x}+\psi _{n}+l\,w_{n}\right)_{x}-lk_{3}\left( w_{n,x}-l\varphi _{n}\right) +\delta \frac{\partial^k}{\partial x^k}\theta _{n}=f_{2,n} ,\vspace{0.2cm}\\
i\lambda _{n}\psi _{n}-{\tilde{\psi}}_{n} =f_{3,n} ,\vspace{0.2cm}\\
i\lambda _{n}{\tilde{\psi}}_{n}-k_2 \psi _{n,xx}+k_1\left( \varphi _{n,x}+\psi_{n}+l\,w_{n}\right) =f_{4,n} ,\vspace{0.2cm}\\
i\lambda _{n}w_{n}-{\tilde{w}}_{n}=f_{5,n} ,\vspace{0.2cm}\\
i\lambda _{n}{\tilde{w}}_{n}-k_{3}\left( w_{n,x}-l\varphi _{n}\right)
_{x}+lk_1\left( \varphi _{n,x}+\psi _{n}+l\,w_{n}\right) =f_{6,n} ,\vspace{0.2cm}\\
i\lambda _{n}\theta_n -\displaystyle\int_0^{\infty}g\eta_{n,xx}ds-(-1)^k\delta \frac{\partial^k}{\partial x^k}{\tilde{\varphi}}_{n}=f_{7,n} ,\vspace{0.2cm}\\
i\lambda _{n}\eta_n +\eta_{n,s} -\theta_n =f_{8,n} .  
\end{array}
\right.  \label{eq_2_5}
\end{equation}
We put $N=(n+1)\pi$ and we choose
\begin{equation}
\left\{
\begin{array}{l}
{\tilde{\varphi}}_{n} =i\lambda_{n}\varphi_{n} ,\quad {\tilde{\psi}}_{n} =i\lambda_{n}\psi_{n} ,\quad {\tilde{w}}_{n} =i\lambda _{n} w_n ,\quad
\eta_n (x,s)=\frac{i}{\lambda _{n}}\left(e^{-i\lambda _{n} s}-1\right)\theta_n (x),\vspace{0.2cm}\\
f_{1,n}=f_{3,n}=f_{5,n} =f_{8,n} =0. 
\end{array}
\right.  \label{fn0}
\end{equation}
Then \eqref{eq_2_5}$_1$, 
\eqref{eq_2_5}$_3$, \eqref{eq_2_5}$_5$ and 
\eqref{eq_2_5}$_8$ are satisfied. On the other hand, choosing  
\begin{equation}
\left\{
\begin{array}{l}
\varphi_{n} (x)=\alpha_{1,n}\cos\,(Nx) ,\quad \psi_{n} (x) =-\alpha_{2,n}\sin\,(Nx) ,\quad w_{n} (x) = -\alpha_{3,n}\sin\,(Nx), \vspace{0.2cm}\\
\theta_n =\alpha_{4,n}\sin\,(Nx),\quad
f_{2,n} (x)=\beta_{2,n}\cos\,(Nx) ,\quad f_{4,n} (x) =-\beta_{4,n}\sin\,(Nx) , \vspace{0.2cm}\\ 
f_{6,n} (x)= -\beta_{6,n}\sin\,(Nx),\quad f_{7,n} =0
\end{array}
\right.  \label{fn01}
\end{equation} 
when $k=1$ (that is, for \eqref{syst1}), and
\begin{equation}
\left\{
\begin{array}{l}
\varphi_{n} (x)=\alpha_{1,n}\sin\,(Nx),\quad \psi_{n} (x)=\alpha_{2,n}\cos\,(Nx),\quad w_{n} (x)=\alpha_{3,n}\cos\,(Nx),\vspace{0.2cm}\\
\theta_n =0,\quad f_{2,n} (x)=\beta_{2,n}\sin\,(Nx) ,\quad f_{4,n} (x) =\beta_{4,n}\cos\,(Nx) , \vspace{0.2cm}\\
f_{6,n} (x)= \beta_{6,n}\cos\,(Nx),\quad f_{7,n} (x)= \beta_{7,n}\sin\,(Nx), 
\end{array}
\right.  \label{fn000}
\end{equation} 
when $k=0$ (that is, for \eqref{syst11}), where $\alpha_{j,n}$ and 
$\beta_{j,n}$ are constants. Let us put
\begin{equation}
\mu_{k,n}=\left\{
\begin{array}{ll}
\frac{\delta^2 \lambda _{n}^2}{\lambda _{n}^2 -N^2 (g_0 -\mu_{2,n})} &  \hbox{if}\,\,k=1,\vspace{0.2cm}\\
0 & \hbox{if}\,\,k=0
\end{array}
\right. \quad \hbox{and}\quad\mu_{2,n}=\displaystyle\int_0^{\infty}g(s)
e^{-i\lambda _{n} s}ds \label{mukn2n}
\end{equation}
and choose
\begin{equation}
\alpha_{4,n}=\frac{\delta \lambda _{n}^2 N}{\lambda _{n}^2 -N^2 (g_0 -\mu_{2,n})} \alpha_{1,n} \quad \hbox{and}\quad \beta_{7,n}=-i\delta \lambda _{n}\alpha_{1,n}. \label{fn00}
\end{equation}
The choices \eqref{fn0}-\eqref{fn00} guarantee that \eqref{eq_2_5}$_7$ holds,
$\Phi_{n} \in D\left(\mathcal{A}\right)$ and $F_{n} \in \mathcal{H}$.
Moreover, \eqref{eq_2_5}$_2$, \eqref{eq_2_5}$_4$ and \eqref{eq_2_5}$_6$ are reduced to the algebraic system
\begin{equation}
\left\{
\begin{array}{l}
\left[(k_1 +\mu_{k,n}) N^2 -\lambda_{n}^{2}+l^2 k_3 \right]\alpha_{1,n} +k_1 N\alpha_{2,n} +l\left(k_1 +k_3\right)N\alpha_{3,n} =\beta_{2,n}, \vspace{0.2cm}\\
k_1 N\alpha_{1,n} +\left(k_2 N^2 -\lambda _{n}^{2}+k_1\right)\alpha_{2,n} +k_1 l\alpha_{3,n} =\beta_{4,n},\vspace{0.2cm}\\
l\left(k_1 +k_3\right)N\alpha_{1,n} +lk_1 \alpha_{2,n} +\left(k_3 N^2 -\lambda _{n}^{2}+l^2 k_1\right)\alpha_{3,n} =\beta_{6,n} .  
\end{array}
\right. \label{eq_2_6}
\end{equation}
Now, we addapt some arguments of \cite{afas1, ag1, ag2, ag6} to our models. Let us distinguish the next cases.
\vskip0,1truecm
\subsection{Case $k_2 =k_3$} We choose
\begin{equation}
\left\{
\begin{array}{l}
\alpha_{1,n} =\beta_{2,n} =0,\quad\lambda_{n} =N{\sqrt{k_{3}}},\vspace{0.2cm}\\
\beta_{4,n}=-lk_{3}\alpha_{3,n},\quad \beta_{6,n}=-l^{2}k_{3}\alpha_{3,n} . 
\end{array}
\right. \label{fn2}
\end{equation}
Then \eqref{eq_2_6} is reduced to
\begin{equation}
\left\{
\begin{array}{l}
k_1\alpha_{2,n} +l\left( k_1 +k_{3}\right) \alpha_{3,n} =0, \vspace{0.2cm}\\
k_1 \alpha_{2,n} +lk_1\alpha_{3,n} =-lk_{3} \alpha_{3,n}, \vspace{0.2cm}\\
lk_1\alpha_{2,n} +l^{2}k_{1} \alpha_{3,n} =-l^{2}k_{3} \alpha_{3,n} ,
\end{array}
\right. \label{fn4}
\end{equation}
which is equivalent to $\alpha_{2,n} =-\frac{l(k_1 +k_{3})}{k_1}
\alpha_{3,n}$. Choosing $\alpha_{3,n} =\frac{1}{lk_{3} \sqrt{1 +l^2 }}$ and using \eqref{fn0}, \eqref{fn00} and \eqref{fn2}, we obtain
\begin{equation*}
\left\Vert F_{n}\right\Vert_{\mathcal{H}}^{2} =\left\Vert f_{4,n}\right\Vert^{2}+\left\Vert f_{6,n}\right\Vert^{2} 
\leq \beta_{4,n}^2 +\beta_{6,n}^2 =\left(lk_{3}\right)^{2}\left( 1+ l^{2}\right) \alpha_{3,n}^{2} = 1 ,
\end{equation*}
so \eqref{Fn} is satisfied. On the other hand, we have  
\begin{equation*}
\left\Vert\Phi_{n}\right\Vert_{\mathcal{H}}^{2} \geq k_3 \left\Vert w_{nx} -l\varphi_n\right\Vert^{2} =k_3 \left\Vert w_{nx}\right\Vert^{2} =\frac{k_3\alpha_{3,n}^2}{2} N^{2}\int_{0}^{1} \left[1\pm\cos\,\left(2Nx\right)\right]\,dx=\frac{k_3\alpha_{3,n}^2}{2} N^{2} ,
\end{equation*}
hence \eqref{eq_3}, since 
\begin{equation}
\lim_{n\to\infty}\left\Vert \Phi _{n}\right\Vert _{\mathcal{H}} = \infty. \label{fn69}
\end{equation}
\vskip0,1truecm
\subsection{Case $k_2\ne k_{3}$ and $\left\{
\begin{array}{ll}
\delta^2 \ne\frac{(k_2 -k_1 )(k_2 -g_0 )}{k_2} & \hbox{if}\,\,k=1,\vspace{0.2cm}\\
k_2\ne k_{1} & \hbox{if}\,\,k=0
\end{array}
\right.$} Let pick 
$\beta_{4,n}:=k_1 \beta_{4}$, with 
$\beta_{4}\in (0,\infty)$ not depending on $n$, and choose
\begin{equation}
\beta_{2,n} = \beta_{6,n} =0\quad\hbox{and}\quad
\lambda_n =\sqrt{k_2 N^2 +k_1} .\label{fn0650+0}
\end{equation}
We have, thanks to \eqref{fn0650+0}, \eqref{eq_2_6} is equivalent to
\begin{equation}
\left\{
\begin{array}{l}
\left[\left(k_1 -k_2 +\mu_{k,n}\right)N^2 +l^2 k_3 -k_1 \right]\alpha_{1,n} +k_1 N\alpha_{2,n} +l\left(k_1 +k_3\right)N\alpha_{3,n} =0, \vspace{0.2cm}\\
k_1 N\alpha_{1,n} +k_1 l\alpha_{3,n} =k_1 \beta_4 ,\vspace{0.2cm}\\
l\left(k_1 +k_3\right)N\alpha_{1,n} +lk_1\alpha_{2,n} +\left[\left(k_3 -k_2 \right)N^2 +l^2 k_1 -k_1 \right]\alpha_{3,n} =0.
\end{array}
\right.  \label{eq_2_610}
\end{equation}
It is clear that \eqref{eq_2_610}$_2$ and \eqref{eq_2_610}$_3$ are equivalent to
\begin{equation}
\alpha_{3,n} =-\frac{1}{l} N\alpha_{1,n} +\dfrac{\beta_{4}}{l} \quad\hbox{and}\quad \alpha_{2,n} =a_{1,n} N\alpha_{1,n} +\beta_{4} a_{2,n} ,\label{eq_2_630}
\end{equation}
where
\begin{equation}
a_{1,n} =\frac{1}{l^2 k_1}\left[\left(k_3 -k_2\right)N^2 +l^2 k_1 -k_1 \right] -\frac{k_1 +k_3}{k_1} \quad\hbox{and}\quad
a_{2,n} =\frac{1}{l^2 k_1}\left[\left(k_2 -k_3\right)N^2 +k_1 -l^2 k_1\right] , \label{a12n0}
\end{equation}
and then \eqref{eq_2_610}$_1$ is satisfied if and only if
\begin{equation}
\alpha_{1,n} = \frac{\beta_{4}\left[k_1 a_{2,n} +k_1 +k_3 \right]N}{\left(k_2 +k_3 -\mu_{k,n} -k_1 a_{1,n} \right)N^2 +k_1 -l^2 k_3},\label{eq_2_670}
\end{equation} 
therefore, \eqref{eq_2_630} and \eqref{eq_2_670} mean that
\begin{equation}
\alpha_{2,n} = \beta_{4}\frac{\left[(k_1 +k_3 )a_{1,n} +(k_2 +k_3 -\mu_{k,n})a_{2,n}\right]N^2 +(k_1 -l^2 k_3 )a_{2,n}}{\left(k_2 +k_3 -\mu_{k,n} -k_1 a_{1,n} \right)N^2 +k_1 -l^2 k_3}.\label{alpha2n0}
\end{equation}
On the other hand, because $\lambda_n\in \mathbb{R}$ and according to \eqref{g} and \eqref{gexp}$_1$, we have
\begin{equation*}
\lim_{s\to\infty} g(s)=\lim_{s\to\infty} g^{\prime}(s)=0\quad\hbox{and}\quad\vert e^{-i\lambda_n s} \vert=1,
\end{equation*} 
then, integrating by parts, we get
\begin{equation*}
\mu_{2,n} = \dfrac{1}{i\lambda_n}\left(g(0)+\displaystyle\int_0^{\infty} g^{\prime} (s)e^{-i\lambda_n s} \,ds\right)=\frac{1}{i\lambda_n}\left(g(0)+\frac{1}{i\lambda_n} g^{\prime} (0)+\dfrac{1}{i\lambda_n}\displaystyle\int_0^{\infty} g^{\prime\prime} (s)e^{-i\lambda_n s} \,ds\right),
\end{equation*} 
so, thanks to the second property in \eqref{g0}, we find that ($\sim$ means: asymptotically equal, for $n$ large) 
\begin{equation}
\mu_{2,n}\sim \frac{g(0)}{i\lambda_n}, \label{mu2n}  
\end{equation}
and therefore
\begin{equation}
\mu_{1,n}\sim \left\{
\begin{array}{ll} 
\frac{\delta^2 k_2}{k_2 -g_0} & \hbox{if}\quad g_0\ne k_2, \vspace{0.2cm}\\
\frac{i\delta^2 k_2{\sqrt{k_2}}N}{g(0)} & \hbox{if}\quad g_0 = k_2 .
\end{array}
\right. \label{mun00}  
\end{equation}
Thus, we conclude from \eqref{alpha2n0} that
\begin{equation}
\vert\alpha_{2,n}\vert \sim \left\{
\begin{array}{ll}
\frac{\vert k_2 -k_1\vert \beta_{4}}{k_1} & \hbox{if}\quad
k=0, \vspace{0.2cm}\\ 
\frac{\delta^2 k_2{\sqrt{k_2}}\beta_{4} N}{k_1 g(0)} & \hbox{if}\quad k=1
\,\,\hbox{and}\,\,g_0 = k_2 , \vspace{0.2cm}\\ 
\frac{\left\vert(k_2 -k_1 )(k_2 -g_0 )-\delta^2 k_2\right\vert\beta_4}{k_1 \vert k_2 -g_0 \vert } & \hbox{if}\quad k=1\,\,\hbox{and}\,\,g_0 \ne k_2 .
\end{array}
\right.  \label{eq_2_71000} 
\end{equation}
Now, we observe that 
\begin{equation}
\left\Vert \Phi _{n}\right\Vert _{\mathcal{H}}^{2} \geq k_2 \left\Vert \psi_{n,x}\right\Vert^{2} = \frac{k_2 N^2}{2}\vert\alpha_{2,n}\vert^{2}\displaystyle\int_{0}^{1} \left[1\pm\cos\,\left(2Nx\right)\right]\,dx =\frac{k_2 N^2}{2}\vert\alpha_{2,n}\vert^{2} , \label{eq_2_710006}
\end{equation}
then, by \eqref{eq_2_71000}, we get \eqref{fn69}. Moreover, for $k=1$, we take 
$\beta_4 =1$ and we find
\begin{equation}
\left\Vert F_{n}\right\Vert _{\mathcal{H}}^{2} =\left\Vert f_{4,n}\right\Vert^{2} =\beta_{4}^2 \displaystyle\int_{0}^{1} \sin^2 \,(Nx) dx\leq 1 ,
\label{eq_2_7100066}
\end{equation}  
which implies \eqref{Fn}. For $k=0$, we have, using \eqref{fn00} and \eqref{eq_2_670}, 
\begin{equation}
\left\Vert F_{n}\right\Vert_{\mathcal{H}}^{2} =\left\Vert f_{4,n}
\right\Vert^{2} +\left\Vert f_{7,n}\right\Vert^{2} = \beta_{4}^2 \displaystyle\int_{0}^{1} \cos^2 \,(Nx) dx+\vert\beta_{7,n}\vert^2 \displaystyle\int_{0}^{1} \sin^2 \,(Nx) dx\leq \beta_{4}^2 (1+a_{3,n} ), \label{fn065+0}
\end{equation} 
where
\begin{equation*}
a_{3,n} =\frac{\delta^2\left( k_1 a_{2,n} +k_1 +k_3 \right)^2 N^2 \lambda_n^2}{\left[ \left(k_2 +k_3 -k_1 a_{1,n} \right)N^2 +k_1 -l^2 k_3\right]^2}.
\end{equation*}
According to \eqref{a12n0}, we see that the sequence $(a_{3,n})_n$ is bounded, and then we can choose    
\begin{equation}
\beta_4 =\frac{1}{\sqrt{1+\sup_{n\in \mathbb{N}}a_{3,n}}} \label{beta4fn065+0}
\end{equation}
and conclude from \eqref{fn065+0} that \eqref{Fn} is valid.
\vskip0,1truecm
\subsection{Case $k_2\ne k_{3}$ and $\left\{
\begin{array}{ll}
\delta^2 \ne \frac{(k_3 -k_1 )(k_3 -g_0 )}{k_3} & \hbox{if}\,\,k=1,\vspace{0.2cm}\\
k_2 = k_{1} & \hbox{if}\,\,k=0
\end{array}
\right.$} Let pick $\beta_{6,n}:=\beta_{6}\in (0,\infty)$ not depending on 
$n$ and choose
\begin{equation}
\beta_{2,n} = \beta_{4,n} =0\quad\hbox{and}\quad
\lambda_n =\sqrt{k_3 N^2 +l^2 k_1} .\label{fn0650+}
\end{equation}
We have, thanks to \eqref{fn0650+}, \eqref{eq_2_6} is equivalent to
\begin{equation}
\left\{
\begin{array}{l}
\left[\left(k_1 -k_3 +\mu_{k,n}\right)N^2 +l^2 (k_3 -k_1 )\right]\alpha_{1,n} +k_1 N\alpha_{2,n} +l\left(k_1 +k_3\right)N\alpha_{3,n} =0, \vspace{0.2cm}\\
k_1 N\alpha_{1,n} +\left[(k_2 -k_3 ) N^2 +k_1 (1-l^2 )\right]\alpha_{2,n} +k_1 l\alpha_{3,n} =0,\vspace{0.2cm}\\
l\left(k_1 +k_3\right)N\alpha_{1,n} +lk_1\alpha_{2,n} =\beta_{6} .
\end{array}
\right.  \label{eq_2_61}
\end{equation}
We see that \eqref{eq_2_61}$_2$ and \eqref{eq_2_61}$_3$ hold if
\begin{equation}
\alpha_{2,n} =-\frac{k_1 +k_3}{k_1} N\alpha_{1,n} +\dfrac{\beta_{6}}{lk_1} \quad\hbox{and}
\quad \alpha_{3,n} =a_{1,n} N\alpha_{1,n} +\beta_{6} a_{2,n} ,\label{eq_2_63}
\end{equation}
where
\begin{equation}
a_{1,n} =\frac{k_1 +k_3}{lk_1^2}\left(k_2 -k_3\right)N^2 +\frac{k_3}{lk_1}-\frac{l (k_1 +k_3)}{k_1} \quad\hbox{and}\quad
a_{2,n} =\frac{1}{(lk_1)^2}\left[\left(k_3 -k_2\right)N^2 +l^2 k_1 -k_1\right] , \label{a12n}
\end{equation}
and \eqref{eq_2_61}$_1$ is satisfied if
\begin{equation}
\alpha_{1,n} = \frac{\beta_{6}\left[l\left(k_1 +k_3\right)a_{2,n} +\frac{1}{l}\right]N}{\left[2k_3 -\mu_{k,n} -l\left(k_1 +k_3\right)a_{1,n} \right]N^2 +l^2 \left(k_1 -k_3\right)}.\label{eq_2_67}
\end{equation} 
As in the previous case, using \eqref{mu2n}, we have
\begin{equation}
\mu_{1,n}\sim \left\{
\begin{array}{ll} 
\frac{\delta^2 k_3}{k_3 -g_0} & \hbox{if}\quad g_0\ne k_3, \vspace{0.2cm}\\
\frac{i\delta^2 k_3{\sqrt{k_3}}N}{g(0)} & \hbox{if}\quad g_0 = k_3.
\end{array}
\right. \label{mun}  
\end{equation}
Noting that \eqref{a12n}-\eqref{eq_2_67} imply that
\begin{equation}
\vert\alpha_{1,n}\vert \sim \frac{\beta_{6}}{l(k_1 +k_3 )N}, \label{eq_2_7100}
\end{equation}
and \eqref{eq_2_63} and \eqref{eq_2_7100} lead to 
\begin{equation}
\vert\alpha_{3,n}\vert \sim \left\{
\begin{array}{ll} 
\frac{\vert k_3 -k_1\vert\beta_{6}}{l^2 (k_1 +k_3 )^2} & \hbox{if}\quad k=0, \vspace{0.2cm}\\ 
\frac{\delta^2 k_3{\sqrt{k_3}}\beta_{6} N}{l^2 (k_1 +k_3 )^2 g(0)} & \hbox{if}\quad k=1\,\,\hbox{and}\,\,g_0 = k_3, \vspace{0.2cm}\\ 
\frac{\left\vert(k_3 -k_1 )(k_3 -g_0 )-\delta^2 k_3\right\vert\beta_6}{l^2 (k_1 +k_3 )^2 \vert k_3 -g_0 \vert} & \hbox{if}\quad k=1 \,\hbox{and}\,\,g_0 \ne k_3 .\end{array}
\right.  \label{eq_2_71} 
\end{equation}
Now, using Young's inequality, we find, for $k=1$, 
\begin{equation*}
\begin{array}{lll}
\left\Vert \Phi _{n}\right\Vert _{\mathcal{H}}^{2} &\geq& k_3 \left\Vert w_{n,x} -l\varphi_n\right\Vert^{2} \\
\\
&\geq& k_3 \displaystyle\int_{0}^{1} \left[\frac{N^2}{2}\vert\alpha_{3,n}\vert^{2}\cos^2\,\left(Nx\right)-l^2\vert\alpha_{1,n}\vert ^{2}\cos^2\,\left(Nx\right)\right]\,dx\\ 
\\
&\geq& k_3 \left(\frac{N^2}{4}\vert\alpha_{3,n}\vert^{2} -\frac{l^2}{2}\vert\alpha_{1,n}\vert^{2}\right)\displaystyle\int_{0}^{1} \left[1+\cos\,\left(2Nx\right)\right]\,dx \\
\\
&\geq& k_3 \left(\frac{N^2}{4}\vert\alpha_{3,n}\vert^{2} -\frac{l^2}{2}\vert\alpha_{1,n}\vert^{2}\right),
\end{array}
\end{equation*}
then, by \eqref{eq_2_7100} and \eqref{eq_2_71}, we get \eqref{fn69}. Similarily (with $\sin^2 \,(Nx)$ instead of $\cos^2 \,(Nx)$), \eqref{fn69} is obtained when $k=0$. On the other hand, for $k=1$, we take $\beta_6 =1$ and we find
\begin{equation*}
\left\Vert F_{n}\right\Vert _{\mathcal{H}}^{2} =\left\Vert f_{6,n}\right\Vert^{2} =\beta_{6}^2 \displaystyle\int_{0}^{1} \sin^2 \,(Nx) dx\leq 1 ,
\end{equation*}  
which implies \eqref{Fn}. For $k=0$, we have, using \eqref{fn00} and \eqref{eq_2_67}, 
\begin{equation}
\left\Vert F_{n}\right\Vert_{\mathcal{H}}^{2} =\left\Vert f_{6,n}
\right\Vert^{2}
+\left\Vert f_{7,n}\right\Vert^{2} = \beta_{6}^2 \displaystyle\int_{0}^{1} \cos^2 \,(Nx) dx+\vert\beta_{7,n}\vert^2 \displaystyle\int_{0}^{1} \sin^2 \,(Nx) dx\leq \beta_{6}^2 (1+a_{3,n} ), \label{fn065+}
\end{equation} 
where
\begin{equation*}
a_{3,n} =\frac{\delta^2\left[ l\left(k_1 +k_3\right)a_{2,n} +\frac{1}{l}\right]^2 N^2 \lambda_n^2}{\left[ \left[2k_3 -l\left(k_1 +k_3\right)a_{1,n} \right]N^2 +l^2 \left(k_1 -k_3\right)\right]^2}.
\end{equation*}
According to \eqref{a12n}, we notice that the sequence $(a_{3,n})_n$ is bounded, and then we can choose    
\begin{equation*}
\beta_6 =\frac{1}{\sqrt{1+\sup_{n\in \mathbb{N}}a_{3,n}}}
\end{equation*}
and conclude from \eqref{fn065+} that \eqref{Fn} is valid.
\vskip0,1truecm
Consequently, we deduce from these three cases that the second condition in \eqref{expon} is not satisfied in case $k=0$, and in case $k=1$ if \eqref{k2k3} or \eqref{k23delta} does not hold, hence the exponential stability does not hold for system \eqref{cdt_10}-\eqref{cdt_1}, and for 
system \eqref{syst1}-\eqref{cdt_10} if \eqref{k2k3} or \eqref{k23delta} is not satisfied. This ends the proof of Theorem \ref{Theorem 3.1}.
\end{proof}

\section{Exponential stability of system \eqref{syst1}-\eqref{cdt_10}}

In this section, we show that system \eqref{syst1}-\eqref{cdt_10} is exponentailly stable if \eqref{lpi}-\eqref{k23delta} are satisfied. 
\vskip0,1truecm
\begin{theorem}\label{Theorem 3.11}
Assume that \eqref{l}-\eqref{f+} and \eqref{lpi}-\eqref{k23delta} hold. Then system \eqref{syst_2} in case $k=1$ is exponentially stable.    
\end{theorem}
\vskip0,1truecm
\begin{proof} 
Because the exponential stability is equivalent to \eqref{expon} (see \cite{huan} and \cite{prus}), and because \eqref{lpi} is equivalent to the first condition in \eqref{expon} (see section 3), then it will be enough to prove that, in case $k=1$, \eqref{k2k3}-\eqref{k23delta} imply the second condition in \eqref{expon}. We assume by contradiction that the second condition in \eqref{expon} is false. Then there is sequences $\left( \lambda _{n}\right) _{n}\subset\mathbb{R}$ and $\left( \Phi _{n}\right) _{n}\subset D\left( \mathcal{A}\right)$ such that 
\begin{equation}
\left\Vert \,\Phi _{n}\right\Vert _{\mathcal{H}}\,=\,1,\quad\forall \,n\in \mathbb{N},  \label{eq_4_4}
\end{equation}
\begin{equation}
\lim_{n\rightarrow \infty }\left\vert \lambda _{n}\right\vert =\infty \label{eq_4_5}
\end{equation} 
and 
\begin{equation}
\underset{n\rightarrow \infty }{\lim }\left\Vert \left( i\lambda _{n}I-
\mathcal{A}\right) \Phi _{n}\right\Vert _{\mathcal{H}}=0.  \label{eq_3_18}
\end{equation}
Defining $\Phi_n$ by 
\begin{equation}
\Phi_{n} =\left(\varphi_{n} ,\overset{\sim }{\varphi}_n ,\psi_n ,\overset{\sim }{\psi}_{n} ,w_n, \overset{\sim }{w}_{n} ,\theta_n ,\eta_n \right)^T\in D\left(\mathcal{A}\right). \label{Phin} 
\end{equation}
We will prove that  
\begin{equation}
\left\Vert\,\Phi_{n}\right\Vert_{\mathcal{H}}\,\rightarrow 0,\label{limPhin}
\end{equation}
which is a contradiction with \eqref{eq_4_4}. The limit \eqref{eq_3_18} in case $k=1$ implies the following convergences: 
\begin{equation}
\left\{
\begin{array}{ll}
i\lambda _{n}\varphi _{n}-\overset{\sim }{\varphi }_{n} \rightarrow 0\,\,&\text{in}\,\,H_1,\vspace{0.2cm}\\
i\lambda _{n}\overset{\sim }{\varphi }
_{n}-k_1 \left( \varphi _{n,x}+\psi _{n}+lw_{n}\right) _{x}-lk_{3}\left(
w_{n,x}-l\varphi _{n}\right) +\delta\theta_{n,x} \rightarrow 0\,\,&\text{in}\,\,L_1,\vspace{0.2cm}\\
i\lambda _{n}\psi _{n}-\overset{\sim }{\psi }_{n}
\rightarrow 0\,\,&\text{in}\,\,H_{0},\vspace{0.2cm}\\
i\lambda _{n}\overset{\sim }{\psi }
_{n}-k_2\psi _{n,xx}+k_1\left( \varphi _{n,x}+\psi _{n}+lw_{n}\right)
\rightarrow 0\,\,&\text{in}\,\,L^{2}\left( 0,1\right),\vspace{0.2cm}\\
i\lambda _{n}w_{n}-\overset{\sim }{w}_{n}
\rightarrow 0\,\,&\text{in}\,\,H_{0},\vspace{0.2cm}\\
i\lambda _{n}\overset{\sim }{w}
_{n}-k_{3}\left( w_{n,x}-l\varphi _{n}\right) _{x}+lk_1\left( \varphi
_{n,x}+\psi _{n}+lw_{n}\right) \rightarrow 0\,\,&\text{in}
\,\,L^{2}\left( 0,1\right),\vspace{0.2cm}\\
i\lambda _{n}\theta _{n} -\displaystyle\int_0^{\infty}g\eta_{n,xx}ds+\delta \overset{\sim }{\varphi }_{n,x}\rightarrow 0\,\,&\text{in}\,\,L^{2}\left( 0,1\right) ,\vspace{0.2cm}\\
i\lambda_{n}\eta_{n} -\theta_n +\eta_{n,s}\rightarrow 0\,\,&\text{in}\,\,L_g .
\end{array}
\right.  \label{eq_3_19}
\end{equation}
Taking the inner product of $
i\,\lambda _{n}\,I\,-\,\mathcal{A} \,\,\Phi _{n}$ with $\Phi _{n}$
in $\mathcal{H}$ and using \eqref{g} and \eqref{dissp}, we get
\begin{equation*}
Re\left\langle \left(i\lambda _{n} I-\mathcal{A}\right)\Phi _{n},\Phi _{n}\right\rangle_{\mathcal{H}} =Re\left\langle -\mathcal{A}\Phi _{n},\Phi _{n}\right\rangle_{\mathcal{H}}=-\frac{1}{2}\int_0^{\infty} g^{\prime}\Vert\eta_{n,x}\Vert^2 ds\geq\frac{\mu_1}{2} \Vert\eta_{n}\Vert_{L_g}^2 .
\end{equation*}
So, \eqref{eq_4_4} and \eqref{eq_3_18} imply that
\begin{equation}
\eta _{n}\longrightarrow 0\,\,\text{in}\,\,L_g .  \label{eq_3_20}
\end{equation}
We put
\begin{equation}
{\hat\theta}_n (x)=\int_0^x\int_0^y \theta_n (\tau)d\tau dy-\left(\int_0^1\int_0^y \theta_n (\tau)d\tau dy\right)x. \label{hatthetan}
\end{equation} 
We see that 
\begin{equation*}
{\hat\theta}_{n} (0)={\hat\theta}_{n} (1)=0\quad\hbox{and}\quad{\hat\theta}_{n,xx} =\theta_n .
\end{equation*} 
Moreover, using Young's and H\"older's inequalities, we get
\begin{equation*}
\Vert {\hat\theta}_{n}\Vert_{L_g}^2 =g_0 \Vert {\hat\theta}_{n,x}\Vert^2
\leq 4g_0 \Vert \theta_{n}\Vert^2,
\end{equation*}
then ${\hat\theta}_{n}\in L_g$ and $\left(\Vert {\hat\theta}_{n}\Vert_{L_g}\right)_n$ is bounded, since $\theta_{n} \in L^2 (0,1)$ and \eqref{eq_4_4}. 
Therefore, taking the inner product of \eqref{eq_3_19}$_8$ with 
${\hat\theta}_{n}$ in $L_g$ and using \eqref{eq_4_4}, we entail 
\begin{equation*}
i\lambda _{n}\int_0^{\infty} g\left\langle \eta _{n,x},{\hat\theta}_{n,x}\right\rangle ds-g_0 \left\langle \theta_{n,x},{\hat\theta}_{n,x}\right\rangle+\int_0^{\infty} g\left\langle \eta _{n,sx},{\hat\theta}_{n,x}\right\rangle ds \longrightarrow 0,
\end{equation*}
then, integrating with respect to $x$ and $s$, using \eqref{gexp}$_1$ and noticing that $\eta_n (s=0)=0$ (definition of $D(\mathcal{A})$), we find 
\begin{equation}
-i\lambda _{n} \int_0^{\infty}g\left\langle \eta_{n},\theta_{n}\right \rangle ds+g_0 \Vert\theta_{n}\Vert^2 +\int_0^{\infty}g^{\prime}\left\langle \eta_{n},\theta_{n}\right \rangle ds\longrightarrow 0. \label{HolderPoincare0}
\end{equation}
Moreover, applying Cauchy-Schwartz, H\"older's and Poincar\'e's inequalities, and using \eqref{g}, we get
\begin{equation*}
\left\vert\int_0^{\infty}g^{\prime}\left\langle \eta_{n},\theta_{n}\right \rangle ds\right\vert\leq \mu_2\Vert\theta_{n}\Vert \int_0^{\infty}{\sqrt{g}}{\sqrt{g}}\Vert\eta_{n}\Vert ds\leq \mu_2{\sqrt{g_0}}\Vert\theta_{n}\Vert
\left(\int_0^{\infty}g\Vert\eta_{n}\Vert^2 ds\right)^{\frac{1}{2}} \leq \mu_2 c_0 {\sqrt{g_0}}\Vert\theta_{n}\Vert \Vert\eta_{n}\Vert_{L_g} ,
\end{equation*}
where $c_0$ is the Poincar\'e's constant, so, using \eqref{eq_4_4} and \eqref{eq_3_20}, we observe that
\begin{equation}
\int_0^{\infty}g^{\prime}\left\langle \eta_{n},\theta_{n}\right \rangle ds\longrightarrow 0.\label{HolderPoincare}
\end{equation}
Then, combining \eqref{HolderPoincare0} and \eqref{HolderPoincare},  
we conlude that 
\begin{equation}
g_0\Vert\theta _{n}\Vert^2 -i\lambda _{n}\int_0^{\infty} g\left\langle \eta _{n},\theta_{n}\right\rangle ds\longrightarrow 0.\label{eq_4_100}
\end{equation}
On the other hand, similarily, we see that 
\begin{equation*}
\left\Vert\int_0^{\infty}g\eta_{n} ds\right\Vert\leq \int_0^1\left(\int_0^{\infty}{\sqrt{g}}{\sqrt{g}}\vert\eta_{n}\vert ds\right)^2 dx\leq g_0\int_0^{\infty}g\Vert\eta_{n}\Vert^2 ds\leq c_0^2  g_0\Vert\eta_{n}\Vert_{L_g}^2 ,
\end{equation*}
thus 
\begin{equation}
\int_0^{\infty}g\eta_{n} ds\in L^2 (0,1). \label{intgetan}
\end{equation}
Taking the inner product of \eqref{eq_3_19}$_7$ with $\int_0^{\infty} g\eta_n ds$ in $L^2 (0,1)$, integrating by parts and using \eqref{eq_4_4} and the boundary conditions, we entail 
\begin{equation}
-i\lambda _{n}\int_0^{\infty} g\left\langle \eta _{n},\theta_{n}\right\rangle ds+\left\langle \int_0^{\infty} g\eta _{n,x} ds,\int_0^{\infty} g\eta _{n,x} ds\right\rangle -\delta\int_0^{\infty} g\left\langle \eta _{n,x},\overset{\sim }{\varphi}_{n}\right\rangle ds\longrightarrow 0,\label{intgetands}
\end{equation}
then, using \eqref{eq_4_4} and \eqref{eq_3_20}, it is clear that the last two terms in \eqref{intgetands} converge to zero, so we get
\begin{equation}
\lambda _{n}\int_0^{\infty} g\left\langle \eta _{n},\theta_{n}\right\rangle ds \longrightarrow 0, \label{6eq_3_24}
\end{equation}
so, by combinig \eqref{eq_4_100} and \eqref{6eq_3_24}, we find 
\begin{equation}
\theta_{n}\rightarrow 0\,\,\text{in}\,\,L^{2} \left( 0,1\right). \label{eq_3_24}
\end{equation} 
Multiplying \eqref{eq_3_19}$_1$, \eqref{eq_3_19}$_3$ and \eqref{eq_3_19}$_5$ by $\lambda _{n}^{-1}$, and using \eqref{eq_4_4} and \eqref{eq_4_5}, we obtain
\begin{equation}
\varphi _{n}\longrightarrow 0\text{ in}\,\,L_1\quad\hbox{and}\quad
\psi _{n} ,\,w_{n}\longrightarrow 0\text{ in}\,\,L^2 (0,1).\label{eq_3_21}
\end{equation}
Multiplying \eqref{eq_3_19}$_2$ by $\lambda _{n}^{-1}$ and using \eqref{eq_4_5}, we find
\begin{equation*}
i\overset{\sim }{\varphi }_{n}-\frac{k_1}{\lambda _{n}}\left( \varphi
_{n,x}+\psi _{n}+lw_{n}\right) _{x}-\frac{lk_{3}}{\lambda _{n}}\left(
w_{n,x}-l\varphi _{n}\right) + \frac{\delta}{\lambda _{n}}\theta_{n,x}
\longrightarrow 0 \text{ in}\,\,L^{2}\left( 0,1\right).
\end{equation*}
Using \eqref{eq_4_4} and \eqref{eq_4_5}, we conclude that
\begin{equation}
\left( \lambda _{n}^{-1}(k_1 \varphi _{n,xx}-\delta\theta_{n,x})\right) _{n}\,\,\text{is bounded in}\,\,L^{2}\left( 0,1\right). \label{eq_3_25}
\end{equation} 
Taking the inner product of \eqref{eq_3_19}$_7$
with $\lambda _{n}^{-1}(k_1 \varphi _{n,x}-\delta\theta_{n})$ in 
$L^{2}\left( 0,1\right)$, integrating by parts and using the boundary conditions, \eqref{eq_4_4} and \eqref{eq_4_5}, we entail
\begin{equation*}
i\left\langle \theta _{n},k_1 \varphi _{n,x}-\delta\theta_{n}\right\rangle +\int_0^{\infty}g\left\langle \eta_{n,x},
\lambda_{n}^{-1} (k_1 \varphi _{n,xx}-\delta\theta_{n,x})\right\rangle
\end{equation*}
\begin{equation*}
-\delta \left\langle i\lambda _{n}\varphi _{n,x}-\overset{\sim }{\varphi }_{n,x} ,\lambda_{n}^{-1}(k_1 \varphi _{n,x}-\delta\theta_{n}) \right\rangle -i\delta^2 \left\langle \varphi _{n,x} ,\theta_{n} \right\rangle +i\delta k_1\left\Vert \varphi _{n,x}\right\Vert ^{2}\longrightarrow 0,
\end{equation*}
Combining \eqref{eq_4_4}, \eqref{eq_4_5}, \eqref{eq_3_19}$_1$, \eqref{eq_3_20}, \eqref{eq_3_24}, \eqref{eq_3_25} and the above limit, it follows that
\begin{equation}
\varphi _{n,x}\longrightarrow 0\,\,\text{in}\,\,L^{2}\left(0,1\right). \label{eq_3_27}
\end{equation} 
Moreover, by \eqref{eq_3_19}$_1$ and \eqref{eq_3_27}, we see that
\begin{equation}
\lambda _{n}^{-1}\overset{\sim }{\varphi }_{n,x}\rightarrow 0\text{ in}\,\,L^{2}\left(0,1\right), \label{eq_3_28}
\end{equation}
and so, thanks to the definition of $L_1$ and Poincar\'e's inequality, \eqref{eq_3_28} leads to
\begin{equation}
\lambda _{n}^{-1}\overset{\sim }{\varphi }_{n}\longrightarrow 0\,\,\text{in}\,\,L_1 .\label{eq_3_29}
\end{equation}
Taking the inner product of \eqref{eq_3_19}$_2$
with $\lambda _{n}^{-1}\overset{\sim }{\varphi}_{n}$ in $L^{2}\left( 0,1\right)$, integrating by parts and using \eqref{eq_3_29} and the boundary conditions, we get
\begin{equation*}
i\left\Vert \overset{\sim }{\varphi }_{n}\right\Vert^{2}
+\left\langle \lambda _{n}^{-1}\left( k_1\varphi _{n,xx}-\delta\theta_{n,x}\right) ,i\lambda_n \varphi_n -\overset{\sim }{\varphi }_{n}\right\rangle 
\end{equation*}
\begin{equation*}
-i\left\langle k_1\varphi _{n,x}-\delta\theta_{n} ,\varphi_{n,x} \right\rangle-\left\langle k_1 (\psi _{n,x}+lw_{n,x} )+lk_3 (w_{n,x}-l\varphi_n ) ,\lambda _{n}^{-1}\overset{\sim }{\varphi }_{n}\right\rangle\rightarrow 0.
\end{equation*}
So, using \eqref{eq_4_4}, \eqref{eq_3_19}$_1$, \eqref{eq_3_25},  \eqref{eq_3_27} and \eqref{eq_3_29}, we deduce that
\begin{equation}
\overset{\sim }{\varphi }_{n}\longrightarrow 0\,\,\text{in}\,\,L_1 ,  \label{eq_3_30}
\end{equation}
and by \eqref{eq_3_19}$_1$, we find
\begin{equation}
\lambda _{n}\varphi _{n}\longrightarrow 0\;\text{in}\,\,L_1 . \label{eq_3_31}
\end{equation}  
Taking the inner product of \eqref{eq_3_19}$_{4}$ with $w_n$ in $L^{2}\left( 0,1\right)$, integration by parts and using the boundary conditions, we get 
\begin{equation*}
-\left\langle \overset{\sim }{\psi}_{n},i\lambda_{n} w_{n}-\overset{\sim }{w}_{n} \right\rangle -\left\langle \overset{\sim }{\psi}_{n},\overset{\sim }{w}_{n}\right\rangle +k_2\left\langle \psi _{n,x},w_{n,x}\right\rangle +k_1 \left\langle\varphi_{n,x}+\psi_{n}+l\,w_{n} ,w_{n}\right\rangle\rightarrow 0,
\end{equation*}
then, using \eqref{eq_4_4}, \eqref{eq_3_19}$_{5}$ and  
\eqref{eq_3_21}, we deduce that 
\begin{equation}
k_2\left\langle \psi _{n,x},w_{n,x}-l\varphi _{n}\right\rangle -\left\langle \overset{\sim}{\psi}_{n},\overset{\sim}{w}_{n}\right\rangle \rightarrow 0.  \label{eq_3_32}
\end{equation}
Taking the inner product of \eqref{eq_3_19}$_{2}$ with 
$w_{n,x}-l\varphi _{n}$ in $L^{2}\left( 0,1\right)$, we obtain 
\begin{equation*}
-\frac{1}{k_1}\left\langle \overset{\sim}{\varphi}_{n},i\lambda
_{n}w_{n,x}-\overset{\sim}{w}_{n,x}\right\rangle -\frac{1}{k_1}\left\langle
\overset{\sim}{\varphi}_{n},\overset{\sim}{w}_{n,x}\right\rangle +\frac{l}{k_1}\left\langle \overset{\sim}{\varphi}_{n},i\lambda_{n}\varphi _{n}\right\rangle 
\end{equation*}
\begin{equation*}
-\left\langle \left( \varphi _{n,x}+\psi _{n}+l\,w_{n}\right) _{x},
w_{n,x}-l\varphi _{n} \right\rangle -\frac{lk_{3}}{k_1}\left\Vert
w_{n,x}-l\varphi _{n} \right\Vert ^{2}+\frac{\delta }{k_1}
\left\langle \theta _{n,x},w_{n,x}-l\varphi _{n}\right\rangle
\rightarrow 0.
\end{equation*}
By using \eqref{eq_4_4}, \eqref{eq_3_19}$_{5}$ and \eqref{eq_3_31}, we get
\begin{equation}
-\frac{1}{k_1}\left\langle \overset{\sim}{\varphi}_{n},\overset{\sim}{w}_{n,x}\right\rangle -\left\langle \left( \varphi _{n,x}+\psi _{n}+l\,w_{n}\right)_{x},w_{n,x}-l\varphi _{n} \right\rangle \label{eq33+}
\end{equation}
\begin{equation*}  
-\frac{lk_{3}}{k_1}\left\Vert w_{n,x}-l\varphi _{n}
\right\Vert ^{2}+\frac{\delta}{k_1}\left\langle \theta _{n,x},
w_{n,x}-l\varphi _{n} \right\rangle \rightarrow 0. 
\end{equation*}
Taking the inner product of \eqref{eq_3_19}$_{6}$ with $\varphi
_{n,x}+\psi _{n}+lw_{n}$ in $L^{2}\left( 0,1\right)$, integration by parts and using the boundary conditions, we find 
\begin{equation*}
-\left\langle \overset{\sim}{w}_{n},i\lambda _{n}\varphi _{n,x}\right\rangle
-\left\langle \overset{\sim}{w}_{n},i\lambda _{n}\psi _{n}-\overset{\sim}{\psi}_{n}\right\rangle -\left\langle \overset{\sim}{w}_{n},\overset{\sim}{\psi}_{n}\right\rangle -l\left\langle \overset{\sim}{w}_{n},i\lambda_{n}\,w_{n}-\overset{\sim}{w}_{n}\right\rangle 
\end{equation*}
\begin{equation*}
-l\left\Vert \overset{\sim}{w}_{n}\right\Vert^{2}+k_{3}\left\langle 
w_{n,x}-l\varphi _{n} ,\left( \varphi _{n,x}+\psi
_{n}+l\,w_{n}\right) _{x}\right\rangle +lk_1\left\Vert \varphi
_{n,x}+\psi _{n}+l\,w_{n} \right\Vert ^{2}\rightarrow 0,  
\end{equation*}
therefore, using \eqref{eq_4_4}, \eqref{eq_3_19}$_{3}$, \eqref{eq_3_19}$_{5}$, \eqref{eq_3_21}, \eqref{eq_3_27} and \eqref{eq_3_32}, we deduce that 
\begin{equation}
-\frac{1}{k_{3}}\left\langle \overset{\sim}{w}_{n},i\lambda _{n}\varphi
_{n,x}\right\rangle -\frac{k_2}{k_{3}}\left\langle
w_{n,x}-l\varphi _{n},\psi _{n,x}\right\rangle  -\frac{l}{k_{3}}\left\Vert \overset{\sim}{w}_{n}\right\Vert ^{2}+
\left\langle w_{n,x}-l\varphi _{n} ,\left( \varphi _{n,x}+\psi
_{n}+l\,w_{n}\right) _{x}\right\rangle \rightarrow 0,  \label{eq34}
\end{equation}
combining \eqref{eq33+} and \eqref{eq34}, we find 
\begin{equation*}
-\frac{1}{k_{3}}\left\langle \overset{\sim}{w}_{n},i\lambda _{n}\varphi
_{n,x}\right\rangle -\frac{k_2}{k_{3}}\left\langle
w_{n,x}-l\varphi _{n},\psi _{n,x}\right\rangle -\frac{l}{k_{3}}
\left\Vert \overset{\sim}{w}_{n}\right\Vert^{2} 
\end{equation*}
\begin{equation*}
-\frac{1}{k_1}\left\langle \overset{\sim}{w}_{n,x} ,\overset{\sim}{\varphi}_{n}\right\rangle -\frac{lk_{3}}{k_1}\left\Vert w_{n,x}-l\varphi _{n}\right\Vert^{2}+\frac{\delta}{k_1}\left\langle w_{n,x}-l\varphi _{n} ,\theta _{n,x}\right\rangle \rightarrow 0,
\end{equation*}
then, with integration by parts and using the bounday conditions, we obtain 
\begin{equation*}
-\frac{1}{k_{3}}\left\langle \overset{\sim}{w}_{n},i\lambda _{n}\varphi
_{n,x}\right\rangle -\frac{k_2}{k_{3}}\left\langle
w_{n,x}-l\varphi _{n},\psi _{n,x}\right\rangle -\frac{l}{k_{3}}
\left\Vert \overset{\sim}{w}_{n}\right\Vert ^{2} -\frac{1}{k_1}\left\langle \overset{\sim}{w}_{n} , i\lambda _{n}\varphi_{n,x}-\overset{\sim}{\varphi}_{n,x} \right \rangle
\end{equation*}
\begin{equation*}
+\frac{1}{k_1}\left\langle \overset{\sim}{w}_{n}, i\lambda _{n}\varphi _{n,x}\right\rangle -\frac{lk_{3}}{k_1}\left\Vert w_{n,x}-l\varphi _{n}
\right\Vert ^{2}+\frac{\delta }{k_1}\left\langle
w_{n,x}-l\varphi _{n} ,\theta _{n,x}\right\rangle \rightarrow 0,
\end{equation*}
using \eqref{eq_4_4} and \eqref{eq_3_19}$_{1}$, we arrive at
\begin{equation}
\frac{1}{k_{3}}\left( \frac{k_{3}}{k_1}-1\right) 
\left\langle \overset{\sim}{w}_{n},i\lambda _{n}\varphi _{n,x}\right\rangle -\frac{k_2}{k_{3}}\left\langle w_{n,x}-l\varphi _{n},\psi_{n,x}\right\rangle -\frac{l}{k_{3}}\left\Vert \overset{\sim}{w}_{n}\right\Vert ^{2} \label{eq_3.35} 
\end{equation}
\begin{equation*}
-\frac{lk_{3}}{k_1}\left\Vert w_{n,x}-l\varphi _{n}
\right\Vert ^{2}+\frac{\delta }{k_1}\left\langle 
w_{n,x}-l\varphi _{n} ,\theta _{n,x}\right\rangle \rightarrow 0. 
\end{equation*}
From \eqref{eq_4_4}, \eqref{eq_3_19}$_{3}$ and \eqref{eq_3_19}$_{5}$, we observe that 
\begin{equation}
\left( \Vert \lambda _{n}\psi _{n}\Vert \right)_{n\in \mathbb{N}} \,\,\hbox{and} \,\, \left(\Vert \lambda _{n} w_{n}\Vert \right)_{n\in 
\mathbb{N}}\,\,\hbox{are bounded.}
\label{psiwlambdanbound}
\end{equation}
We have, by integrating by parts and using the boundary conditions, 
\begin{eqnarray*}
\left\langle \lambda _{n}^{2}\psi _{n}+i\lambda _{n} \overset{\sim}{\psi}_{n},i\theta _{n}\right\rangle &=&-i\left\langle 
i\lambda _{n}\psi _{n}-\overset{\sim}{\psi}_{n} ,i\lambda _{n}\theta _{n}\right\rangle \\
&=&-i\left\langle i\lambda _{n}\psi_{n}-\overset{\sim}{\psi}_{n} ,i\lambda _{n}\theta _{n}-\int_0^{\infty}g\eta_{n,xx} ds+\delta
\overset{\sim}{\varphi}_{n,x} \right\rangle \\
&&-i\delta \left\langle i\lambda _{n}\psi _{n,x}-\overset{\sim}{\psi}_{n,x},\overset{\sim}{\varphi}_{n}\right\rangle +i\int_0^{\infty}g\left\langle i\lambda _{n}\psi _{n,x}-\overset{\sim}{\psi}_{n,x},\eta_{n,x}\right\rangle ds,
\end{eqnarray*}
by using \eqref{eq_4_4}, \eqref{eq_3_19}$_{3}$ and \eqref{eq_3_19}$_{7}$, we deduce that 
\begin{equation}
\left\langle \lambda _{n}^{2}\psi _{n}+i\lambda _{n} \overset{\sim}{\psi}_{n},i\theta _{n}\right\rangle \rightarrow 0.  \label{eq_3_36}
\end{equation}
Also, we have 
\begin{equation}
\left\langle \lambda _{n}\psi _{n},\overset{\sim}{\varphi}_{n,x}\right\rangle =-\left\langle \lambda _{n}\psi _{n,x},\overset{\sim}{\varphi}_{n}\right\rangle =-\left\langle i\lambda _{n}\psi _{n,x}-\overset{\sim}{\psi}_{n,x},i\overset{\sim}{\varphi}_{n}\right\rangle +\left\langle \overset{\sim}{\psi}_{n},i\overset{\sim}{\varphi}_{n,x}\right\rangle .
\label{eq37}
\end{equation}
Using again integration by parts and the boundary conditions, we have
\begin{equation*}
\lambda _{n}\int_0^{\infty}g\left\langle \psi _{n,x},\eta_{n,x}\right\rangle ds =-\lambda_{n}\int_0^{\infty}g\left\langle \psi _{n},\eta_{n,xx}\right\rangle ds =-\lambda_{n}\left\langle \psi _{n},\int_0^{\infty}g\eta_{n,xx} ds\right\rangle
\end{equation*}
\begin{eqnarray*}
&=&\lambda _{n}\left\langle \psi _{n},i\lambda _{n}\theta
_{n}-\int_0^{\infty}g\eta_{n,xx}ds+\delta \overset{\sim}{\varphi}_{n,x}\right\rangle -\lambda _{n}\left\langle\psi _{n},i\lambda _{n}\theta _{n}\right\rangle -\delta \lambda_{n}\left\langle \psi _{n},\overset{\sim}{\varphi}_{n,x}\right\rangle \\
&=& \left\langle \lambda _{n}\psi _{n},i\lambda _{n}\theta
_{n}-\int_0^{\infty}g\eta_{n,xx}ds+\delta \overset{\sim}{\varphi}_{n,x} \right\rangle -\left\langle \lambda _{n}^{2}\psi _{n}+i\lambda _{n}\overset{\sim}{\psi}_{n} ,i\theta _{n}\right\rangle \\
&&+\left\langle i\lambda _{n} \overset{\sim}{\psi}_{n}-k_2\psi _{n,xx}+k_1 \left( \varphi _{n,x}+\psi _{n}+lw_{n}\right)
,i\theta _{n}\right\rangle \\
&&+k_2\left\langle \psi _{n,xx},i\theta_{n}\right\rangle -k_1\left\langle \varphi_{n,x}+\psi _{n}+lw_{n} ,i\theta _{n}\right\rangle -\delta
\left\langle\lambda_{n}\psi _{n},\overset{\sim}{\varphi}_{n,x}\right\rangle ,
\end{eqnarray*}
then, by \eqref{eq37}, integration by parts and using the boundary conditions, we obtain 
\begin{equation}
\lambda _{n}\int_0^{\infty}g\left\langle \psi _{n,x},\eta_{n,x}\right\rangle ds = \left\langle \lambda _{n}\psi _{n},i\lambda _{n}\theta
_{n}-\int_0^{\infty}g\eta_{n,xx}ds+\delta \overset{\sim}{\varphi}_{n,x} \right\rangle -\left\langle \lambda _{n}^{2}\psi _{n}+i\lambda _{n}\overset{\sim}{\psi}_{n} ,i\theta _{n}\right\rangle \label{eq_3.38}
\end{equation}
\begin{eqnarray}
&+&\left\langle i\lambda _{n} \overset{\sim}{\psi}_{n}-k_2\psi _{n,xx}+k_1 \left( \varphi _{n,x}+\psi _{n}+lw_{n}\right),i\theta _{n} \right\rangle -k_2\left\langle \psi _{n,x},i\theta_{n,x}\right\rangle \notag\\
&-& k_1\left\langle \varphi_{n,x}+\psi _{n}+lw_{n} ,i\theta _{n}\right\rangle + \delta \left\langle i\lambda _{n}\psi _{n,x}-\overset{\sim}{\psi}_{n,x},i\overset{\sim}{\varphi}_{n}\right\rangle -\delta\left\langle \overset{\sim}{\psi}_{n},i\overset{\sim}{\varphi}_{n,x}\right\rangle , \notag
\end{eqnarray}
using \eqref{eq_4_4}, \eqref{eq_3_19}$_{3}$, \eqref{eq_3_19}$_{4}$, \eqref{eq_3_19}$_{7}$, \eqref{eq_3_24}, \eqref{eq_3_27}, \eqref{psiwlambdanbound} and \eqref{eq_3_36}, we deduce from \eqref{eq_3.38} that 
\begin{equation}
\lambda _{n}\int_0^{\infty}g\left\langle \psi _{n,x},\eta_{n,x}\right\rangle ds +k_2\left\langle \psi _{n,x},i\theta _{n,x}\right\rangle +\delta
\left\langle \overset{\sim}{\psi}_{n},i\overset{\sim}{\varphi}_{n,x}\right\rangle \rightarrow 0.  \label{eq_3_39}
\end{equation}
Also, by integrating with respect to $s$ and using \eqref{gexp}$_1$ and 
$\eta_n (s=0)=0$, we have  
\begin{eqnarray*}
\lambda _{n}\int_0^{\infty}g\left\langle \psi_{n,x},\eta_{n,x}\right\rangle ds &=& i\int_0^{\infty}g\left\langle \psi_{n,x},i\lambda_{n}\eta_{n,x}\right\rangle ds \\ 
&=& i\int_0^{\infty}g\left\langle \psi _{n,x},i\lambda _{n}\eta_{n,x} -\theta_{n,x} +\eta_{n,xs} \right\rangle ds \\
&&+i\int_0^{\infty}g^{\prime}\left\langle \psi _{n,x},\eta_{n,x}\right\rangle ds-g_0 \left\langle \psi _{n,x}, i\theta_{n,x}\right\rangle,   
\end{eqnarray*}
therefore, by using \eqref{g}, \eqref{eq_4_4}, \eqref{eq_3_19}$_{8}$, 
\eqref{eq_3_20}, \eqref{eq_3_39} and the above identity, we obtain 
\begin{equation*}
\left( k_2 -g_0 \right) \left\langle \psi_{n,x},\theta _{n,x}\right\rangle +\delta \left\langle\overset{\sim}{\psi}_{n},\overset{\sim}{\varphi}_{n,x}\right\rangle \rightarrow 0,
\end{equation*}
and so
\begin{equation*}
\left( k_2 -g_0 \right) \left\langle \psi_{n,x},\theta _{n,x}\right\rangle -\delta \left\langle\overset{\sim}{\psi}_{n},i\lambda_n \varphi_{n,x} -\overset{\sim}{\varphi}_{n,x}\right\rangle + 
\delta \left\langle\overset{\sim}{\psi}_{n},i\lambda_n \varphi_{n,x}\right\rangle\rightarrow 0,
\end{equation*}
and moreover, by \eqref{eq_4_4} and \eqref{eq_3_19}$_{1}$, we find 
\begin{equation}
\left( k_2 -g_0\right) \left\langle \psi
_{n,x},\theta _{n,x}\right\rangle +\delta \left\langle\overset{\sim}{\psi}_{n},i\lambda_n \varphi_{n,x}\right\rangle \rightarrow 0.  \label{eq_3_43}
\end{equation}
Taking the inner product of \eqref{eq_3_19}$_{4}$ with 
$\varphi _{n,x}+\psi _{n}+l\,w_{n}$ in $L^2 (0,1)$, integration by parts and using the boundary conditions, \eqref{eq_3_21} and \eqref{eq_3_27}, we get 
\begin{equation*}
-\left\langle \overset{\sim}{\psi}_{n},i\lambda _{n}\varphi _{n,x}\right\rangle -\left\langle \overset{\sim}{\psi}_{n}, i\lambda _{n}\psi _{n}-\overset{\sim}{\psi}_{n}\right\rangle -\left\Vert \overset{\sim}{\psi}_{n}\right\Vert ^{2} 
\end{equation*}
\begin{equation*}
-l\left\langle \overset{\sim}{\psi}_{n}, i\lambda _{n}\,w_{n}-\overset{\sim}{w}_{n}\right\rangle -l\left\langle \overset{\sim}{\psi}_{n},\overset{\sim}{w}_{n}\right\rangle +k_2\left\langle \psi _{n,x},\left( \varphi _{n,x}+\psi _{n}+lw_{n}\right)_{x}\right\rangle \rightarrow 0,
\end{equation*}
using \eqref{eq_4_4}, \eqref{eq_3_19}$_{3}$ and 
\eqref{eq_3_19}$_{5}$, we deduce that 
\begin{equation*}
-\left\langle \overset{\sim}{\psi}_{n},i\lambda _{n}\varphi _{n,x}\right\rangle -\left\Vert \overset{\sim}{\psi}_{n}\right\Vert ^{2}-l\left\langle
\overset{\sim}{\psi}_{n},\overset{\sim}{w}_{n}\right\rangle +\frac{k_2}{k_1}\left\langle \psi _{n,x},i\lambda _{n} \overset{\sim}{\varphi}_{n}-lk_{3}\left( w_{n,x}-l\varphi _{n}\right) +\delta \theta_{n,x}\right\rangle 
\end{equation*}
\begin{equation*}
-\frac{k_2}{k_1}\left\langle \psi _{n,x},i\lambda _{n} \overset{\sim}{\varphi}_{n}-k_1\left( \varphi _{n,x}+\psi _{n}+lw_{n}\right) _{x}-lk_{3}\left(w_{n,x}-l\varphi_{n}\right) +\delta \theta_{n,x} \right\rangle \rightarrow 0,
\end{equation*}
using \eqref{eq_4_4} and \eqref{eq_3_19}$_{2}$, we have 
\begin{equation*}
-\left\langle \overset{\sim}{\psi}_{n},i\lambda _{n}\varphi _{n,x}\right\rangle -\left\Vert \overset{\sim}{\psi}_{n}\right\Vert ^{2}-l\left\langle
\overset{\sim}{\psi}_{n},\overset{\sim}{w}_{n}\right\rangle -\frac{k_2}{k_1}\left\langle i\lambda_{n}\psi _{n,x}-\overset{\sim}{\psi}_{n,x} ,\overset{\sim}{\varphi}_{n}\right\rangle 
\end{equation*}
\begin{equation*}
-\frac{k_2}{k_1}\left\langle \overset{\sim}{\psi}_{n,x},\overset{\sim}{\varphi}_{n}\right\rangle -\frac{lk_2 k_{3}
}{k_1}\left\langle \psi _{n,x}, w_{n,x}-l\varphi _{n}\right\rangle +\frac{k_2\delta }{k_1}\left\langle \psi _{n,x},\theta
_{n,x}\right\rangle \rightarrow 0. \label{eq_3_44+}
\end{equation*}
As, by integrating by parts and using the boundary conditions, 
\begin{equation*}
\left\langle \overset{\sim}{\psi}_{n,x},\overset{\sim}{\varphi}_{n}\right\rangle =-\left\langle \overset{\sim}{\psi}_{n},\overset{\sim}{\varphi}_{n,x}\right\rangle =\left\langle \overset{\sim}{\psi}_{n},i\lambda
_{n}\varphi _{n,x}-\overset{\sim}{\varphi}_{n,x}\right\rangle -\left\langle \overset{\sim}{\psi}_{n},i\lambda_{n}\varphi _{n,x}\right\rangle ,
\end{equation*}
and with \eqref{eq_4_4}, \eqref{eq_3_19}$_{1}$, \eqref{eq_3_19}$_{3}$ and 
\eqref{eq_3_44+}, we see that 
\begin{equation*}
\left( \frac{k_2}{k_1}-1\right) \left\langle \overset{\sim}{\psi}_{n},i\lambda_{n}\varphi _{n,x}\right\rangle -\left\Vert \overset{\sim}{\psi}_{n}\right\Vert^{2}-l\left\langle \overset{\sim}{\psi}_{n},\overset{\sim}{w}_{n}\right\rangle -\frac{lk_2 k_{3}}{k_1}\left\langle \psi _{n,x},w_{n,x}-l\varphi_{n}\right\rangle +\frac{k_2\delta }{k_1}\left\langle \psi _{n,x},\theta_{n,x}\right\rangle \rightarrow 0,
\end{equation*}
thus, combining this limit with \eqref{eq_3_32} and \eqref{eq_3_43}, we obtain 
\begin{equation}
\frac{1}{\delta k_1}\left[ k_2\delta ^{2}-\left( k_2 -k_1
\right) \left( k_2 -g_0\right) \right]\left\langle \psi _{n,x},\theta _{n,x}\right\rangle -\left\Vert \overset{\sim}{\psi}_{n}\right\Vert ^{2} -lk_2 \left( 1+\frac{k_{3}}{k_1}\right) \left\langle \psi
_{n,x},w_{n,x}-l\varphi _{n} \right\rangle \rightarrow 0.
\label{eq_3_44}
\end{equation}
We have, by integrating by parts (with respect to $x$ and $s$) and using the boundary conditions, 
\begin{eqnarray}
\int_0^{\infty}g\left\langle \overset{\sim}{w}_{n},\eta_{n,xx}\right\rangle ds &=& \int_0^{\infty}g \left\langle i\lambda _{n}w_{n,x} -\overset{\sim}{w}_{n,x} ,\eta_{n,x}\right\rangle ds -\int_0^{\infty}g\left\langle i\lambda _{n}w_{n,x},\eta_{n,x}\right\rangle ds \notag \\
&=&\int_0^{\infty}g\left\langle i\lambda _{n} w_{n,x}-\overset{\sim}{w}_{n,x} 
,\eta_{n,x}\right\rangle  ds+\int_0^{\infty}g\left\langle w_{n,x},i\lambda
_{n}\eta_{n,x} -\theta _{n,x}+\eta_{n,xs} \right\rangle ds
\label{eq_3.45} \\
&&+\int_0^{\infty}g^{\prime}\left\langle w_{n,x},\eta_{n,x}\right\rangle ds +g_0 \left\langle w_{n,x},\,\theta _{n,x}\right\rangle .  \notag
\end{eqnarray}
Also, we see that
\begin{eqnarray*}
\left\langle i\lambda _{n} \overset{\sim}{w}_{n},\theta _{n}\right\rangle &=&-\left\langle \overset{\sim}{w}_{n},i\lambda _{n}\,\theta _{n}\right\rangle \\
&=&-\left\langle \overset{\sim}{w}_{n},i\lambda _{n} \theta _{n} -\int_0^{\infty}g\eta_{n,xx} ds+\delta \overset{\sim}{\varphi}_{n,x}\right\rangle -\int_0^{\infty}\left\langle \overset{\sim}{w}_{n},g\eta_{n,xx}\right\rangle ds \\
&&-\delta\left\langle \overset{\sim}{w}_{n},i\lambda
_{n}\varphi _{n,x}-\overset{\sim}{\varphi}_{n,x} \right\rangle +\delta\left\langle \overset{\sim}{w}_{n},i\lambda _{n}\varphi _{n,x}\right\rangle ,
\end{eqnarray*}
by using \eqref{eq_3.45}, we obtain 
\begin{eqnarray}
\left\langle i\lambda _{n} \overset{\sim}{w}_{n},\,\theta _{n}\right\rangle &=& -\left\langle \overset{\sim}{w}_{n},i\lambda _{n} \theta _{n} -\int_0^{\infty}g\eta_{n,xx} ds+\delta \overset{\sim}{\varphi}_{n,x}\right\rangle -\delta\left\langle \overset{\sim}{w}_{n},i\lambda
_{n}\varphi _{n,x}-\overset{\sim}{\varphi}_{n,x} \right\rangle \notag \\
&& -\int_0^{\infty} g\left\langle i\lambda_n w_{n,x} -\overset{\sim}{w}_{n,x},\eta_{n,x}\right\rangle ds -\int_0^{\infty} g\left\langle w_{n,x},i\lambda_n\eta_{n,x} -\theta_{n,x}+\eta_{n,xs}\right\rangle ds \notag \\
&&+\delta\left\langle \overset{\sim}{w}_{n},i\lambda_n \varphi_{n,x}\right\rangle -\int_0^{\infty}g^{\prime}\left\langle w_{n,x},\eta_{n,x}\right\rangle ds -g_0 \left\langle w_{n,x},\theta_{n,x}\right\rangle . 
\label{eq_3.46}
\end{eqnarray}
Taking the inner product of \eqref{eq_3_19}$_{6}$ with 
$\theta _{n}$ in $L^2 (0,1)$, integration by parts and using the boundary conditions and \eqref{eq_4_4}, we find
\begin{equation*}
\left\langle i\lambda _{n} \overset{\sim}{w}_{n},\theta _{n}\right\rangle
+k_{3}\left\langle w_{n,x}-l\varphi _{n} ,\theta
_{n,x}\right\rangle +lk_1 \left\langle \varphi _{n,x}+\psi
_{n}+lw_{n} ,\theta _{n}\right\rangle \rightarrow 0,
\end{equation*}
then, exploiting \eqref{eq_4_4}, \eqref{eq_3_19}$_{1}$, 
\eqref{eq_3_19}$_{5}$, \eqref{eq_3_19}$_{7}$, \eqref{eq_3_19}$_{8}$, 
\eqref{eq_3_20}, \eqref{eq_3_24} and \eqref{eq_3.46}, we obtain 
\begin{equation*}
-g_0 \left\langle w_{n,x},\theta_{n,x}\right\rangle +\delta\left\langle
\overset{\sim}{w}_{n},i\lambda _{n}\varphi _{n,x}\right\rangle +k_{3}\left
\langle w_{n,x}-l\varphi _{n},\theta _{n,x}\right\rangle\rightarrow 0,
\end{equation*}
and according to \eqref{eq_4_4} and \eqref{eq_3_24}, we observe that  
\begin{equation*}
\left\langle \varphi _{n},\theta _{n,x}\right\rangle =-\left\langle \varphi _{n,x},\theta _{n} \right\rangle\rightarrow 0,
\end{equation*} 
then, by combining the above two limits, we get 
\begin{equation}
\left( k_{3} -g_0\right) \left\langle w_{n,x}-l\varphi _{n} ,\theta _{n,x}\right\rangle +\delta\left\langle \overset{\sim}{w}_{n},i\lambda _{n}\varphi
_{n,x}\right\rangle \rightarrow 0.  \label{eq_3_47}
\end{equation}
By combining \eqref{eq_3.35} and \eqref{eq_3_47}, we observe that 
\begin{equation}
\frac{1}{k_1\delta }\left[ k_3\delta^{2}-\left( k_3 -k_1\right)
\left( k_{3}-g_0\right) \right] \left\langle w_{n,x}-l\varphi _{n} ,\theta_{n,x}\right\rangle \label{eq_3.48}
\end{equation}
\begin{equation*}
-k_2\left\langle w_{n,x}-l\varphi
_{n},\psi _{n,x}\right\rangle -l\left\Vert
\overset{\sim}{w}_{n}\right\Vert ^{2}-\frac{lk_{3}^2}{k_1}\left\Vert w_{n,x}-l\varphi_{n}\right\Vert ^{2}\rightarrow 0.  
\end{equation*}
Taking the inner product in $L^2 (0,1)$ of \eqref{eq_3_19}$_4$ with $w_{n}$, and \eqref{eq_3_19}$_{6}$ with $\psi _{n}$, we get, respectively, 
\begin{equation*}
-\left\langle \overset{\sim}{\psi}_{n},i\lambda _{n}w_{n}-\overset{\sim}{w}_{n} \right\rangle -\left\langle \overset{\sim}{\psi}_{n},\overset{\sim}{w}_{n}\right\rangle -k_2\left\langle\psi _{n,xx},w_{n}\right\rangle +k_1\left\langle \varphi _{n,x}+\psi _{n}+l\,w_{n},w_{n}\right\rangle\rightarrow 0
\end{equation*}
and
\begin{equation*}
-\left\langle \overset{\sim}{w}_{n},i\lambda _{n}\psi _{n}-\overset{\sim}{\psi}_{n} \right\rangle-\left\langle \overset{\sim}{w}_{n},\overset{\sim}{\psi}_{n}\right\rangle -k_{3}\left\langle\left( w_{n,x}-l\varphi _{n}\right) _{x},\psi _{n}\right\rangle +lk_{1}\left\langle \varphi _{n,x}+\psi _{n}+l\,w_{n} ,\psi_{n}\right\rangle \rightarrow 0, 
\end{equation*}
by integration by parts and using the boundary conditions, \eqref{eq_4_4}, 
\eqref{eq_3_19}$_{3}$, \eqref{eq_3_19}$_{5}$ and \eqref{eq_3_21}, we obtain 
\begin{equation*}
-\left\langle \overset{\sim}{\psi}_{n},\overset{\sim}{w}_{n}\right\rangle +k_2\left\langle \psi _{n,x},w_{n,x}\right\rangle \rightarrow 0\quad\hbox{and}\quad -\left\langle \overset{\sim}{\psi}_{n} , \overset{\sim}{w}_{n}\right\rangle +k_{3}\left\langle \psi _{n,x} ,  w_{n,x}-l\varphi _{n} \right\rangle \rightarrow 0, 
\end{equation*}
so, using \eqref{eq_4_4} and \eqref{eq_3_21}, we find  
\begin{equation*}
\left( k_{2}-k_{3}\right) \left\langle \psi _{n,x}, w_{n,x}-l\varphi _{n} \right\rangle\rightarrow 0.
\end{equation*}
Because $k_{2}-k_{3} \ne 0$, then we obtain 
\begin{equation}
\left\langle \psi_{n,x} ,w_{n,x}-l\varphi _{n} \right\rangle \rightarrow 0.  \label{eq_3_49}
\end{equation}
Because $k_2 \delta ^{2}-\left(k_2 -k_1\right) \left(k_2 -g_0\right)=0$, then \eqref{eq_3_44} and \eqref{eq_3_49} imply that 
\begin{equation}
\overset{\sim}{\psi}_{n}\rightarrow 0\,\,\text{in}\,\,L^{2}\left( 0,1\right) .  \label{eq_3_50}
\end{equation}
By \eqref{eq_3_19}$_{3}$\ and \eqref{eq_3_50}, we have 
\begin{equation}
\lambda _{n}\psi _{n}\rightarrow 0\,\,\text{in}\,\,L^{2}\left( 0,1\right) .
\label{eq_3_51}
\end{equation}
Taking the inner product in $L^2 (0,1)$ of \eqref{eq_3_19}$_4$ with
$\psi _{n}$, integrating by parts and using the boundary conditions, we remark that  
\begin{equation}
\left\langle i\overset{\sim}{\psi}_{n},\lambda _{n}\psi _{n}\right\rangle +\frac{k_2}{2}\left\Vert \psi _{n,x}\right\Vert ^{2}+k_1\left\langle \varphi _{n,x}+\psi _{n}+l\,w_{n} ,\psi_{n}\right\rangle \rightarrow 0.  \label{eq_3_52}
\end{equation}
By using \eqref{eq_4_4}, \eqref{eq_3_21}, \eqref{eq_3_51} and \eqref{eq_3_52}, we arrive at
\begin{equation}
\psi_{n,x}\rightarrow 0\,\,\text{in}\,\,L^{2}\left( 0,1\right) .  \label{eq_3_53}
\end{equation}
Because $k_3 \delta ^{2}-\left(k_3 -k_1\right) \left(k_3 -g_0\right)=0$, then, using \eqref{eq_3.48} and \eqref{eq_3_49}, we deduce that 
\begin{equation}
\overset{\sim}{w}_{n}\rightarrow 0\,\,\text{in}\,\,L^{2}\left( 0,1\right) 
\quad\hbox{and}\quad w_{n,x}-l\varphi_{n}\rightarrow 0\,\,\text{in}\,\,L^{2}\left( 0,1\right) ,\label{eq_3_54}
\end{equation}
and so, using \eqref{eq_3_21},
\begin{equation}
w_{n,x}\rightarrow 0\,\,\text{in}\,\,L^{2}\left( 0,1\right). \label{eq_3_55+}
\end{equation}
Finally, \eqref{eq_3_20}, \eqref{eq_3_24}, \eqref{eq_3_21}, \eqref{eq_3_27}, \eqref{eq_3_30}, \eqref{eq_3_50}, \eqref{eq_3_53}, \eqref{eq_3_54} and \eqref{eq_3_55+} lead to \eqref{limPhin}, which is a contradiction with \eqref{eq_4_4}.  
Hence, the proof of Theorem \ref{Theorem 3.11} is completed.
\end{proof}  

\section{Polynomial stability}

In this section, we prove the polynomial stability of system \eqref{syst_2}.
\vskip0,1truecm
\begin{theorem}\label{Theorem 4.2}
Under the assumptions \eqref{l}-\eqref{f+} and \eqref{lpi},  
and for any $m\in\mathbb{N}^*$, there exists a constant $c_{m}>0$ such that
\begin{equation}
\forall \Phi _{0}\in D\left( \mathcal{A}^m\right) ,\,\,\forall t>2,\,\,\left\Vert \Phi (t)\right\Vert _{\mathcal{H}}\leq c_{m}\left\Vert \Phi _{0}\right\Vert_{D\left( \mathcal{A}^m\right) }\left( \frac{\ln t}{t}\right) ^{\frac{m}{24-16k}}\ln t.  \label{eq_4_1}
\end{equation}
\end{theorem}
\vskip0,1truecm
\begin{proof} It is known (see \cite{liu0}) that \eqref{eq_4_1} holds if 
\begin{equation}
i\mathbb{R} \subset \rho\left( \mathcal{A}\right)\quad\hbox{and}\quad \sup_{\vert\lambda\vert\geq 1}\lambda^{-(24-16k)}\left\Vert \left( i\lambda I-\mathcal{A}\right)^{-1} \right\Vert_{\mathcal{L}\left( \mathcal{H}\right) } <\infty. \label{polyc}
\end{equation}  
We have proved in section 3 that \eqref{lpi} implies the first condition in \eqref{polyc}. So we will prove that the second condition in \eqref{polyc} is also satisfied. This will be done by contradiction using some well chosen multipliers, where some of them were used in \cite{afas1, ag1, ag2, ag6}. Let us assume that the second condition in \eqref{polyc} is false, then, there
exist sequences $\left(\Phi_{n}\right)_n\subset\,D\left(\mathcal{A}\right)$ and $\left(\lambda _{n}\right)_n\subset \,\mathbb{R}$ satisfying \eqref{eq_4_4}, \eqref{eq_4_5} and
\begin{equation}
\lim_{n\rightarrow \infty }\lambda _{n}^{24-16k}\left\Vert \left( i\lambda_{n}\,I\,-\,\mathcal{A}\right) \,\Phi _{n}\right\Vert _{\mathcal{H}}\,=\,0.\label{eq_4_6}
\end{equation}
The contradiction will be obtained by proving \eqref{limPhin}. 
Let define $\Phi_{n}$ by \eqref{Phin}. We get from \eqref{eq_4_6} that
\begin{equation}
\left\{
\begin{array}{ll}
\lambda _{n}^{24-16k}\left[ i\lambda _{n}\varphi _{n}-\overset{\sim }{\varphi }_{n}\right] \rightarrow 0\,\,&\text{in}\,\,H_k,\vspace{0.2cm}\\
\lambda _{n}^{24-16k}\left[ i\lambda _{n}\overset{\sim }{\varphi }
_{n}-k_1 \left( \varphi _{n,x}+\psi _{n}+lw_{n}\right) _{x}-lk_{3}\left(
w_{n,x}-l\varphi _{n}\right) +\delta \frac{\partial^k}{\partial x^k}\theta_{n}\right] \rightarrow 0\,\,&\text{in}\,\,L_k,\vspace{0.2cm}\\
\lambda _{n}^{24-16k}\left[ i\lambda _{n}\psi _{n}-\overset{\sim }{\psi }_{n}
\right] \rightarrow 0\,\,&\text{in}\,\,H_{1-k},\vspace{0.2cm}\\
\lambda _{n}^{24-16k}\left[ i\lambda _{n}\overset{\sim }{\psi }
_{n}-k_2\psi _{n,xx}+k_1\left( \varphi _{n,x}+\psi _{n}+lw_{n}\right) \right]
\rightarrow 0\,\,&\text{in}\,\,L_{1-k},\vspace{0.2cm}\\
\lambda _{n}^{24-16k}\left[ i\lambda _{n}w_{n}-\overset{\sim }{w}_{n}\right]
\rightarrow 0\,\,&\text{in}\,\,H_{1-k},\vspace{0.2cm}\\
\lambda _{n}^{24-16k}\left[ i\lambda _{n}\overset{\sim }{w}
_{n}-k_{3}\left( w_{n,x}-l\varphi _{n}\right) _{x}+lk_1\left( \varphi
_{n,x}+\psi _{n}+lw_{n}\right) \right] \rightarrow 0\,\,&\text{in}
\,\,L_{1-k},\vspace{0.2cm}\\
\lambda _{n}^{24-16k}\left[ i\lambda _{n}\theta _{n} -\displaystyle\int_0^{\infty}g\eta_{n,xx}ds-(-1)^k\delta \frac{\partial^k}{\partial x^k}\overset{\sim }{\varphi }_{n}\right]\rightarrow 0\,\,&\text{in}\,\,L^{2}\left( 0,1\right) ,\vspace{0.2cm}\\
\lambda _{n}^{24-16k}\left[ i\lambda _{n}\eta _{n} -\theta_n +\eta_{n,s}\right]\rightarrow 0\,\,&\text{in}\,\,L_g .
\end{array}
\right.  \label{eq_4_7}
\end{equation}
Taking the inner product of $\lambda _{n}^{24-16k}\left(
i\,\lambda _{n}\,I\,-\,\mathcal{A}\right) \,\,\Phi _{n}$ with $\Phi _{n}$
in $\mathcal{H}$ and using \eqref{g} and \eqref{dissp}, we get
\begin{equation*}
Re\left\langle \lambda _{n}^{24-16k}\left( i\,\lambda _{n}\,I\,-\,\mathcal{A}\right) \,\,\Phi _{n},\Phi _{n}\right\rangle_{\mathcal{H}} =-\frac{1}{2} \lambda _{n}^{24-16k}\int_0^{\infty} g^{\prime}\Vert\eta_{n,x}\Vert^2 ds\geq \frac{\mu_1}{2} \lambda _{n}^{24-16k}\Vert\eta_{n}\Vert_{L_g}^2.
\end{equation*}
So, \eqref{eq_4_4} and \eqref{eq_4_6} imply that
\begin{equation}
\lambda_{n}^{12-8k}\eta _{n}\longrightarrow 0\,\,\text{in}\,\,L_g .  \label{eq_4_9}
\end{equation}
We consider the function ${\hat\theta}_n$ defined in \eqref{hatthetan}.
Taking the inner product of \eqref{eq_4_7}$_8$ with $\lambda _{n}^{-(12-8k)} {\hat\theta}_{n}$ in $L_g$ and using \eqref{eq_4_4} and \eqref{eq_4_5}, we entail 
\begin{equation*}
i\left\langle \lambda _{n}^{13-8k}\eta _{n},{\hat\theta}_{n}\right\rangle_{L_g}-\lambda _{n}^{12-8k} g_0 \left\langle \theta_{n,x},{\hat\theta}_{n,x}\right\rangle+\lambda _{n}^{12-8k} \left\langle \eta_{n,s},{\hat\theta}_{n}\right\rangle_{L_g}\longrightarrow 0,
\end{equation*}
then, integrating with respect to $x$ and $s$, and using \eqref{gexp}$_1$, 
$\eta_n (s=0)=0$ and ${\hat\theta}_{n,xx} ={\theta}_{n}$, we find 
\begin{equation*}
-i\lambda _{n}^{13-8k}\int_0^{\infty}g\left\langle\eta _{n},\theta_{n}\right\rangle ds+\lambda _{n}^{12-8k} g_0 \Vert\theta_{n}\Vert^2
+\lambda _{n}^{12-8k} \int_0^{\infty}g^{\prime}\left\langle \eta_{n},\theta_{n}\right \rangle ds\longrightarrow 0,
\end{equation*}
so, using \eqref{g}, \eqref{eq_4_4}, \eqref{eq_4_5} and \eqref{eq_4_9}, we conlude that 
\begin{equation}
-i\lambda _{n}^{13-8k}\int_0^{\infty}g\left\langle\eta _{n},\theta_{n}\right\rangle ds+\lambda _{n}^{12-8k} g_0 \Vert\theta_{n}\Vert^2
\longrightarrow 0.\label{Esta}
\end{equation}
Taking the inner product of \eqref{eq_4_7}$_7$ with $\lambda _{n}^{-(12-8k)}\int_0^{\infty}g\eta_n ds$ in $L^2 (0,1)$, integrating by parts and using the boundary conditions, \eqref{eq_4_4} and \eqref{intgetan}, we get
\begin{equation*}
-i\lambda_{n}^{13-8k}\int_0^{\infty}g\left\langle\eta _{n},\theta_{n}\right\rangle ds +\lambda _{n}^{12-8k} \left\langle\int_0^{\infty}g\eta _{n,x}ds,\int_0^{\infty}g\eta_{n,x}ds\right\rangle -\delta\lambda_{n}^{12-8k}\int_0^{\infty}g\left\langle\frac{\partial^k}{\partial x^k}\eta _{n},\overset{\sim }{\varphi}_{n}\right\rangle ds \longrightarrow 0,
\end{equation*}
therefore, using \eqref{eq_4_9}, we remark that 
\begin{equation}
\lambda _{n}^{13-8k}\int_0^{\infty}g\left\langle\eta _{n},\theta_{n}\right\rangle ds\longrightarrow 0.\label{Estb}
\end{equation}
Then we deduce from \eqref{Esta} and \eqref{Estb} that 
\begin{equation}
\lambda _{n}^{6-4k}\theta _{n}\longrightarrow 0\text{ in}\,\,L^{2}\left( 0,1\right).\label{eq_4_10}
\end{equation}
Multiplying \eqref{eq_4_7}$_1$, \eqref{eq_4_7}$_3$ and \eqref{eq_4_7}$_5$ by $\lambda _{n}^{-(25-16k)}$, and using \eqref{eq_4_4} and \eqref{eq_4_5}, we obtain
\begin{equation}
\varphi _{n}\longrightarrow 0\text{ in}\,\,L_k\quad\hbox{and}\quad
\psi _{n} ,\,w_{n}\longrightarrow 0\text{ in}\,\,L_{1-k}.\label{eq_4_11}
\end{equation}
Multiplying \eqref{eq_4_7}$_3$ and \eqref{eq_4_7}$_5$ by $\lambda _{n}^{-(24-16k)}$, and using \eqref{eq_4_4} and \eqref{eq_4_5}, we have
\begin{equation}
\left( \lambda _{n}\psi _{n}\right)_{n} \,\,\hbox{and}\,\,
\left( \lambda _{n}w_{n}\right) _{n}\,\,\text{are bounded in}\,\,L_{1-k}.\label{eq_4_25}
\end{equation}
For the limits of $\varphi _{n,x}$, $\varphi _{n}$ and $\overset{\sim }{\varphi}_{n}$, we distiguish the cases $k=0$ and $k=1$.
\vskip0,1truecm
{\bf Case $k=0$}: taking the inner product of \eqref{eq_4_7}$_7$ with 
$\lambda _{n}^{-19} \varphi_{n}$ in $L^2 (0,1)$ and using \eqref{eq_4_4} and \eqref{eq_4_5}, we entail 
\begin{equation*}
i\left\langle \lambda _{n}^6\theta _{n},\varphi_{n}\right\rangle
-\lambda _{n}^5 \int_0^{\infty}g\left\langle \eta_{n,xx} ,\varphi_{n}\right\rangle -\delta\lambda _{n}^5\left\langle\overset{\sim }{\varphi}_{n},\varphi_{n}\right\rangle\longrightarrow 0,
\end{equation*}
then, integrating the second term of the above limit with respect to $x$, we arrive at 
\begin{equation*}
i\left\langle \lambda _{n}^6\theta _{n},\varphi_{n}\right\rangle
+\lambda _{n}^5 \left\langle \eta_{n} ,\varphi_{n}\right\rangle_{L_g} +\delta\lambda _{n}^5\left\langle i\lambda _{n} {\varphi}_{n} -\overset{\sim }{\varphi}_{n} ,{\varphi }_{n}\right\rangle
-i\delta\lambda _{n}^6 \Vert {\varphi}_{n}\Vert^2
\longrightarrow 0,
\end{equation*}
so, using \eqref{eq_4_4}, \eqref{eq_4_5}, \eqref{eq_4_7}$_1$, \eqref{eq_4_9} and \eqref{eq_4_10}, we conlude that 
\begin{equation}
\lambda _{n}^3\varphi _{n}\longrightarrow 0\,\,\text{in}\,\,L^2 (0,1) , \label{eq_4_16}
\end{equation}
and then, multiplying \eqref{eq_4_7}$_1$ by $\lambda _{n}^{-22}$ and using \eqref{eq_4_16}, 
\begin{equation}
\lambda _{n}^2 \overset{\sim }{\varphi }_{n}\longrightarrow 0\,\,\text{in}\,\,L^2 (0,1) . \label{eq_4_17}
\end{equation}
Taking the inner product of \eqref{eq_4_7}$_2$ with 
$\lambda _{n}^{-22} \varphi_{n}$ in $L^2 (0,1)$ and using \eqref{eq_4_5}, we find
\begin{equation*}
i\left\langle \overset{\sim}{\varphi}_{n},\lambda _{n}^{3} {\varphi}_{n}\right\rangle +k_1 \lambda _{n}^{2}\Vert{\varphi}_{n,x}\Vert^2 +  
\left\langle -k_1 \left( \psi _{n}+lw_{n}\right) _{x}-lk_{3}\left(
w_{n,x}-l\varphi _{n}\right) +\delta\theta_{n},\lambda _{n}^{2} \varphi_{n}
\right\rangle\longrightarrow 0,
\end{equation*}
using \eqref{eq_4_4} and \eqref{eq_4_16} and integrating with respect to $x$, we obtain
\begin{equation}
\lambda _{n}\varphi _{n,x}\longrightarrow 0\,\,\text{in}\,\,L^{2}\left(0,1\right).  \label{eq_4_226}
\end{equation}
Moreover, by \eqref{eq_4_7}$_1$, we see that
\begin{equation}
\overset{\sim }{\varphi }_{n,x}\rightarrow 0\text{ in}\,\,L^{2}\left(
0,1\right), \label{eq_4_21}
\end{equation}
and so, thanks to the definition of $H_0$ and Poincar\'e's inequality, \eqref{eq_4_226} and \eqref{eq_4_21} lead to
\begin{equation}
\lambda _{n}\varphi _{n}\longrightarrow 0\,\,\text{in}\,\,L^2 (0,1) \label{eq_4_166}
\end{equation}
and 
\begin{equation}
\overset{\sim }{\varphi }_{n}\longrightarrow 0\,\,\text{in}\,\,L^2 (0,1) . \label{eq_4_1760+}
\end{equation}
\vskip0,1truecm
{\bf Case $k=1$}: multiplying \eqref{eq_4_7}$_2$ by $\lambda _{n}^{-9}$ and using \eqref{eq_4_5}, we find
\begin{equation*}
i\overset{\sim }{\varphi }_{n}-\frac{k_1}{\lambda _{n}}\left( \varphi
_{n,x}+\psi _{n}+lw_{n}\right) _{x}-\frac{lk_{3}}{\lambda _{n}}\left(
w_{n,x}-l\varphi _{n}\right) + \frac{\delta}{\lambda _{n}}\theta_{n,x}
\longrightarrow 0 \text{ in}\,\,L^{2}\left( 0,1\right).
\end{equation*}
Using \eqref{eq_4_4} and \eqref{eq_4_5}, we conclude that
\begin{equation}
\left( \lambda _{n}^{-1}(k_1 \varphi _{n,xx}-\delta\theta_{n,x})\right) _{n}\,\,\text{is bounded in}\,\,L^{2}\left( 0,1\right).
\label{eq_4_13}
\end{equation} 
Taking the inner product of \eqref{eq_4_7}$_7$
with $\lambda _{n}^{-7}(k_1 \varphi _{n,x}-\delta\theta_{n})$ in $L^{2}\left( 0,1\right)$, integrating by parts and using \eqref{eq_4_4} and \eqref{eq_4_5}, we entail
\begin{equation*}
i\left\langle \lambda _{n}^2\theta _{n},k_1 \varphi _{n,x}-\delta\theta_{n}\right\rangle +\int_0^{\infty}g\left\langle \lambda _{n}^2 \eta_{n,x},
\frac{1}{\lambda _{n}}(k_1 \varphi _{n,xx}-\delta\theta_{n,x})\right\rangle
\end{equation*}
\begin{equation*}
-\delta \left\langle \lambda _{n}\left( i\lambda _{n}\varphi _{n,x}-\overset{\sim }{\varphi }_{n,x}\right) ,k_1 \varphi _{n,x}-\delta\theta_{n} \right\rangle -i\delta^2 \left\langle \varphi _{n,x} ,\lambda _{n}^2\theta_{n} \right\rangle +i\delta k_1\lambda_{n}^2\left\Vert \varphi _{n,x}\right\Vert ^{2}\longrightarrow 0,
\end{equation*}
Combination of \eqref{eq_4_4}, \eqref{eq_4_5}, \eqref{eq_4_7}$_1$, \eqref{eq_4_9}, \eqref{eq_4_10} and \eqref{eq_4_13} gives \eqref{eq_4_226}. Similarily to case $k=0$ and according to the definition of $H_1$, we see that the limit \eqref{eq_4_21}, and the limits \eqref{eq_4_166} and \eqref{eq_4_1760+} in $L_1$ (instead of $L^2 (0,1)$) hold true in case $k=1$. Taking the inner product of \eqref{eq_4_7}$_2$ with $i\lambda _{n}^{-7}\overset{\sim }{\varphi}_{n}$ in $L^{2}\left( 0,1\right)$, integrating by parts and using the boundary conditions, we get
\begin{equation*}
\left\Vert \lambda _{n}\overset{\sim }{\varphi }_{n}\right\Vert^{2}
+k_1\left\langle \lambda _{n}\left( \varphi _{n,x}+\psi _{n}+lw_{n}\right) ,i\overset{\sim }{\varphi }_{n,x}\right\rangle +lk_{3}\left\langle \lambda _{n}w_{n},i\overset{\sim }{\varphi }_{n,x}\right\rangle +l^2 k_{3}\left\langle \lambda_{n}\varphi _{n},i\overset{\sim }{\varphi }_{n}\right\rangle 
-\delta\left\langle \lambda _{n}\theta _{n},i\overset{\sim }{\varphi }_{n,x}\right\rangle\rightarrow 0.
\end{equation*}
So, using \eqref{eq_4_4}, \eqref{eq_4_5}, \eqref{eq_4_10}, \eqref{eq_4_25}, \eqref{eq_4_226}, \eqref{eq_4_21} and \eqref{eq_4_166}, we deduce that
\begin{equation}
\lambda _{n}\overset{\sim }{\varphi }_{n}\longrightarrow 0\,\,\text{in}\,\,L_1 ,  \label{eq_4_26}
\end{equation}
and by \eqref{eq_4_5} and \eqref{eq_4_7}$_1$, we find
\begin{equation}
\lambda _{n}^{2}\varphi _{n}\longrightarrow 0\;\text{in}\,\,L_1 . \label{eq_4_27}
\end{equation}
We notice that, because of \eqref{eq_4_5}, \eqref{eq_4_16} and \eqref{eq_4_17}, the limits \eqref{eq_4_26} and \eqref{eq_4_27} in 
$L^2 (0,1)$ (instead of $L_1$) are satisfied also in case $k=0$. Hence 
\begin{equation}
\lambda _{n}^2\varphi _{n} ,\,\,\lambda _{n}\overset{\sim }{\varphi }_{n}\longrightarrow 0\,\,\text{in}\,\,L_k . \label{eq_4_176}
\end{equation}
\vskip0,1truecm
Now, we ccontinue to look for \eqref{limPhin} in both cases $k=0$ and $k=1$.
Multiplying \eqref{eq_4_7}$_4$ and \eqref{eq_4_7}$_6$ by $\lambda_{n}^{-(25-16k)}$ and exploiting \eqref{eq_4_4} and \eqref{eq_4_5}, it appears that
\begin{equation}
\left(\lambda_{n}^{-1}\psi _{n,xx}\right)_{n}\,\,\hbox{and} \,\, \left(\lambda_{n}^{-1} w_{n,xx}\right)_{n}\,\,\text{are bounded in}\,\,L^{2}\left( 0,1\right).\label{eq_4_30}
\end{equation}
Taking the inner product of \eqref{eq_4_7}$_2$
with $\lambda _{n}^{-(24-16k)}\left[ k_1\psi _{n,x}+l\left( k_1 +k_{3}\right) w_{n,x}\right]$ in $L^{2}\left( 0,1\right)$, we arrive at
\begin{equation}
\left\langle i\lambda _{n}\overset{\sim }{\varphi }_{n},
k_1\psi _{n,x}+l\left( k_1 +k_{3}\right) w_{n,x} \right\rangle
+\left\langle \delta\frac{\partial^k}{\partial x^k} \theta_{n} -k_1 \varphi _{n,xx}, k_1\psi _{n,x}+l\left( k_1 +k_{3}\right)w_{n,x}\right\rangle
\label{eq_4_28}
\end{equation}
\begin{equation*}
-\left\Vert k_1\psi _{n,x}+l\left( k_1 +k_{3}\right) w_{n,x}\right\Vert^{2}
+l^{2}k_{3}\left\langle \varphi _{n},k_1\psi _{n,x}+l\left(
k_1+k_{3}\right) w_{n,x} \right\rangle\rightarrow 0. 
\end{equation*}
Again, integrating by parts and using the boundary conditions, we see that
\begin{equation}
-k_1\left\langle\varphi _{n,xx},k_1\psi _{n,x}+l\left( k_1 +k_{3}\right)w_{n,x}\ \right\rangle =k_1\left\langle \lambda _{n}\varphi _{n,x},\frac{k_1}{\lambda _{n}}\psi _{n,xx}+\frac{l\left( k_1 +k_{3}\right)}{\lambda _{n}}w_{n,xx}\right\rangle   \label{eq_4_29}
\end{equation}
and, if $k=1$,
\begin{equation}
\delta\left\langle \theta_{n,x} ,k_1\psi _{n,x}+l\left( k_1 +k_{3}\right)w_{n,x}\right\rangle =-\delta\left\langle \lambda _{n}\theta_{n},\frac{k_1}{\lambda _{n}}\psi _{n,xx}+\frac{l\left( k_1 +k_{3}\right)}{\lambda _{n}}w_{n,xx}\right\rangle  \label{eq_4_296}
\end{equation} 
(if $k=0$, no integrating by parts is needed). Then, using \eqref{eq_4_4}, \eqref{eq_4_5}, \eqref{eq_4_10}, \eqref{eq_4_226}, \eqref{eq_4_30}, \eqref{eq_4_29} and \eqref{eq_4_296}, we deduce that
\begin{equation}
\left\langle \delta\frac{\partial^k}{\partial x^k}\theta_n -k_1 \varphi _{n,xx},k_1\psi _{n,x}+l\left( k_1 +k_{3}\right)
w_{n,x}\right\rangle \rightarrow 0,  \label{eq_4_31}
\end{equation}
so, exploiting \eqref{eq_4_4}, \eqref{eq_4_5}, \eqref{eq_4_11}, \eqref{eq_4_176}, \eqref{eq_4_28} and \eqref{eq_4_31}, we entail 
\begin{equation}
k_1\psi _{n,x}+l\left( k_1 +k_{3}\right) w_{n,x}\rightarrow 0\text{ in}\,\, L^{2}\left( 0,1\right). \label{eq_4_32}
\end{equation}
Taking the inner product of \eqref{eq_4_7}$_4$ with $\lambda _{n}^{-(24-16k)}\psi _{n}$ in $L^{2}\left( 0,1\right)$, using \eqref{eq_4_4} and \eqref{eq_4_5}, integrating by parts and using the boundary conditions, we obtain
\begin{equation*}
-\left\langle \overset{\sim }{\psi }_{n},i\lambda _{n}\psi
_{n}-\overset{\sim }{\psi }_{n}\right\rangle -\left\Vert
\overset{\sim }{\psi }_{n}\right\Vert^{2}
+k_2 \left\Vert \psi_{n,x}\right\Vert^{2}+k_1 \left\langle \varphi _{n,x}+\psi_{n}+lw_{n},\psi _{n}\right\rangle \rightarrow 0,
\end{equation*}
then, using \eqref{eq_4_4}, \eqref{eq_4_5}, \eqref{eq_4_7}$_3$ and 
\eqref{eq_4_11}, we find
\begin{equation}
k_2 \left\Vert \psi _{n,x}\right\Vert^{2}-\left\Vert \overset{\sim }{\psi }_{n}\right\Vert^{2}\rightarrow 0. \label{eq_4_33}
\end{equation}
On the other hand, taking the inner product of \eqref{eq_4_7}$_6$ with 
$\lambda _{n}^{-(24-16k)}w_{n}$ in $L^{2}\left( 0,1\right)$, using \eqref{eq_4_4} and \eqref{eq_4_5},
integrating by parts and using the boundary conditions, we observe that
\begin{equation*}
-\left\langle \overset{\sim }{w}_{n},i\lambda _{n}w_{n}-\overset{\sim }{w}_{n}\right\rangle -\left\Vert \overset{
\sim }{w}_{n}\right\Vert +k_{3}\left\Vert w_{n,x}\right\Vert^{2}
-lk_{3}\left\langle \varphi_{n},w_{n,x}\right\rangle +lk_1\left\langle  \varphi _{n,x}+\psi_{n}+lw_{n} ,w_{n}\right\rangle\rightarrow 0.
\end{equation*}
By \eqref{eq_4_4}, \eqref{eq_4_5}, \eqref{eq_4_7}$_5$ and \eqref{eq_4_11}, it follows that
\begin{equation}
k_{3}\left\Vert w_{n,x}\right\Vert^{2}-\left\Vert \overset{\sim }{w
}_{n}\right\Vert^{2}\rightarrow 0.  \label{eq_4_34}
\end{equation}
Taking the inner product in $L^2 (0,1)$ of \eqref{eq_4_7}$_4$ with 
$\lambda_{n}^{-(24-16k)}w_{n}$ and of \eqref{eq_4_7}$_6$ with $\lambda _{n}^{-(24-16k)}\psi _{n}$, and using \eqref{eq_4_4} and \eqref{eq_4_5},  integrating by parts and using the boundary conditions, it appears that
\begin{equation*}
-\left\langle \overset{\sim }{\psi }_{n},i\lambda _{n}w_{n}-
\overset{\sim }{w}_{n} \right\rangle -\left\langle \overset{\sim }{\psi }_{n},\overset{\sim }{w}_{n}\right\rangle
+k_2\left\langle \psi_{n,x},w_{n,x}\right\rangle +k_1\left\langle \varphi _{n,x}+\psi_{n}+lw_{n} ,w_{n}\right\rangle \rightarrow 0
\end{equation*}
and
\begin{equation*}
-\left\langle \overset{\sim }{w}_{n},i\lambda _{n}\psi _{n}-
\overset{\sim }{\psi }_{n} \right\rangle -\left\langle
\overset{\sim }{w}_{n},\overset{\sim }{\psi }_{n}\right\rangle
+k_{3}\left\langle w_{n,x}-l\varphi _{n} ,\psi
_{n,x}\right\rangle +lk_1\left\langle \varphi _{n,x}+\psi
_{n}+lw_{n} ,\psi _{n}\right\rangle \rightarrow 0,
\end{equation*}
then, using \eqref{eq_4_4}, \eqref{eq_4_5}, \eqref{eq_4_7}$_3$, 
\eqref{eq_4_7}$_5$ and \eqref{eq_4_11}, we obtain
\begin{equation*}
k_2\left\langle \psi _{n,x},w_{n,x}\right\rangle -\left\langle \overset{\sim }{\psi }_{n},\overset{\sim }{w}_{n}\right\rangle \rightarrow 0\quad\hbox{and}\quad k_{3}\left\langle \psi _{n,x} ,w_{n,x}\right\rangle-\left\langle \overset{\sim }{\psi }_{n} ,\overset{\sim }{w}_{n}\right\rangle\rightarrow 0,
\end{equation*}
which implies that
\begin{equation}
\left(k_2 -k_3\right)\left\langle \overset{\sim }{\psi }_{n},\overset{\sim }{w}_{n}\right\rangle \rightarrow 0\quad \hbox{and}\quad
\left(k_2 -k_{3}\right)\left\langle \psi _{nx},w_{nx}\right\rangle \rightarrow 0. \label{eq_4_36}
\end{equation}
Now, we consider two cases.
\vskip0,1truecm
{\bf Case 1:} $k_2 \ne k_{3}$. From \eqref{eq_4_36}, we see that
\begin{equation}
\left\langle \overset{\sim }{\psi }_{n},\overset{\sim }{w}_{n}\right\rangle \rightarrow 0 \quad\hbox{and}\quad\left\langle \psi _{n,x},w_{n,x}\right\rangle \rightarrow 0. \label{eq_4_360}
\end{equation}
Therefore, taking the inner product in $L^2 (0,1)$ of \eqref{eq_4_32}, first, with $\psi_{n,x}$, and second, with $w_{n,x}$, and using \eqref{eq_4_4} and \eqref{eq_4_360}, we remark that
\begin{equation}
\psi _{n,x}\rightarrow 0\quad\hbox{and}\quad w_{n,x}\rightarrow 0 \text{ in}\,\, L^{2}\left( 0,1\right),\label{eq_4_37}
\end{equation}
and then, by \eqref{eq_4_33}, \eqref{eq_4_34} and \eqref{eq_4_37},
\begin{equation}
\overset{\sim }{\psi }_{n}\rightarrow 0\quad\hbox{and}\quad
\overset{\sim }{w}_{n}\rightarrow 0\text{ in}\,\, L_{1-k}.\label{eq_4_38}
\end{equation}
Finally, combining \eqref{eq_4_9}, \eqref{eq_4_10}, \eqref{eq_4_11},  \eqref{eq_4_226}, \eqref{eq_4_176}, \eqref{eq_4_37} and \eqref{eq_4_38}, we get \eqref{limPhin}, which is a contradiction with \eqref{eq_4_4}, so the second condition in \eqref{polyc} holds. 
\vskip0,1truecm
{\bf Case 2:} $k_2 =k_{3}$. Multiplying \eqref{eq_4_7}$_4$ and 
\eqref{eq_4_7}$_6$ by $\lambda_{n}^{-(21-16k)}$ and using \eqref{eq_4_4}, 
\eqref{eq_4_5}, \eqref{eq_4_7}$_3$ and \eqref{eq_4_7}$_5$, we obtain 
\begin{equation}
\left\{
\begin{array}{ll}
\lambda _{n}^{3}\left[ -\frac{1}{k_2}\lambda _{n}^{2}\psi _{n}-\psi
_{n,xx}+\frac{k_1}{k_2}\left( \varphi _{n,x}+\psi _{n}+lw_{n}\right) \right]
\rightarrow 0\,\,&\text{in}\,\,L^{2}\left( 0,1\right), \vspace{0.2cm}\\
\lambda _{n}^{3}\left[ -\frac{1}{k_2}\lambda _{n}^{2}w_{n}-\left(
w_{n,x}-l\varphi _{n}\right) _{x}+\dfrac{lk_1}{k_{3}}\left( \varphi
_{n,x}+\psi _{n}+lw_{n}\right) \right] \rightarrow 0\,\,&\text{in}\,\,L^{2}\left( 0,1\right).
\end{array}
\right.  \label{eq_4_39}
\end{equation}
Multiplying \eqref{eq_4_39}$_1$ and \eqref{eq_4_39}$_2$ by 
$\lambda_{n}^{-3}$ and using \eqref{eq_4_5}, \eqref{eq_4_11} and  
\eqref{eq_4_226}, we find
\begin{equation}
\frac{1}{k_2}\lambda _{n}^{2}\psi _{n}+\psi _{n,xx}\rightarrow 0\,\,\text{in}\,\,L^{2}\left( 0,1\right)\quad\hbox{and}\quad
\frac{1}{k_2}\lambda _{n}^{2}w_{n}+w_{n,xx}\rightarrow 0\,\,\text{in}\,\,L^{2}\left( 0,1\right).\label{eq_4_40}
\end{equation}
Multiplying \eqref{eq_4_40}$_1$ by $k_1$ and \eqref{eq_4_40}$_2$ by 
$l(k_1 +k_{3})$ and adding the obtained limits, and multiplying 
\eqref{eq_4_40}$_1$ by $k_1$ and \eqref{eq_4_40}$_2$ by $-l(k_1 +k_{3})$ and adding the obtained limits, we entail
\begin{equation}
\left\{
\begin{array}{l}
\frac{1}{k_2}\lambda _{n}^{2}\left[ k_1\psi _{n}+l(k_1 +k_{3})w_{n}\right]
+\left[ k_1\psi _{n,xx}+l(k_1 +k_{3})w_{n,xx}\right] \rightarrow 0\,\,\text{in}\,\,L^{2}\left( 0,1\right), \vspace{0.2cm}\\
\frac{1}{k_2}\lambda _{n}^{2}\left[ k_1\psi _{n}-l(k_1 +k_{3})w_{n}\right]
+\left[ k_1\psi _{n,xx}-l(k_1 +k_{3})w_{n,xx}\right] \rightarrow 0\,\,
\text{in}\,\,L^{2}\left( 0,1\right).
\end{array}
\right. \label{eq_4_41}
\end{equation}
Taking the inner product in $L^2 (0,1)$ of \eqref{eq_4_41}$_1$ and \eqref{eq_4_41}$_2$ with $k_1\psi _{n}+l(k_1 +k_{3})w_{n}$, integrating by parts and using the boundary conditions, we infer that
\begin{equation*}
\frac{1}{k_2}\left\Vert k_1\lambda _{n}\psi _{n}+l(k_1 +k_{3})\lambda
_{n}w_{n}\right\Vert^{2} -\Vert k_1\psi _{n,x}+l(k_1 +k_{3})w_{n,x}\Vert^2
\rightarrow 0
\end{equation*}
and
\begin{equation*}
\frac{1}{k_2}\left\langle \lambda_{n}^{2}\left[ k_1\psi
_{n}-l (k_1 +k_{3})w_{n}\right] ,k_1\psi _{n}+l(k_1 +k_{3})w_{n}
\right\rangle 
\end{equation*}
\begin{equation*}
-\left\langle k_1\psi _{n,x}-l(k_1 +k_{3})w_{n,x} , k_1\psi_{n,x}+l(k_1 +k_{3})w_{n,x} \right\rangle\rightarrow 0,
\end{equation*}
so, using \eqref{eq_4_4} and \eqref{eq_4_32}, it follows that
\begin{equation}
k_1\lambda _{n}\psi _{n}+l(k_1 +k_{3})\lambda _{n}w_{n}\rightarrow 0\,\,\text{in}\,\,L_{1-k}\quad\hbox{and}\quad 
k_1^{2}\left\Vert \lambda _{n}\psi _{n}\right\Vert^{2}-l^{2}(k_1 +k_{3})^{2}\left\Vert \lambda_{n}w_{n}\right\Vert^{2}\rightarrow 0.\label{eq_4_42}
\end{equation}
Taking the inner product in $L^2 (0,1)$ of \eqref{eq_4_39}$_1$ with $\lambda_{n}^{-1} w_{n}$ and \eqref{eq_4_39}$_2$ with $\lambda_{n}^{-1} \psi_{n}$, using \eqref{eq_4_4} and \eqref{eq_4_5}, integrating by parts and using the boundary conditions, we arrive at
\begin{equation}
-\frac{1}{k_2}\lambda _{n}^{4}\left\langle \psi
_{n},w_{n}\right\rangle +\lambda _{n}^{2} \left\langle \psi
_{n,x},w_{n,x}\right\rangle -\frac{k_1}{k_2}\left\langle \lambda_{n}^2\varphi_{n},w_{n,x}\right\rangle +\frac{k_1}{k_2}\left\langle \lambda _{n}\psi _{n},\lambda _{n}w_{n}\right\rangle+\frac{lk_1
}{k_2}\left\Vert \lambda _{n}w_{n}\right\Vert^{2}\rightarrow 0 \label{eq_4_43}
\end{equation}
and
\begin{equation}
-\frac{1}{k_2}\lambda _{n}^{4}\left\langle \psi_{n} ,w_{n}\right\rangle +\lambda _{n}^{2}\left\langle\psi_{n,x} ,w_{n,x}\right\rangle-l\left( 1+\frac{k_1}{k_{3}}\right)\left\langle\psi _{n,x} ,\lambda_{n}^2\varphi _{n}\right\rangle \label{eq_4_430}
\end{equation}
\begin{equation*}
+\frac{lk_1}{k_{3}}\left\Vert \lambda _{n}\psi _{n}\right\Vert^{2}+\frac{l^{2} k_1}{k_{3}}\left\langle\lambda _{n}\psi _{n}, \lambda _{n}w_{n}\right\rangle \rightarrow 0,
\end{equation*}
therefore, multiplying \eqref{eq_4_43} by $\frac{k_2 k_3}{k_1}$ and  \eqref{eq_4_430} by $-\frac{k_2 k_3}{k_1}$, adding the obtained limits and using \eqref{eq_4_4} and \eqref{eq_4_27}, it appears that
\begin{equation}
lk_{3}\left\Vert \lambda _{n}w_{n}\right\Vert^{2}-lk_2\left\Vert \lambda
_{n}\psi _{n}\right\Vert^{2} +\left( k_{3}-l^{2} k_2\right)
\left\langle \lambda _{n}\psi_{n},\lambda _{n} w_{n}\right\rangle
\rightarrow 0.  \label{eq_4_44}
\end{equation}
By taking the inner product of \eqref{eq_4_42}$_1$ with $\lambda_n \psi_n$ in $L^2 (0,1)$ and using \eqref{eq_4_25}, we have  
\begin{equation}
k_1\left\Vert \lambda _{n}\psi _{n}\right\Vert^{2}+l(k_1 +k_{3})
\left\langle \lambda _{n}w_{n},\lambda _{n}\psi _{n}\right\rangle
\rightarrow 0.\label{eq_4_45}
\end{equation}
Combining \eqref{eq_4_42}$_2$ and \eqref{eq_4_44}, we get
\begin{equation}
\frac{1}{l(k_1 +k_{3})^{2}}\left[ k_{3}k_1^{2}-k_2 l^{2}(k_1 +k_{3})^{2}\right]\left\Vert \lambda _{n}\psi _{n}\right\Vert^{2}+\left( k_{3}-l^{2} k_2\right)\left\langle \lambda _{n}w_{n},\lambda _{n}\psi _{n}\right\rangle
\rightarrow 0,\label{eq_4_450+}
\end{equation}
so, multiplying \eqref{eq_4_45} by $\frac{\left(k_1 +k_3\right)\left(k_3 -l^2 k_2\right)}{k_3}$ and \eqref{eq_4_450+} by $-\frac{l\left(k_1 +k_3 \right)^2}{k_3}$ and adding the obtained limits, we find
\begin{equation*}
\left[ k_1 k_{3}+k_2 l^{2}(k_1 +k_{3})\right] \,\left\Vert \lambda _{n}\psi
_{n}\right\Vert^{2}\rightarrow 0,
\end{equation*}
thus,
\begin{equation}
\lambda _{n}\psi_{n}\rightarrow 0\,\,\text{in}\,\, L_{1-k}\label{eq_4_46}
\end{equation}
and, using \eqref{eq_4_42}$_1$,
\begin{equation}
\lambda _{n}w_{n}\rightarrow 0\,\,\text{in}\,\, L_{1-k}. \label{eq_4_47}
\end{equation}
Using \eqref{eq_4_5}, \eqref{eq_4_7}$_3$, \eqref{eq_4_7}$_5$,
\eqref{eq_4_46} and \eqref{eq_4_47}, we deduce that
\begin{equation}
\overset{\sim }{\psi }_{n} ,\,\,\overset{\sim }{w}_{n}\rightarrow 0\,\,\text{in}\,\,L_{1-k}.\label{eq_4_48}
\end{equation}
Taking the inner product in $L^2 (0,1)$ of \eqref{eq_4_40}$_1$ with $\psi_{n}$ and  \eqref{eq_4_40}$_2$ with $w_{n}$, integrating by parts and using the boundary conditions, we entail
\begin{equation*}
\frac{1}{k_2}\left\Vert \lambda _{n}\psi _{n}\right\Vert^{2}-\left\Vert \psi _{n,x}\right\Vert^{2}\rightarrow 0 \quad\hbox{and}\quad
\frac{1}{k_2}\left\Vert \lambda _{n}w_{n}\right\Vert^{2}-\left\Vert
w_{n,x}\right\Vert^{2}\rightarrow 0,
\end{equation*}
then, from \eqref{eq_4_46} and \eqref{eq_4_47}, we find
\begin{equation}
\psi _{n,x}\rightarrow 0\,\,\text{in}\,\, L^{2}\left( 0,1\right)\quad\hbox{and}\quad w_{n,x}\rightarrow 0\,\,\text{in}\,\, L^{2}\left( 0,1\right).
\label{eq_4_49}
\end{equation}
Finally, \eqref{eq_4_9}, \eqref{eq_4_10}, \eqref{eq_4_11}, \eqref{eq_4_226},  \eqref{eq_4_176}, \eqref{eq_4_48} and \eqref{eq_4_49} imply \eqref{limPhin},
which is a contradiction with \eqref{eq_4_4}. Consequentely, Theorem \ref{Theorem 4.2} is satisfied.
\end{proof}
\vskip0,1truecm
\begin{remark}\label{remark51}
The estimate \eqref{eq_4_1} does not lead to any stability propoerty when $m=0$; that is, when the initial data $\Phi_{0}$ is in 
$D\left(\mathcal{A}^0\right)=\mathcal{H}$ (case of weak solutions). But in fact, under the assumptions of Theorem \ref{Theorem 4.2} and for $\Phi_{0} \in\mathcal{H}$, the solution $\Phi$ of system \eqref{syst_2} satisfies (strong stability)
\begin{equation}
\lim_{t\to\infty} \Vert \Phi(t)\Vert_{\mathcal{H}} =0. \label{SS51}
\end{equation}   
Indeed, a $C_{0}$-semigroup of contractions $e^{t\mathcal{A}}$ generated by an operator $\mathcal{A}$ on a Hilbert space $\mathcal{H}$ is strogly stable if $\mathcal{A}$ has no imaginary eigenvalues and 
$\sigma (\mathcal{A})\cap i\mathbb{R}$ is countable, where 
$\sigma (\mathcal{A})$ is the spectrum set of $\mathcal{A}$ (see \cite{arba}). According to the fact that $0\in \rho \left( \mathcal{A}\right)$ (proved in section 2) and since $D(\mathcal{A})$ has a compact embedding into $\mathcal{H}$, the linear bounded operator $\mathcal{A}^{-1}$ is a bijection between $\mathcal{H}$ and $D(\mathcal{A})$, and $\mathcal{A}^{-1}$ is a compact operator, which implies that $\sigma (\mathcal{A})$ is discrete and has only eigenvalues. Consequently, to get \eqref{SS51}, we have only to prove that there is no imaginary eigenvalues for $\mathcal{A}$, which was proved in section 3 under condition \eqref{lpi}. 
\end{remark}

\section{Timoshenko type thermoelastic systems}

In this last section, we consider the Timoshenko type thermoelastic systems
corresponding to the particular case \eqref{wl0}; that is
\begin{equation}
\left\{
\begin{array}{ll}
\varphi _{tt}-k_1\left( \varphi _{x}+\psi \right)
_{x} +\delta \frac{\partial^k}{\partial x^k}\theta =0 & \text{in }
\left( 0,1\right) \times \left( 0,\infty \right) , \vspace{0.2cm}\\
\psi _{tt}-k_2 \psi _{xx}+k_1\left( \varphi _{x}+\psi \right) =0 &\text{in }\left( 0,1\right) \times \left( 0,\infty \right) , \vspace{0.2cm}\\
\theta _{t} -\displaystyle\int_0^{\infty}f(s)\theta _{xx} (x,t-s)ds -(-1)^k\delta \frac{\partial^k}{\partial x^k}\varphi _{t} =0 & \text{in }\left(
0,1\right) \times \left( 0,\infty \right)
\end{array}
\right. \label{syst16}
\end{equation}
along with the homogeneous Dirichlet-Neumann boundary conditions and the initial data
\begin{equation}
\left\{
\begin{array}{ll}
\frac{\partial^k}{\partial x^k}\varphi \left( 0,t\right) =\,\frac{\partial^{1-k}}{\partial x^{1-k}}\psi \left( 0,t\right)=\theta \left( 0,t\right) =0 & \text{in }\left( 0,\infty \right),
\vspace{0.2cm}\\
\frac{\partial^k}{\partial x^k}\varphi \left( 1,t\right) =\,\frac{\partial^{1-k}}{\partial x^{1-k}}\psi \left( 1,t\right)=\theta\left( 1,t\right) =0 & \text{in }\left( 0,\infty \right),
\vspace{0.2cm}\\
\varphi \left( x,0\right) =\varphi _{0}\left( x\right) ,\,\varphi _{t}\left(
x,0\right) =\varphi _{1}\left( x\right),\,\psi \left( x,0\right) =\psi _{0}\left( x\right) ,\,\psi _{t}\left(x,0\right) =\psi _{1}\left( x\right) & \text{in }\left( 0,1\right) , \vspace{0.2cm}\\
\,\theta \left( x,-t\right) =\theta_{0}\left( x,t\right) & \text{in }\left( 0,1 \right)\times \left( 0,\infty \right),
\end{array}
\right.  \label{cdt_106}
\end{equation}
where $k\in\{0,1\}$ and, without lose of generality, \eqref{rhojL} is considered.   
\vskip0,1truecm
\subsection{Well-posedness} As in section 2, system 
\eqref{syst16}-\eqref{cdt_106} can be formulated in the form \eqref{syst_2}, where $\eta$ and its initial data $\eta_0$ are still defined by \eqref{etaname}, the components $w$ and $\tilde w$, theirs initial data
$w_0$ and $w_1$, theirs spaces $H_{1-k}$ and $L_{1-k}$ and theirs lines \eqref{A}$_5$-\eqref{A}$_6$ are not considered in the definitions of 
$\mathcal{H}$, $\Phi$, $\Phi_0$, $\mathcal{A}$ and $D(\mathcal{A})$, and 
$w_j ={\tilde w}_j =l=0$ in \eqref{innerproductH}. We remark that, if 
$(\varphi ,\psi)\in H_k\times H_{1-k}$ satisfying (instead of \eqref{Nvarphipsiw}) 
\begin{equation*}
\mathcal{N} (\varphi,\psi):=k_1\left\Vert \varphi_{x}+\psi \right\Vert^2 +k_2\left\Vert \psi _{x} \right\Vert^2 =0,
\end{equation*}
then $\varphi_{x}+\psi =\psi _{x} =0.$ Therefore, the definition of 
$H_{1-k}$ implies that $\psi =0$, so $\varphi_{x}=0$, thus, using
the definition of $H_{k}$, we get $\varphi =0$. Hence 
$(\mathcal{H},\left\langle\cdot,\cdot\right\rangle_{\mathcal{H}})$
is a Hilbert space. Then, using the same arguments as in section 2, we 
get \eqref{dissp} and $0\in \rho(\mathcal{A})$. Consequently, $\mathcal{A}$
generates a semigroup of contractions on $\mathcal{H}$. Finally, the next  theorem is valid.    
\vskip0,1truecm
\begin{theorem}\label{Theorem 6.1}
Under assumptions \eqref{f}-\eqref{f+}, and for any $m\in \mathbb{N}$ and 
$\Phi_0 \in D(\mathcal{A}^m)$, system \eqref{syst_2} in case 
\eqref{syst16}-\eqref{cdt_106} admits a unique solution $\Phi$ satisfying \eqref{exist}. 
\end{theorem}
\vskip0,1truecm
For the stability problem, we prove the following stability results: 
\vskip0,1truecm
\begin{theorem}\label{Theorem 6.2}
Under the assumptions \eqref{f}-\eqref{f+}, the following stability results
are valid:
\vskip0,1truecm
1. System \eqref{syst_2} in case \eqref{syst16}-\eqref{cdt_106} with $k=1$ is  exponentially stable if and only if 
\begin{equation}
\delta^2 =\frac{(k_2 -k_1)(k_2 -g_0 )}{k_2} .\label{chi10}
\end{equation}
\vskip0,1truecm
2. System \eqref{syst_2} in case \eqref{syst16}-\eqref{cdt_106} with $k=0$ is  exponentially stable if and only if 
\begin{equation}
k_1 =k_2 =g_0 .\label{chi100}
\end{equation}
\vskip0,1truecm
3. System \eqref{syst_2} in case \eqref{syst16}-\eqref{cdt_106} satisfies the 
polynomial and strong stability estimates \eqref{eq_4_1} and \eqref{SS51}.    
\end{theorem}
\vskip0,1truecm
\subsection{Lack of exponential stability} In this subsection, we prove that \eqref{syst_2} in case \eqref{syst16}-\eqref{cdt_106} is not exponentially stable if [$k=1$ and \eqref{chi10} does not hold] or [$k=0$ and \eqref{chi100} does not hold]. It is enough to prove that the second condition in \eqref{expon} is not satisfied. The proof is very similar to the one given in section 3. We prove the existence of sequences 
\begin{equation*}
(\lambda_{n})_n \subset \mathbb{R}, \quad (\Phi_{n})_n \subset D\left(\mathcal{A}\right)\quad\hbox{and}\quad (F_{n})_n \subset \mathcal{H}
\end{equation*}
satisfying \eqref{Fn} and \eqref{eq_4}, where
\begin{equation}
w_n = \overset{\sim }{w}_{n} =l =0 \label{wntildewnl}
\end{equation} 
and 
\begin{equation}
f_{5,n}=f_{6,n}=0; \label{wntildewnl+}
\end{equation} 
that is, in case \eqref{syst16}-\eqref{cdt_106}, \eqref{eq_2_5}$_5$ and \eqref{eq_2_5}$_6$ are empty. We distinguish two cases. 
\vskip0,1truecm
{\bf Case 1: [$k=1$ and \eqref{chi10} does not hold] or [$k=0$ and $k_1\ne k_2$]}. Considering the choices \eqref{fn0}-\eqref{fn00} (with 
$\alpha_{3,n}=\beta_{6,n}=0$), \eqref{eq_2_6} is reduced to the system
\begin{equation}
\left\{
\begin{array}{l}
\left[(k_1 +\mu_{k,n}) N^2 -\lambda_{n}^{2}\right]\alpha_{1,n} +k_1 N\alpha_{2,n} =\beta_{2,n}, \vspace{0.2cm}\\
k_1 N\alpha_{1,n} +\left(k_2 N^2 -\lambda _{n}^{2}+k_1\right)\alpha_{2,n} =\beta_{4,n} .  
\end{array}
\right. \label{eq_2_66}
\end{equation}
As in subsection 3.2, we pick 
$\beta_{4,n}:=\beta_{4}\in (0,\infty)$ not depending on $n$, and choose
\begin{equation}
\beta_{2,n} =0\quad\hbox{and}\quad
\lambda_n =\sqrt{k_2 N^2 +k_1} ,\label{fn0650+06}
\end{equation}
then \eqref{mun00} is satisfied. It is clear that \eqref{eq_2_66} is equivalent to
\begin{equation}
\alpha_{1,n} =\frac{\beta_4}{k_1 N} \quad\hbox{and}\quad 
\alpha_{2,n} =\frac{\beta_4 \left[(k_2 -k_1 -\mu_{k,n}) N^2 +k_1 \right]}{k_1^2 N^2} ,\label{eq_2_6306}
\end{equation}
therefore, because \eqref{chi10} does not hold in case $k=1$, and $k_1 \ne k_2$ in case $k=0$, \eqref{mun00} and \eqref{eq_2_6306} lead to
\begin{equation}
\vert\alpha_{2,n}\vert \sim \left\{
\begin{array}{ll}
\frac{\vert k_2 -k_1\vert \beta_{4}}{k_1^2} & \hbox{if}\quad
k=0, \vspace{0.2cm}\\ 
\frac{\delta^2 k_2{\sqrt{k_2}}\beta_{4} N}{k_1^2 g(0)} & \hbox{if}\quad k=1
\,\,\hbox{and}\,\,g_0 = k_2 , \vspace{0.2cm}\\ 
\frac{\left\vert(k_2 -k_1 )(k_2 -g_0 )-\delta^2 k_2\right\vert\beta_4}{k_1^2 \vert k_2 -g_0 \vert } & \hbox{if}\quad k=1\,\,\hbox{and}\,\,g_0 \ne k_2 .
\end{array}
\right.  \label{eq_2_71000666} 
\end{equation}
Using \eqref{eq_2_71000666}, \eqref{eq_2_710006} implies \eqref{fn69}. Moreover, for $k=1$, we take $\beta_4 =1$ and we find \eqref{eq_2_7100066}, which implies \eqref{Fn}. For $k=0$, we see that \eqref{fn00}, \eqref{fn0650+06} and \eqref{eq_2_6306} imply \eqref{fn065+0}, where 
\begin{equation*}
a_{3,n} =\frac{\delta^2 (k_2 N^2 +k_1 )}{k_1^2 N^2}. 
\end{equation*}
The fact that $(a_{3,n})_n$ is bounded allows us to define $\beta_4$ by \eqref{beta4fn065+0}, and then \eqref{fn065+0} leads to \eqref{Fn}. This proves the lack of exponential stability for \eqref{syst16}-\eqref{cdt_106} when $k=1$ and \eqref{chi10} does not hold.  
\vskip0,1truecm
{\bf Case 2: $k=0$, $k_1 = k_2$ and $k_1 \ne g_0$}. We consider \eqref{fn0} and
\begin{equation}
\left\{
\begin{array}{l}
\varphi_{n} (x)=\alpha_{1,n}\sin\,(Nx),\quad \psi_{n} (x)=\alpha_{2,n}\cos\,(Nx),\quad \theta_n =\alpha_{4,n}\sin\,(Nx) ,\vspace{0.2cm}\\
f_{2,n} (x)=\sin\,(Nx) ,\quad f_{4,n} (x) = f_{7,n} = 0. 
\end{array}
\right.  \label{62fn000}
\end{equation}
These choices imply that \eqref{eq_2_5}$_1$, \eqref{eq_2_5}$_3$ and \eqref{eq_2_5}$_8$ are satisfied, \eqref{eq_2_5}$_7$ is equivalent to
\begin{equation*}
\alpha_{4,n} =\frac{\delta\lambda_n^2}{\lambda_n^2 - N^2 (g_0 -\mu_{2,n} )} \alpha_{1,n}  
\end{equation*} 
($\mu_{2,n}$ is defined in \eqref{mukn2n}) and \eqref{eq_2_5}$_2$ and \eqref{eq_2_5}$_4$ hold if and only if
\begin{equation}
\left\{
\begin{array}{l}
\left[k_1 N^2 -\lambda_{n}^{2} +\frac{\delta^2\lambda_n^2}{\lambda_n^2 - N^2 (g_0 -\mu_{2,n} )}\right]\alpha_{1,n} +k_1 N\alpha_{2,n} =1, \vspace{0.2cm}\\
k_1 N\alpha_{1,n} +\left(k_1 N^2 -\lambda _{n}^{2}+k_1\right)\alpha_{2,n} =0.  
\end{array}
\right. \label{62eq_2_66}
\end{equation}  
We take
\begin{equation}
\lambda_n =\sqrt{k_1 N^2 +{\sqrt{k_1}}N+\frac{1}{2}k_1} .\label{62fn0650+06}
\end{equation}  
We see that
\begin{equation*}
\left\Vert F_{n}\right\Vert _{\mathcal{H}}^{2} =\left\Vert f_{2,n}\right\Vert^{2} =\displaystyle\int_{0}^{1} \sin^2 \,(Nx) dx\leq 1 ,
\end{equation*}  
which gives \eqref{Fn}. On the other hand, \eqref{62eq_2_66} is equivalent to
\begin{equation}
\alpha_{2,n} =\frac{k_1 N}{{\sqrt{k_1}}N-\frac{1}{2}k_1}\alpha_{1,n} \quad\hbox{and}\quad 
\alpha_{1,n} =\frac{{\sqrt{k_1}}N-\frac{1}{2}k_1}{\left({\sqrt{k_1}}N-\frac{1}{2}k_1 \right)\frac{\delta^2\lambda_n^2}{\lambda_n^2 - N^2 (g_0 -\mu_{2,n} )}+\frac{1}{4}k_1^2 } ,\label{62eq_2_6306}
\end{equation}
then, because $k_1 \ne g_0$ and according to \eqref{mu2n}, we have 
\begin{equation*}
\alpha_{1,n} \sim \frac{k_1 -g_0}{k_1 \delta^2} \quad\hbox{and}\quad 
\alpha_{2,n} \sim \frac{{\sqrt{k_1}}(k_1 -g_0)}{k_1 \delta^2},
\end{equation*}
so the equivalence of $\alpha_{2,n}$ implies \eqref{fn69} because 
\begin{equation*}
\left\Vert \Phi_{n}\right\Vert _{\mathcal{H}}^{2} \geq k_2\Vert\psi_{n,x}\Vert^2 =\frac{1}{2}k_1 \vert\alpha_{2,n}\vert^2 N^2\displaystyle\int_{0}^{1} (1-\cos \,(2Nx) )dx= \frac{1}{2}k_1 \vert\alpha_{2,n}\vert^2 N^2 .
\end{equation*} 
Consequently, the proof of the lack of exponential stability for \eqref{syst16}-\eqref{cdt_106} when $k=0$ and \eqref{chi100} does not hold is achieved.
\vskip0,1truecm
\subsection{Exponential stability} We prove here that \eqref{syst_2} in case \eqref{syst16}-\eqref{cdt_106} is exponentially stable if [$k=1$ and \eqref{chi10} is satisfied] or [$k=0$ and \eqref{chi100} is satisfied]. As in section 4, we prove that \eqref{expon} holds. 
\vskip0,1truecm
The first condition in \eqref{expon} is staisfied. Indeed, it is enough to prove that the unique solution of \eqref{eq_3_4*} (with 
$\lambda\in \mathbb{R}^*$) is $\Phi=0$. In case \eqref{syst16}-\eqref{cdt_106}, \eqref{eq_3_5} is reduced to 
\begin{equation}
\left\{
\begin{array}{l}
{\tilde{\varphi}}-i\lambda\varphi = {\tilde{\psi}}-i\lambda\psi =0, \vspace{0.2cm}\\
k_1\left( \varphi _{x}+\psi \right) _{x} -\delta\frac{\partial^k}{\partial x^k}\theta =i\lambda {\tilde{\varphi}} ,\vspace{0.2cm}\\
k_2\psi _{xx}-k_1\left( \varphi _{x}+\psi
\right) =i\lambda {\tilde{\psi}}, \vspace{0.2cm}\\
\displaystyle\int_0^{\infty}g\eta_{xx}ds+(-1)^k \delta \frac{\partial^k}{\partial x^k}{\tilde{\varphi}}=i\lambda\theta ,\vspace{0.2cm}\\
\theta -\eta_s =i\lambda\eta.
\end{array}
\right. \label{eq_3_56}
\end{equation}
As in section 3, using \eqref{dissp}, the definition of $H_1$ (in case $k=1$), the first equation and the last two ones in \eqref{eq_3_56}, we find $\varphi={\tilde{\varphi}}=\eta =\theta =0$, and therefore \eqref{eq_3_56} becomes 
\begin{equation}
{\tilde{\psi}}-i\lambda\psi =\psi _{x}=
k_2\psi _{xx} +\left(\lambda^{2} -k_1 \right)\psi  =0. \label{eq_3_11006}
\end{equation}
From the definition of $H_{1-k}$ and Poincar\'e's inequality, we see that \eqref{eq_3_11006} implies $\psi ={\tilde{\psi}} =0$. Consequently $\Phi=0$.  
\vskip0,1truecm
Now, we prove that the second condition in \eqref{expon} is satisfied under conditions \eqref{chi10} and \eqref{chi100} for $k=1$ and $k=0$, respectively. We distinguish the two cases $k=1$ and $k=0$, where we have \eqref{wntildewnl} and we apply similar arguments of proof as in section 4. 
\vskip0,1truecm
{\bf Case 1: $k=1$}. Because of \eqref{wntildewnl}, the computations used in section 4 are still satisfied (where \eqref{eq_3_49} holds without any restriction on the parameters, and \eqref{eq_3_50} is satisfied thanks to \eqref{chi10}), this proves the exponential stability of \eqref{syst16}-\eqref{cdt_106} if $k=1$ and \eqref{chi10} is valid.
\vskip0,1truecm
{\bf Case 2: $k=0$}. As in section 4, let us consider \eqref{eq_4_4}-\eqref{Phin} (with \eqref{wntildewnl}) and prove \eqref{limPhin}. The limit \eqref{eq_3_18} in case $k=0$ and \eqref{wntildewnl} implies the following convergences: 
\begin{equation}
\left\{
\begin{array}{ll}
i\lambda _{n}\varphi _{n}-\overset{\sim }{\varphi }_{n} \rightarrow 0\,\,&\text{in}\,\,H_0,\vspace{0.2cm}\\
i\lambda _{n}\overset{\sim }{\varphi }
_{n}-k_1 \left( \varphi _{n,x}+\psi _{n}\right) _{x} +\delta\theta_{n} \rightarrow 0\,\,&\text{in}\,\,L^2 (0,1),\vspace{0.2cm}\\
i\lambda _{n}\psi _{n}-\overset{\sim }{\psi }_{n}
\rightarrow 0\,\,&\text{in}\,\,H_{1},\vspace{0.2cm}\\
i\lambda _{n}\overset{\sim }{\psi }
_{n}-k_2\psi _{n,xx}+k_1\left( \varphi _{n,x}+\psi _{n}\right)
\rightarrow 0\,\,&\text{in}\,\,L_{1},\vspace{0.2cm}\\
i\lambda _{n}\theta _{n} -\displaystyle\int_0^{\infty}g\eta_{n,xx}ds-\delta \overset{\sim }{\varphi }_{n}\rightarrow 0\,\,&\text{in}\,\,L^{2}\left( 0,1\right) ,\vspace{0.2cm}\\
i\lambda_{n}\eta_{n} -\theta_n +\eta_{n,s}\rightarrow 0\,\,&\text{in}\,\,L_g .
\end{array}
\right.  \label{6eq_3_19}
\end{equation}
Taking the inner product of $i\,\lambda _{n}\,I\,-\,\mathcal{A} \,\,\Phi _{n}$ with $\Phi _{n}$ in $\mathcal{H}$ and using \eqref{g} and \eqref{dissp}, we get \eqref{eq_3_20}. Using \eqref{hatthetan}, we obtain \eqref{eq_4_100} as in section 4. On the other hand, taking the inner product of \eqref{6eq_3_19}$_5$ with $\int_0^{\infty} g\eta_n ds$ in $L^2 (0,1)$, integrating by parts and using \eqref{eq_4_4}, \eqref{intgetan} and the boundary conditions, we entail 
\begin{equation*}
-i\lambda _{n}\int_0^{\infty} g\left\langle \eta _{n},\theta_{n}\right\rangle ds+\left\langle \int_0^{\infty} g\eta _{n,x} ds,\int_0^{\infty} g\eta _{n,x} ds\right\rangle -\delta\int_0^{\infty} g\left\langle \eta _{n},\overset{\sim }{\varphi}_{n}\right\rangle ds\longrightarrow 0,
\end{equation*}
then, using \eqref{eq_4_4} and \eqref{eq_3_20}, we get \eqref{6eq_3_24},
so, by combinig \eqref{eq_4_100} and \eqref{6eq_3_24}, we find \eqref{eq_3_24}. Similarily to \eqref{eq_3_21} and \eqref{eq_3_25}, by multiplying \eqref{6eq_3_19}$_1$, \eqref{6eq_3_19}$_2$ and \eqref{6eq_3_19}$_3$ by $\lambda_n^{-1}$, we have 
\begin{equation}
\varphi _{n}\longrightarrow 0\,\,\text{in}\,\,L^2 (0,1) , \quad
\psi _{n} \longrightarrow 0\,\,\text{in}\,\,L_{1} \label{6eq_3_21}
\end{equation}
and 
\begin{equation}
\left( \lambda_n^{-1} \varphi_{n,xx}\right)_{n}\,\,\hbox{is bounded in}\,L^2 (0,1).\label{6eq_3_25}
\end{equation}
Thanks to \eqref{eq_4_4} and \eqref{eq_4_5}, we observe that  
\begin{equation}
\lambda _{n}^{-1}\overset{\sim }{\varphi}_{n}\longrightarrow 0\,\, \text{in}\,\,L^2 (0,1).\label{70}
\end{equation}
Taking the inner product of \eqref{6eq_3_19}$_2$
with $\lambda _{n}^{-1}\overset{\sim }{\varphi}_{n}$ in $L^{2}\left( 0,1\right)$, integrating by parts and using \eqref{70} and the boundary conditions, we get
\begin{equation*}
i\left\Vert \overset{\sim }{\varphi }_{n}\right\Vert^{2}
+\left\langle \lambda _{n}^{-1}\left( k_1\varphi _{n,xx}-\delta\theta_{n}\right) ,i\lambda_n \varphi_n -\overset{\sim }{\varphi }_{n}\right\rangle 
-ik_1 \left\Vert \varphi_{n,x}\right\Vert^{2} -i\delta\left\langle \theta_{n} ,\varphi_{n} \right\rangle-\left\langle k_1 \psi _{n,x} ,\lambda _{n}^{-1}\overset{\sim }{\varphi }_{n}\right\rangle\rightarrow 0,
\end{equation*}
thus, using \eqref{eq_4_4}, \eqref{6eq_3_19}$_1$, \eqref{6eq_3_21}, \eqref{6eq_3_25} and \eqref{70}, we deduce that
\begin{equation}
\left\Vert \overset{\sim }{\varphi }_{n}\right\Vert^{2}
-k_1 \left\Vert \varphi_{n,x}\right\Vert^{2}\longrightarrow 0. \label{6eq_3_30}
\end{equation}
From \eqref{eq_4_4}, \eqref{6eq_3_19}$_{1}$ and \eqref{6eq_3_19}$_{3}$, we observe that 
\begin{equation}
\left( \Vert \lambda _{n}\varphi _{n}\Vert \right)_{n}\,\,\hbox{and}\,\, \left( \Vert \lambda _{n}\psi _{n}\Vert \right)_{n} \,\,\hbox{are bounded.}\label{6psiwlambdanbound}
\end{equation}
We have, by integrating by parts and using the boundary conditions, 
\begin{eqnarray*}
\left\langle \lambda _{n}^{2}\varphi_{n}+i\lambda _{n} \overset{\sim}{\varphi}_{n},i\theta _{n}\right\rangle &=&-i\left\langle 
i\lambda _{n}\varphi_{n}-\overset{\sim}{\varphi}_{n} ,i\lambda _{n}\theta _{n}\right\rangle \\
&=&-i\left\langle i\lambda _{n} \varphi_{n}-\overset{\sim}{\varphi}_{n} ,i\lambda _{n}\theta _{n}-\int_0^{\infty}g\eta_{n,xx} ds-\delta
\overset{\sim}{\varphi}_{n} \right\rangle \\
&&-i\delta \left\langle i\lambda _{n} \varphi_{n}-\overset{\sim}{\varphi}_{n},\overset{\sim}{\varphi}_{n}\right\rangle +i\int_0^{\infty}g\left\langle i\lambda _{n} \varphi_{n,x}-\overset{\sim}{\varphi}_{n,x},\eta_{n,x}\right\rangle ds,
\end{eqnarray*}
by using \eqref{eq_4_4}, \eqref{6eq_3_19}$_{1}$ and \eqref{6eq_3_19}$_{5}$, we obtain  
\begin{equation}
\left\langle \lambda _{n}^{2}\varphi_{n}+i\lambda _{n} \overset{\sim}{\varphi}_{n},i\theta _{n}\right\rangle \longrightarrow 0.\label{71}
\end{equation}
Also, we see that 
\begin{equation}
\left\langle \lambda _{n} \varphi_{n},\overset{\sim}{\varphi}_{n}\right\rangle= \left\langle i\lambda _{n} \varphi_{n} -\overset{\sim}{\varphi}_{n},i\overset{\sim}{\varphi}_{n}\right\rangle -i\Vert\overset{\sim}{\varphi}_{n}\Vert^2 .\label{6eq37}
\end{equation}
Using again integration by parts and the boundary conditions, we have
\begin{equation*}
\lambda _{n}\int_0^{\infty}g\left\langle \varphi_{n,x},\eta_{n,x}\right\rangle ds =-\lambda_{n}\int_0^{\infty}g\left\langle \varphi_{n},\eta_{n,xx}\right\rangle ds =-\lambda_{n}\left\langle \varphi_{n},\int_0^{\infty}g\eta_{n,xx} ds\right\rangle
\end{equation*}
\begin{eqnarray*}
&=&\lambda _{n}\left\langle \varphi_{n},i\lambda _{n}\theta
_{n}-\int_0^{\infty}g\eta_{n,xx}ds-\delta \overset{\sim}{\varphi}_{n}\right\rangle -\lambda _{n}\left\langle \varphi_{n},i\lambda _{n}\theta _{n}\right\rangle +\delta  \left\langle\lambda_{n} \varphi_{n} ,\overset{\sim}{\varphi}_{n}\right\rangle \\
&=& \left\langle \lambda _{n} \varphi_{n},i\lambda _{n}\theta
_{n}-\int_0^{\infty}g\eta_{n,xx}ds-\delta \overset{\sim}{\varphi}_{n} \right\rangle -\left\langle \lambda _{n}^{2} \varphi_{n}+i\lambda _{n}\overset{\sim}{\varphi}_{n} ,i\theta _{n}\right\rangle \\
&&+\left\langle i\lambda _{n} \overset{\sim}{\varphi}_{n}-k_1 \left(\varphi_{n,x}+\psi _{n}\right)_x +\delta \theta _{n}
,i\theta _{n}\right\rangle \\
&&+k_1\left\langle \varphi_{n,xx},i\theta_{n}\right\rangle +\left\langle k_1\psi_{n,x}+\delta\theta_{n} ,i\theta_{n}\right\rangle +\delta\left\langle
\lambda_{n} \varphi_{n} ,\overset{\sim}{\varphi}_{n}\right\rangle ,
\end{eqnarray*}
then, by \eqref{6eq37}, integration by parts and using the boundary conditions, we obtain 
\begin{equation}
\lambda _{n}\int_0^{\infty}g\left\langle \varphi_{n,x},\eta_{n,x}\right\rangle ds = \left\langle \lambda _{n} \varphi_{n},i\lambda _{n}\theta
_{n}-\int_0^{\infty}g\eta_{n,xx}ds-\delta \overset{\sim}{\varphi}_{n} \right\rangle -\left\langle \lambda _{n}^{2} \varphi_{n}+i\lambda _{n}\overset{\sim}{\varphi}_{n} ,i\theta _{n}\right\rangle \label{6eq_3.38}
\end{equation}
\begin{eqnarray*}
&&+\left\langle i\lambda _{n} \overset{\sim}{\varphi}_{n}-k_1 \left(\varphi_{n,x}+\psi _{n}\right)_x +\delta \theta _{n}
,i\theta _{n}\right\rangle -k_1 \left\langle\varphi_{n,x} , i\theta _{n,x}\right\rangle\\
&&+\left\langle k_1\psi_{n,x}+\delta \theta _{n} ,i\theta _{n}\right\rangle +\delta\left\langle\lambda_{n} \varphi_{n} -\overset{\sim}{\varphi}_{n} ,i\overset{\sim}{\varphi}_{n}\right\rangle -i\delta\Vert\overset{\sim}{\varphi}_{n}\Vert^2 ,
\end{eqnarray*}
using \eqref{eq_4_4}, \eqref{6eq_3_19}$_{1}$, \eqref{6eq_3_19}$_{2}$, \eqref{6eq_3_19}$_{5}$, \eqref{eq_3_24}, \eqref{6psiwlambdanbound} and \eqref{71}, we deduce from \eqref{6eq_3.38} that 
\begin{equation}
\lambda _{n}\int_0^{\infty}g\left\langle \varphi_{n,x},\eta_{n,x}\right\rangle ds -ik_1\left\langle \varphi_{n,x},\theta _{n,x}\right\rangle +i\delta\Vert\overset{\sim}{\varphi}_{n}\Vert^2 \rightarrow 0.  \label{76eq_3_39}
\end{equation}
Also, by integrating with respect to $s$ and using \eqref{gexp}$_1$ and 
$\eta_n (s=0)=0$, we have  
\begin{equation}
\lambda _{n}\int_0^{\infty}g\left\langle \varphi_{n,x},\eta_{n,x}\right\rangle ds = i\int_0^{\infty}g\left\langle \varphi_{n,x},i\lambda_{n}\eta_{n,x}\right\rangle ds \label{77eq_3_4000}
\end{equation}
\begin{equation*}
= i\int_0^{\infty}g\left\langle \varphi_{n,x},i\lambda _{n}\eta_{n,x} -\theta_{n,x} +\eta_{n,xs} \right\rangle ds 
+i\int_0^{\infty}g^{\prime}\left\langle \varphi_{n,x},\eta_{n,x}\right\rangle ds+ig_0 \left\langle \varphi_{n,x}, \theta_{n,x}\right\rangle,    
\end{equation*}
therefore, by using \eqref{g}, \eqref{eq_4_4}, \eqref{eq_3_20}, \eqref{6eq_3_19}$_{6}$, \eqref{76eq_3_39} and \eqref{77eq_3_4000}, we obtain 
\begin{equation*}
\left( k_1 -g_0 \right) \left\langle \varphi_{n,x},\theta _{n,x}\right\rangle -\delta \Vert\overset{\sim}{\varphi}_{n}\Vert^2 \rightarrow 0,
\end{equation*} 
thus, because $k_1 =g_0$, 
\begin{equation}
\overset{\sim}{\varphi}_{n}\rightarrow 0\,\,\hbox{in}\,\,L^2 (0,1),\label{7aa}
\end{equation}
so \eqref{6eq_3_30} and \eqref{7aa} imply that 
\begin{equation}
\varphi_{n,x}\rightarrow 0\,\,\hbox{in}\,\,L^2 (0,1),\label{7b}
\end{equation}
and therefore, by \eqref{6eq_3_19}$_{1}$ and \eqref{7aa},
\begin{equation}
\lambda_n \varphi_{n}\rightarrow 0\,\,\hbox{in}\,\,L^2 (0,1).\label{7c}
\end{equation}
Taking the inner product of \eqref{6eq_3_19}$_{4}$ with 
$\varphi _{n,x}+\psi _{n}$ in $L^2 (0,1)$, using \eqref{eq_4_4} and  integration by parts and using the boundary conditions, we get 
\begin{equation*}
-\left\langle \overset{\sim}{\psi}_{n},i\lambda _{n}\varphi _{n,x}\right\rangle -\left\langle \overset{\sim}{\psi}_{n}, i\lambda _{n}\psi _{n}-\overset{\sim}{\psi}_{n}\right\rangle -\left\Vert \overset{\sim}{\psi}_{n}\right\Vert ^{2} +k_2\left\langle \psi _{n,x},\left( \varphi _{n,x}+\psi _{n}\right)_{x}\right\rangle +k_1\Vert\varphi_{n,x} +\psi_n\Vert^2\rightarrow 0,
\end{equation*}
using \eqref{eq_4_4}, \eqref{6eq_3_19}$_{3}$, \eqref{6eq_3_21} and \eqref{7b}, we deduce that 
\begin{equation*}
-\left\langle \overset{\sim}{\psi}_{n},i\lambda _{n}\varphi _{n,x}\right\rangle -\left\Vert \overset{\sim}{\psi}_{n}\right\Vert ^{2} +\frac{k_2}{k_1}\left\langle \psi _{n,x},i\lambda _{n} \overset{\sim}{\varphi}_{n} +\delta \theta_{n}\right\rangle -\frac{k_2}{k_1}\left\langle \psi_{n,x},i\lambda _{n} \overset{\sim}{\varphi}_{n} -k_1\left( \varphi_{n,x}+\psi_{n}\right)_{x} +\delta\theta_n\right\rangle\rightarrow 0,
\end{equation*}
using \eqref{eq_4_4}, \eqref{eq_3_24} and \eqref{6eq_3_19}$_{2}$, we have 
\begin{equation}
-\left\langle \overset{\sim}{\psi}_{n},i\lambda _{n}\varphi _{n,x}\right\rangle -\left\Vert \overset{\sim}{\psi}_{n}\right\Vert ^{2} -\frac{k_2}{k_1}\left\langle i\lambda_{n}\psi _{n,x}-\overset{\sim}{\psi}_{n,x} ,\overset{\sim}{\varphi}_{n}\right\rangle -\frac{k_2}{k_1}\left\langle \overset{\sim}{\psi}_{n,x},\overset{\sim}{\varphi}_{n}\right\rangle \rightarrow 0. \label{6eq_3_44+}
\end{equation}
As, by integrating by parts and using the boundary conditions, 
\begin{equation*}
\left\langle \overset{\sim}{\psi}_{n,x},\overset{\sim}{\varphi}_{n}\right\rangle =-\left\langle \overset{\sim}{\psi}_{n},\overset{\sim}{\varphi}_{n,x}\right\rangle =\left\langle \overset{\sim}{\psi}_{n},i\lambda
_{n}\varphi _{n,x}-\overset{\sim}{\varphi}_{n,x}\right\rangle -\left\langle \overset{\sim}{\psi}_{n},i\lambda_{n}\varphi _{n,x}\right\rangle ,
\end{equation*}
and with \eqref{eq_4_4}, \eqref{6eq_3_19}$_{1}$, \eqref{6eq_3_19}$_{3}$ and 
\eqref{6eq_3_44+}, we see that 
\begin{equation*}
\left( \frac{k_2}{k_1}-1\right) \left\langle \overset{\sim}{\psi}_{n},i\lambda_{n}\varphi _{n,x}\right\rangle -\left\Vert \overset{\sim}{\psi}_{n}\right\Vert^{2} \rightarrow 0,
\end{equation*}
so because $k_1 =k_2$, we obtain 
\begin{equation}
\overset{\sim}{\psi}_{n} \rightarrow 0\,\,\hbox{in}\,\,L_1 . \label{7d}
\end{equation}
Therefore, by \eqref{6eq_3_19}$_{3}$ and \eqref{7d}, we find 
\begin{equation}
\lambda_n \psi_{n} \rightarrow 0\,\,\hbox{in}\,\,L_1 . \label{7e}
\end{equation}
Taking the inner product in $L^2 (0,1)$ of \eqref{6eq_3_19}$_4$ with
$\psi _{n}$, using \eqref{eq_4_4}, integrating by parts and using the boundary conditions, we remark that  
\begin{equation}
\left\langle i\overset{\sim}{\psi}_{n},\lambda _{n}\psi _{n}\right\rangle +k_2\left\Vert \psi _{n,x}\right\Vert ^{2}+k_1\left\langle \varphi _{n,x}+\psi _{n} ,\psi_{n}\right\rangle \rightarrow 0.  \label{6eq_3_52}
\end{equation}
By using \eqref{eq_4_4}, \eqref{6eq_3_21}, \eqref{7e} and \eqref{6eq_3_52}, we arrive at \eqref{eq_3_53}. Consequently, \eqref{eq_3_20}, \eqref{eq_3_24}, \eqref{eq_3_53}, \eqref{6eq_3_21}, \eqref{7aa}, \eqref{7b} and \eqref{7d} lead to \eqref{limPhin}, which is a contradiction with \eqref{eq_4_4}.  
\vskip0,1truecm
\subsection{Polynomial and strong stability} The proof of the polynomial stability estimate \eqref{eq_4_1} is identical to the one given in section 5 by proving \eqref{polyc}. The first condition in \eqref{polyc} was proved in subsection 6.3 (which aready implies the strong stability estimate \eqref{SS51}; see Remark \ref{remark51}). The second condition in \eqref{polyc} can be proved using the same contradiction arguments used in section 5, where, in case of \eqref{syst16}-\eqref{cdt_106}, we have \eqref{wntildewnl}.


\begin{thebibliography}{99}
\bibitem{afas1} M. Afilal, A. Guesmia and A. Soufyane, New stability results for a linear thermoelastic Bresse system with second sound, Appl. Math. Optim., DOI: 10.1007/s00245-019-09560-7.

\bibitem{ag1} A. Afilal, A. Guesmia, A. Soufyane and M. Zahir, On the exponential and polynomial stability for a linear Bresse system, Math. Meth. Appl. Scie., 43 (2020), 2615-2625.

\bibitem{arba} W. Arendt and C. J. K. Batty, Tauberian theorems and stability of one one-parameter semigroups, Trans. Amer. Math. Soc., 306 (1988), 837-852. 

\bibitem{bres} J. A. C. Bresse, Cours de M\'{e}canique Appliqu\'{e}e, Mallet
Bachelier, Paris, 1859.

\bibitem{1} P. S. Casas and R. Quintanilla, Exponential decay in 
one-dimensional porousthermo-elasticity, Mech. Resea. Commu., 32 (2005), 
652-658.

\bibitem{chan} D. Chandrasekharaiah, Hyberpolic thermoelasticity: a review of recent literature, Applied Mechanics Reviews, 51 (1998), 705-729.

\bibitem{dafe} C. M. Dafermos, Asymptotic stability in viscoelasticity, Arch. Rational Mech. Anal., 37 (1970), 297-308.

\bibitem{2} F. Dell Oro and V. Pata, On the stability of Timoshenko systems with Gurtin-Pipkin thermal law, J. Diff. Equa., 257 (2014), 523-548.

\bibitem{ag4}  M. De Lima Santos, A. Soufyane and D. Da Silva Almeida 
J\'{u}nior, Asymptotic behavior to Bresse system with past history, Quart.
Appl. Math., 73 (2015), 23-54.

\bibitem{3} L. H. Fatori and J. E. Munoz Rivera, Energy decay for hyperbolic thermoelastic systems of memory type, Quart. Appl. Math., 59 (2001), 441-458.

\bibitem{fato1} L. H. Fatori and J. E. Munoz Rivera, Rates of decay to weak thermoelastic Bresse system, IMA J. Appl. Math., 75 (2010), 881-904.

\bibitem{4} H. D. Fernandez Sare and R. Racke, On the stability of damped Timoshenko system Cattaneo versus Fourier law, Arch. Rat. Mech. Anal., 194 (2009), 221-251.

\bibitem{gree1} A. Green and P. Naghdi, A re-examination of the basic postulates of thermomechanics, Proceedings of the Royal Society of London, Series A: Mathematical and Physical Sciences, 432 (1991), 171-194.

\bibitem{gree2} A. Green and P. Naghdi, On undamped heat waves in an elastic solid, Journal of Thermal Stresses, 15 (1992), 253-264.

\bibitem{8} A. E. Green and P. M. Naghdi, Thermoelasticity without 
energy-dissipation, J. Elast., 31 (1993), 189.

\bibitem{ag2} A. Guesmia, Non-exponential and polynomial stability results of a Bresse system with one infnite memory in the vertical displacement, Nonauton. Dyn. Syst., 4 (2017), 78-97.

\bibitem{ag3} A. Guesmia, The effect of the heat conduction of types I and III on the decay rate of the Bresse system via the longitudinal displacement, Arab. J. Math., 8 (2019), 15-41.

\bibitem{ag6} A. Guesmia, The effect of the heat conduction of types I and III on the decay rate of the Bresse system via the vertical displacement, Applicable Analysis, DOI: 10.1080/00036811.2020.1811974.

\bibitem{9} M. E. Gurtin and A. C. Pipkin, A general theory of heat conduction with finite wave speeds, Arch. Ration. Mech. Anal., 31 (1968),  113-126. 

\bibitem{huan} F. L. Huang, Characteristic condition for exponential stability of linear dynamical systems in Hilbert spaces, Ann. Diff. Equa., 1 (1985), 43-56.

\bibitem{kama} A. Keddi, T. Apalara and S. A. Messaoudi, Exponential and polynomial decay in a thermoelastic-Bresse system with second sound, Appl. Math. Optim., 77 (2018), 315-341.

\bibitem{lagn2} J. E. Lagnese, G. Leugering and J. P. Schmidt, Modelling of dynamic networks of thin thermoelastic beams, Math. Meth. Appl. Scie., 16 (1993), 327-358.

\bibitem{lagn1} J. E. Lagnese, G. Leugering and J. P. Schmidt, Modelling Analysis and Control of Dynamic Elastic Multi-Link Structures, Systems Control Found. Appl., 1994.

\bibitem{liu0} Z. Liu and B. Rao, Characterization of polymomial decay rate for the solution of linear evolution equation, Z. Angew. Math. Phys., 56 (2005), 630-644.

\bibitem{liu2} Z. Liu and B. Rao, Energy decay rate of the thermoelastic Bresse system, Z. Angew. Math. Phys., 60 (2009), 54-69.

\bibitem{liu1} Z. Liu and S. Zheng, Semigroups associated with dissipative systems, 398 Research Notes in Mathematics, Chapman \& Hall CRC, 1999.

\bibitem{10} H. W. Lord and Y. Shulman, A generalized dynamical theory of thermoelasticity, J. Mech. Phys. Sol., 15 (1967), 299-309.

\bibitem{11} A. Magana and R. Quintanilla, On the time decay of solutions in one-dimensional theories of porous materials, Inter. J. Sol. Stru., 43 (2006), 3414-3427.

\bibitem{13} S. A. Messaoudi and A. Fareh, Energy decay in a Timoshenko-type system of thermoelasticity of type III with different wave-propagation speeds, Arab J. Math., 2 (2013), 199-207.

\bibitem{12} S. A. Messaoudi, M. Pokojovy and B. Said-Houari, Nonlinear Damped Timoshenko systems with second: Global existence and exponential stability, Math. Method. Appl. Sci., 32 (2009), 505-534.

\bibitem{14} S. A. Messaoudi and B. Said-Houari, Energy decay in a Timoshenko-type system of thermoelasticity of type III, J. Math. Anal. Appl., 348 (2008), 298-307.

\bibitem{16} J. E. Munoz Rivera and R. Quintanilla, On the time polynomial decay in elastic solids with voids, J. Math. Anal. Appl., 338 (2008), 
1296-1309

\bibitem{15} J. E. Munoz Rivera and R. Racke, Mildly dissipative nonlinear Timoshenko systems - Global existence and exponential stability, J. Math. Anal. Appl., 276 (2002), 248-278.

\bibitem{nawe} N. Najdi and A. Wehbe, Weakly locally thermal stabilization of
Bresse systems, Elec. J. Diff. Equa., 2014 (2014), 1-19.

\bibitem{17} V. Pata and E. Vuk, On the exponential stability of linear thermoelasticity, Contin. Mech. Thermodyn., 12 (2000), 121-130.

\bibitem{pazy} A. Pazy, Semigroups of linear operators and applications to partial differential equations, Springer-Verlag, New York, 1983.

\bibitem{prus} J. Pruss, On the spectrum of $C_0$ semigroups, Trans. Amer. Math. Soc., 284 (1984), 847-857.

\bibitem{19} M. L. Santos, D. S. Almeida Junior and J. E. Munoz Rivera, The stability number of the Timoshenko system with second sound, J. Diff. Equa.,  253 (2012), 2715-2733.

\bibitem{timo} S. Timoshenko, On the correction for shear of the differential equation for transverse vibrations of prismatic bars,
Philosophical Magazine, 41 (1921), 744-746.
\end{thebibliography}
\end{document}